\documentclass[a4paper,11pt]{article}

\usepackage{geometry}
\usepackage{enumitem}
\usepackage{amsmath,amscd}
\usepackage{amsfonts}
\usepackage{amsthm}
\usepackage{amssymb}
\usepackage{upref}
\usepackage{color}
\usepackage{graphicx}
\usepackage{hyperref}
\usepackage{array}
\usepackage{algorithm}
\usepackage{algorithmicx}
\usepackage{algpseudocode}
\usepackage[utf8]{inputenc}
\usepackage{palatino}
\usepackage{multirow}
\graphicspath{{./pics/}}

\theoremstyle{definition}

\setlist[itemize]{noitemsep, topsep=0pt}

\begin{document}

\title{Efficient algorithm for optimizing spectral partitions}
\author{Beniamin Bogosel}
\date{}



\maketitle
\begin{abstract}
We present an amelioration of current known algorithms for optimal spectral partitioning problems. The idea is to use the advantage of a representation using density functions while decreasing the computational time. This is done by restricting the computation to neighbourhoods of regions where the associated densities are above a certain threshold. The algorithm extends and improves known methods in the plane and on surfaces in dimension $3$. It also makes possible to make some of the first computations of volumic 3D spectral partitions on sufficiently large discretizations.
\end{abstract}



\section{Introduction}

Recently many works were dedicated to the numerical study of spectral partitions in the plane and on surfaces in $\Bbb{R}^3$. These studies are motivated by the lack of an exact theoretical description even in the simplest cases. One question which remains open is the study of the spectral honeycomb problem. What is the partition of the plane which asymptotically minimizes the sum of the first eigenvalues of the Dirichlet Laplacian of each cell? It was conjectured in \cite{caflin} that the honeycomb partition is the best one. Big steps towards the proof of this conjecture were made in \cite{BFVV17}, where the authors prove that the spectral honeycomb conjecture is true among convex sets under the hypothesis that Polya's conejcture holds for hexagons and a slightly weaker version holds for pentagons and heptagons. We recall that Polya's conjecture states that among $n$-gons of fixed area, the regular one minimizes the first Dirichlet eigenvalue. Recent numerical simulations which confirm Polya's conjecture can be found in \cite[Chapter 1]{beniphd}.

The numerical aspects of the spectral partitioning problem were approached by Bourdin, Bucur and Oudet in \cite{buboou}. They proposed an algorithm which benefited of the fact that the cells could be represented as density functions on the same grid. This allows a straightforward implementation of a gradient based optimization method for the search of optimal numerical configurations. They managed to study partitions of up to $512$ cells and they noticed that locally partitions seem to be composed of patches of regular hexagons. Due to the size of the computations the $512$ cell computations was done at Texas Advanced Computing Center. Below we propose a modification of this algorithm which does not have any loss on the precision and which decreases in a significant way the cost of the computation. The case of $512$ cells can now be run in a reasonable amount of time on a personal computer. Some simulations for more than $1000$ cells on finer discretizations are possible and some of these are presented in the following.

We recall the following works dealing with numerical computations regarding spectral optimal partitions. In \cite{cybhol,CBH05} the authors present a problem issued from the modelization of chemical reactions for which the stationary state energy is the sum of the Dirichlet Laplace eigenvalues of the cells. They present evidence that the optimal configurations approach the hexagonal one in the for large number of components. In \cite{oow-spectral} Osting, Oudet and White investigate graph partitions using spectral methods. They propose a different optimization method which is gradient free and converges to a local minimum in a finite number of iterations. In \cite{elliott-ranner} the authors investigate partitions minimizing the sum of the Laplace-Beltrami eigenvalues on different surfaces. A different approach is presented in \cite{beni-fsol} in the case of the sphere. An adaptation of the algorithm in \cite{buboou} is presented in \cite{bove-multiphase} in the case of the multiphase problem where the objective functional is the sum of the fundamental eigenvalues and an area penalization. In \cite{BZ15} a method for minimizing the sum of the eigenvalues and the maximal eigenvalue is proposed. In \cite{BoBN16} the authors study the minimization of the largest eigenvalue by minimizing some $p$-norms of eigenvalues for large $p$. They propose a grid restriction procedure in the plane with the purpose of obtaining better precision. This consisted in finding rectangular neighborhoods of the cells on which we restrict the computations. This article proposes further improvements of this procedure by considering even fewer points in the computational region. The grid restriction procedure reduces greatly the computational cost, which allows the study of partitions into a large number of cells with low computational resources. The reduction in computational complexity also allows the study of partitions on surfaces and even volumic partitions for domains in $\Bbb{R}^3$.

A similar approach was devised in \cite{BDLO17} for the study of large partitions which minimize the total perimeter. It is another example where the representation of the sets as density functions helps when dealing with partitions, which is made possible with the use of the Modica-Mortola approximations of the perimeter by $\Gamma$-convergence. For details see \cite{oudet}, \cite{BoOu16}.

This paper focuses on describing how the methods presented in \cite{buboou} can be modified so that we gain in precision and in computational speed. We underline the speed improvements which are obtained and we propose a number of simulations that can be made with this method. Notably, we present some of the first three dimensional simulations for spectral partitioning problems. The 3D problem has also been approached in \cite{OZ16} for a periodic cube, but on a rather small discretization of $26\times 26\times 26$ and up to at most three cells. The method described in this article allows us to work in 3D with over one hundred cells on grids of size $100\times 100 \times 100$.

\section{Numerical algorithm}
\label{sect.num}

Let $D \subset \Bbb{R}^d$ be a bounded open set and $\omega$ one of its subsets. The classical results of \cite{dalmaso-mosco} allow us to find the eigenvalues of the Dirichlet Laplace operator on $\omega$ as solutions of a problem in $D$. More precisely, for a measure $\mu$ we can consider the problem
\begin{equation*}
-\Delta u + \mu u = \lambda(\mu) u \text{ on }D
\end{equation*}
for $u \in H_0^1(D)$. In the case where $\mu$ is given by 
\[ \mu(X) = \begin{cases} 
 0 & \text{ if } \text{cap} (X \cap \omega) = 0 \\
 +\infty & \text{ otherwise}
\end{cases} \]
then $\lambda(\mu)$ is precisely the first eigenvalue of the Dirichlet Laplacian operator on $\omega$. Here $\text{cap}$ denotes the capacity. For more details see \cite{bucurbuttazzo} and \cite{henrot-pierre}. Now, instead of working with the set $\omega$ we can use a density approximation. If $\varphi$ is close to the characteristic function of $\omega$ then solving
\begin{equation} (-\Delta+C(1-\varphi))u = \lambda(C,\varphi) u
\label{penalized}
\end{equation}
for $u \in H_0^1(D)$ and $C\gg 1$ will give us approximations of the eigenvalues and eigenfunctions of $\omega$. In \cite{buboou} it is proved that when $\varphi = \chi_\omega$ and $C\to \infty$ then $\lambda_k(C,\varphi)$ converges to $\lambda_k(\omega)$ as $C \to \infty$. There is also a quantitative convergence result presented in \cite{bove-multiphase} which gives an error estimate in terms of $C$
\[ \frac{|\lambda_k(\Omega)-\lambda_k(\mu_C)|}{\lambda_k(\Omega)} \leq KC^{-1/(N+4)},\]
where $\mu_C = C(1-\chi_\omega)$. This formulation inspired the numerical method presented in \cite{buboou}. 

Consider $D$ a rectangular box in $\Bbb{R}^d$ endowed with a finite differences uniform grid. We consider an ordering for the points of this grid, which will allow us to represent the discretization of a function $u$ on this grid by a column vector  $\bar u$.  Let $L$ be the matrix associated to the Laplace operator on this grid, in the sense that $L\bar u$ computes the discrete laplacian using five point stencils on the difference grid. We consider the discrete version of \eqref{penalized} given by
\begin{equation} (L+C\text{diag}(1-\varphi))\bar u = \lambda \bar u
\label{discretev}
\end{equation}
 In this way we may give an approximation of the eigenvalues and eigenfunctions of a set using a finite difference grid on a larger set and an approximation of its characteristic function. This immediately shows why formulation \eqref{discretev} is so well adapted to study partitioning problems: we can perform computations on a fixed grid in order to find quantities related to some of its subsets. 

Together with this advantage comes a drawback: we have a {\bf fixed} computation grid regardless of the {\bf size} of the cell $\omega$, for which we wish to compute the eigenvalue. For example, if a cell occupies only $5\%$ of the grid we have the same computational cost as if the cell occupied the whole grid. In addition to this we do not gain anything by working on the whole grid, on the contrary. Precision may be lost when penalizing too many points, since we add $C(1-\varphi)$ on the diagonal on the matrix of the laplacian. If $C$ is large and $\varphi$ is close to zero for many points, then the resulting matrix is ill conditioned. We propose that prior to the eigenvalue computation to look if the cell is localized to just one part of the grid and then restrict the computational domain to a suitably sized neighborhood around the cell. In the sequel we call \emph{computational neighborhood} of a cell, the region to which we restrict the computations regarding a particular cell.

\begin{figure}
\centering
\includegraphics[width = 0.2\textwidth]{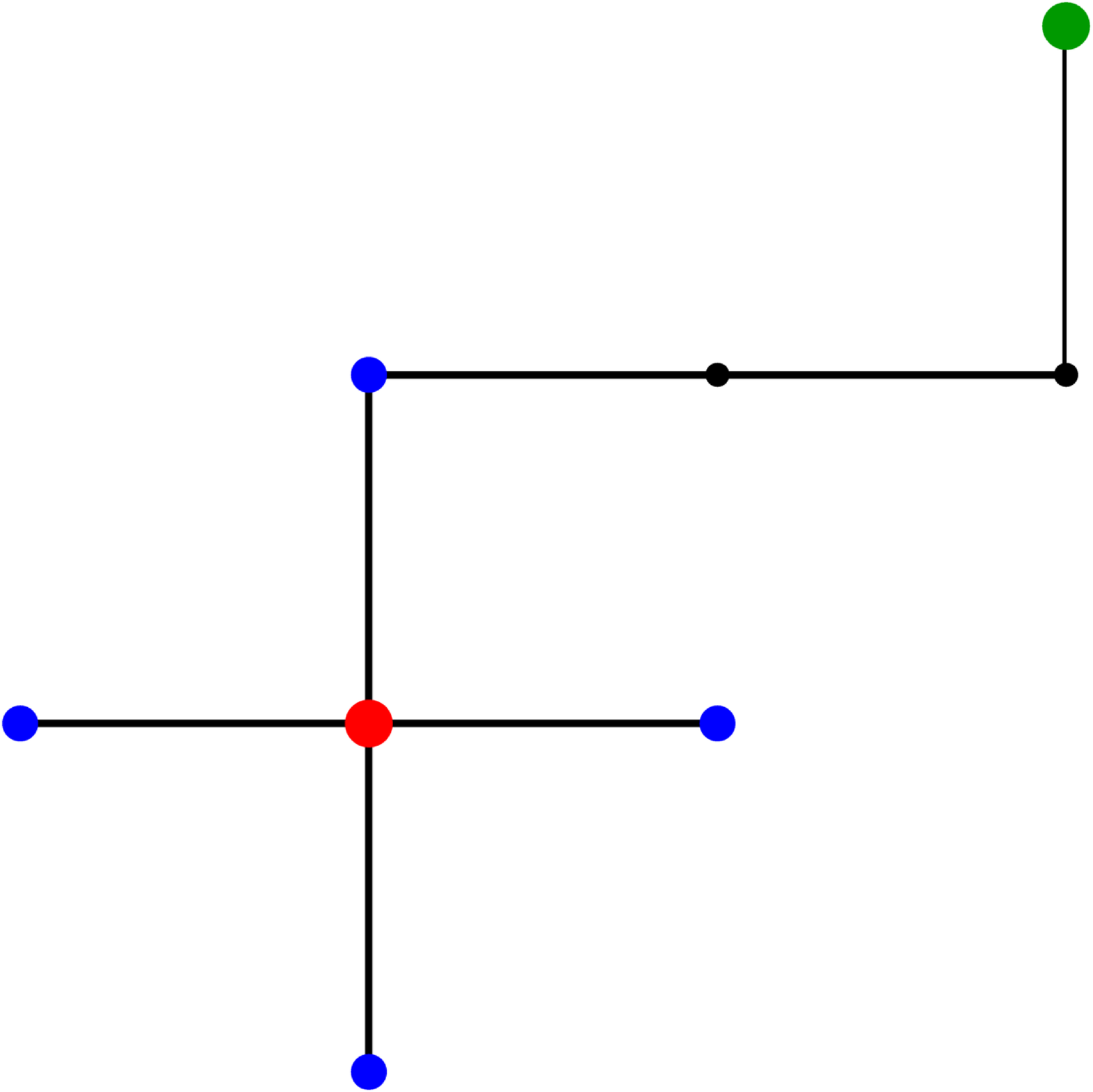}~
\includegraphics[width = 0.2\textwidth]{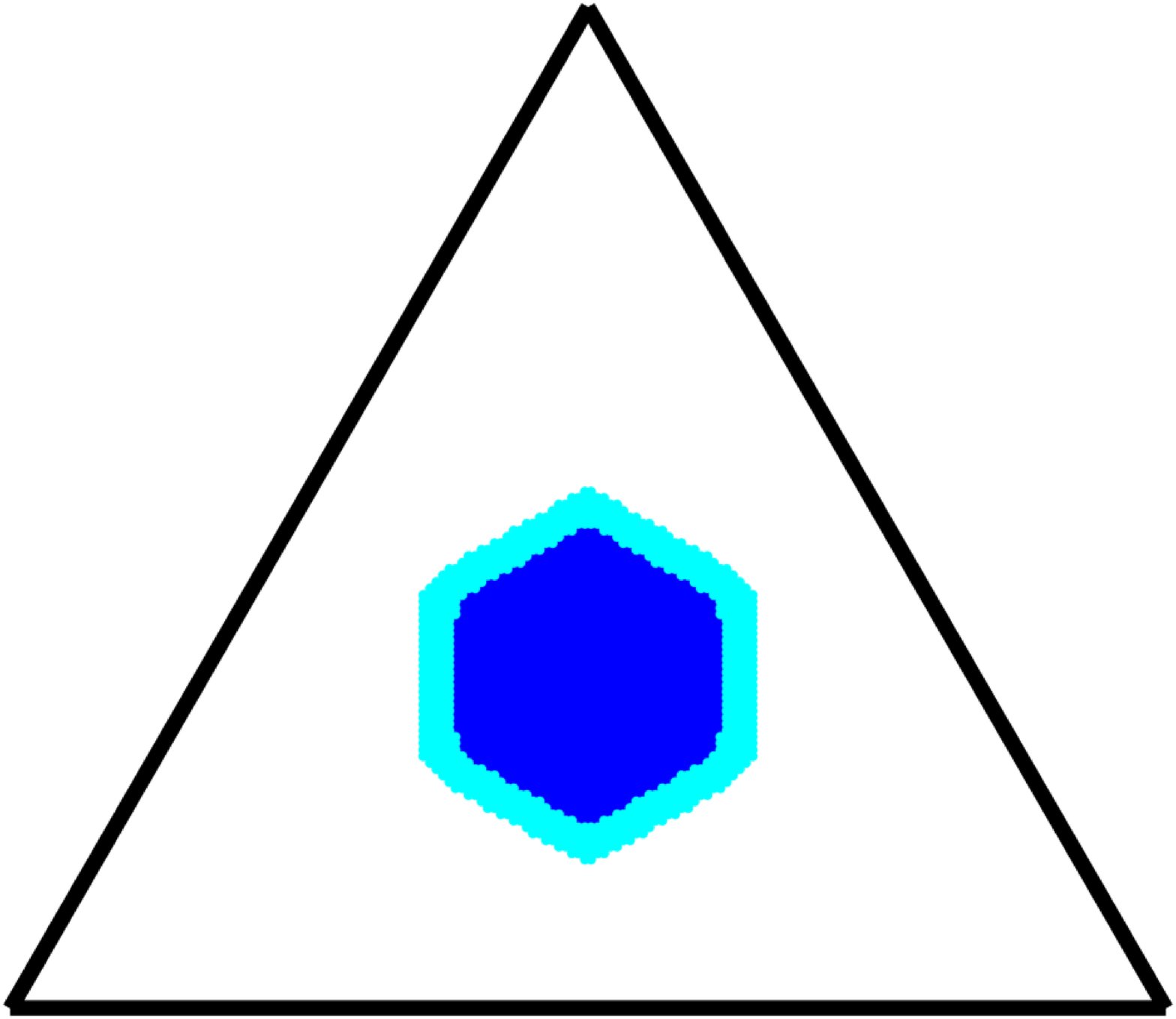}~
\includegraphics[width = 0.2\textwidth]{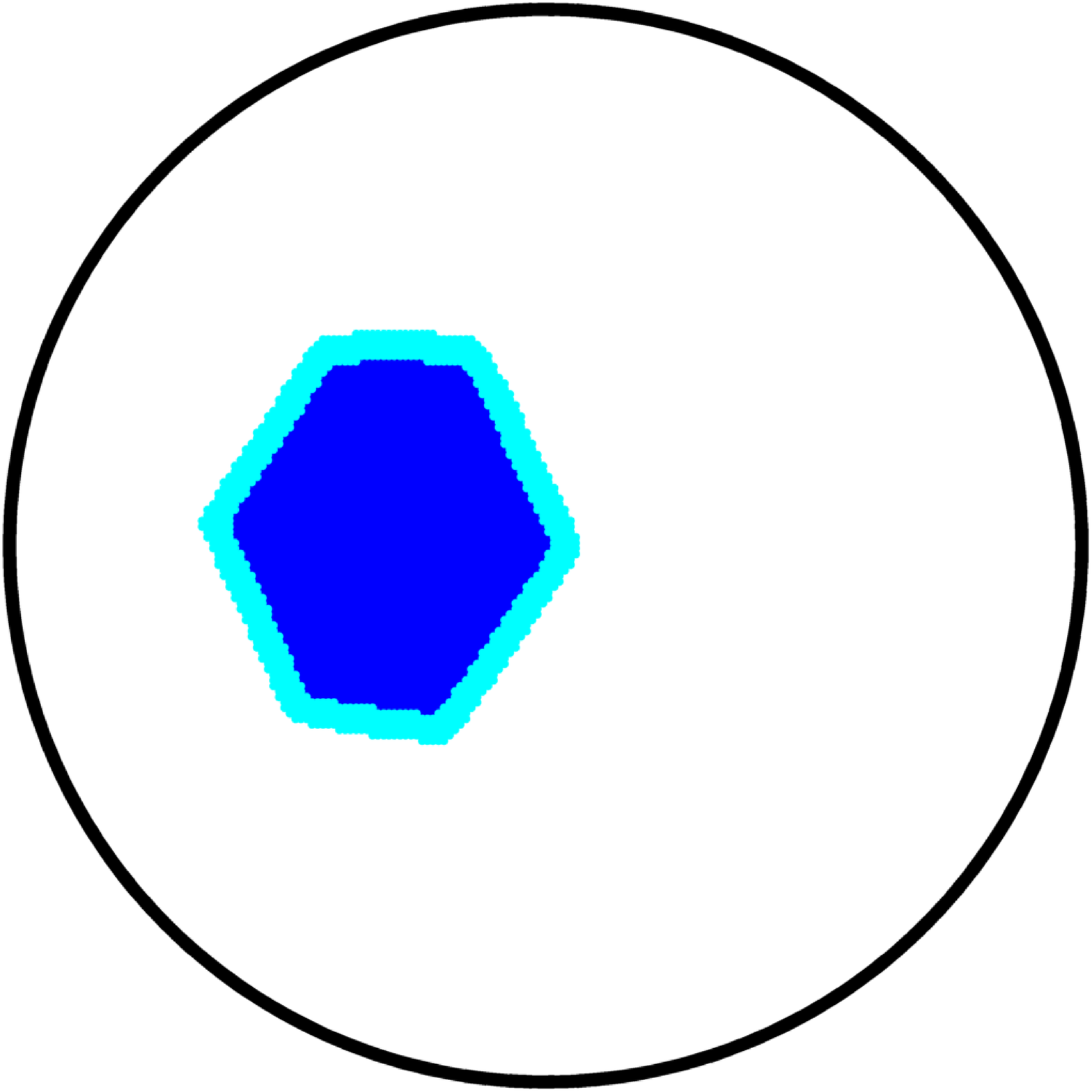}~
\includegraphics[width = 0.2\textwidth]{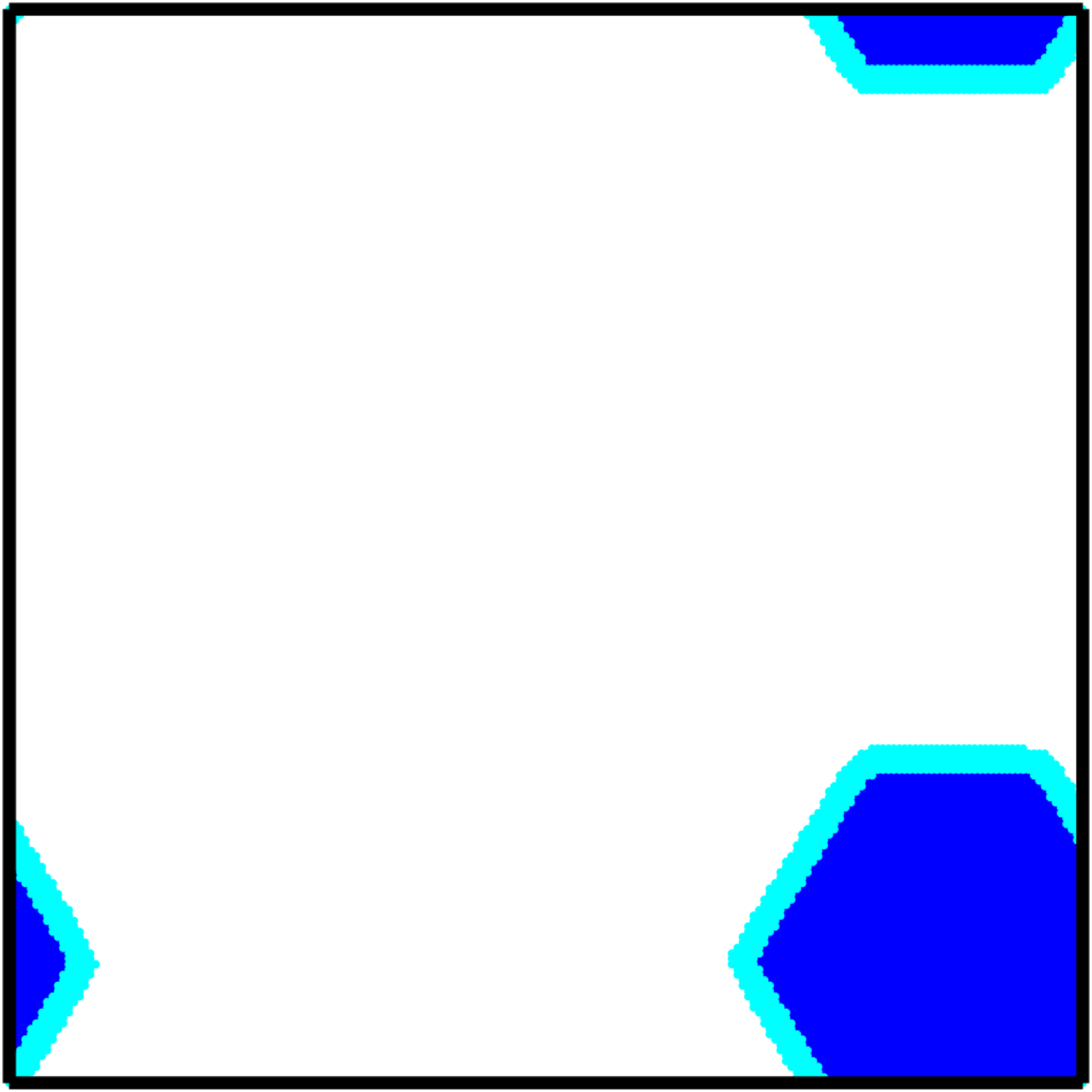}
\caption{Illustration of neighbors on the grid and example of reduced computational grids. Dark points represent points in the actual cell and light points represent additional points along the neighborhood of the cell considered in the computation.}
\label{grid_evol2D}
\end{figure}

First we locate all the points where the current density $\varphi$ is greater than some treshold. In the computations we took the points where $\varphi > 0.01$. In \cite{BoBN16} the authors considered a rectangular neighborhood of the cell with an eventual padding of $5$ to $10$ rows of discretization points to allow enough interaction between cells. A more careful analysis shows that the neighborhood does not need to be rectangular. Indeed, we see that having some exterior points in the computational neighborhood of the cell is necessary, since these points will be penalized in formulation \eqref{discretev} and we recover the Dirichlet boundary conditions for the eigenvalue problem. Moreover, it is enough to choose our computational neighborhood such that every point near the boundary of the cell has some exterior neighbors which are considered in the computation. We will call two points of the finite difference grid \emph{neighbors} if they are situated at a minimal vertical or horizontal distance on the grid. We call two points {\bf x}, {\bf y} on the grid \emph{neighbors of order $p$} if the $\ell^1$-distance between ${\bf x}$ and ${\bf y}$ is equal to $ph$, where $h$ is the discretization parameter. Equivalently, we can reach $\{\bf x\}$ starting from $\{\bf y\}$ by travelling $p$ vertical or horizontal grid segments joining adjacent vertices.  In the first picture of Figure \ref{grid_evol2D} the blue points are all the neighbours of the red point and the red point and green points are neighbors of order $4$. The neighbors of order at most $p$ can be computed starting from an initial adjacency matrix $N$, which contains $1$ on position $i,j$ if nodes $i$ and $j$ are neighbors. In general, we'll fix $p$ between $5$ and $10$, and for cell $i$ we define its associated computational neighborhood to be the family of all points on the grid which are neighbours of order at most $p$ with some of the points in $\{\varphi>0.01\}$. The procedure is presented in detail in Algorithm \ref{multi_neighbors}.

\begin{algorithm}
\caption{Find neighbors on the grid up to a certain order}
\label{multi_neighbors}
\begin{algorithmic}[1]
\Require 
\begin{itemize}
\item Finite difference matrix for the laplacian $L$.
\item $k$: depth of the neighbor search
\end{itemize}
\State construct neighbors adjacency matrix $N$ starting from $L$: $N=(a_{ij})$, $a_{ij} = 1$ if and only if $D_{ij}\neq 0$ and $i \neq j$. Else $a_{ij}=0$.
\State Initialize $N_f=N$. (at this stage we only have neighbors up to order $1$)
\For {$i=1:k$}
  \State $N_f = \min ( N_f+N\cdot A , 1)$
\EndFor

\Return $N_f$
\end{algorithmic}
\end{algorithm}

 We shall call this computational neighborhood $R$ in the sequel. We compute the eigenvalue of the cell corresponding to $\phi$ by using only the grid points inside $R$. In order to do this, we select from the matrix $L$ only the lines and columns corresponding to indices of points which are contained in the computational neighborhood $R$.  In the initial phase of the computation we start from random densities and the proposed algorithm doesn't bring much improvement. However, when the cells become localized the gain in speed is significant. Figure \ref{grid_evol2D} contains a few examples of such reduced computational neighbourhoods. 
If we want to consider the same kind of algorithm on a periodic grid it is enough to modify the adjacency matrix $N$ so that we include neighbors along the sides of the computational domain.

\begin{figure}[!ht]
\centering
\includegraphics[height=0.3\textwidth]{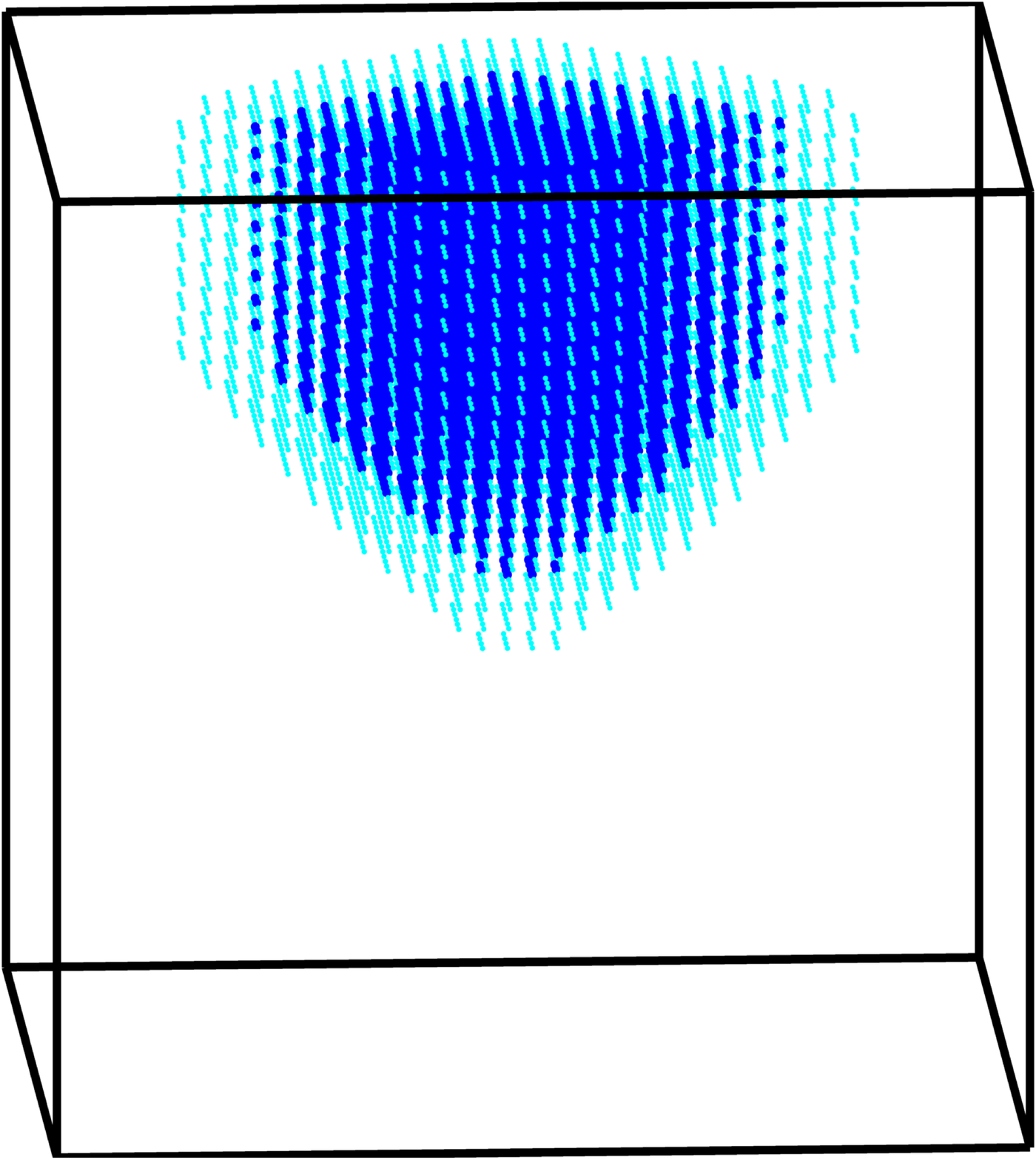}~
\includegraphics[height=0.3\textwidth]{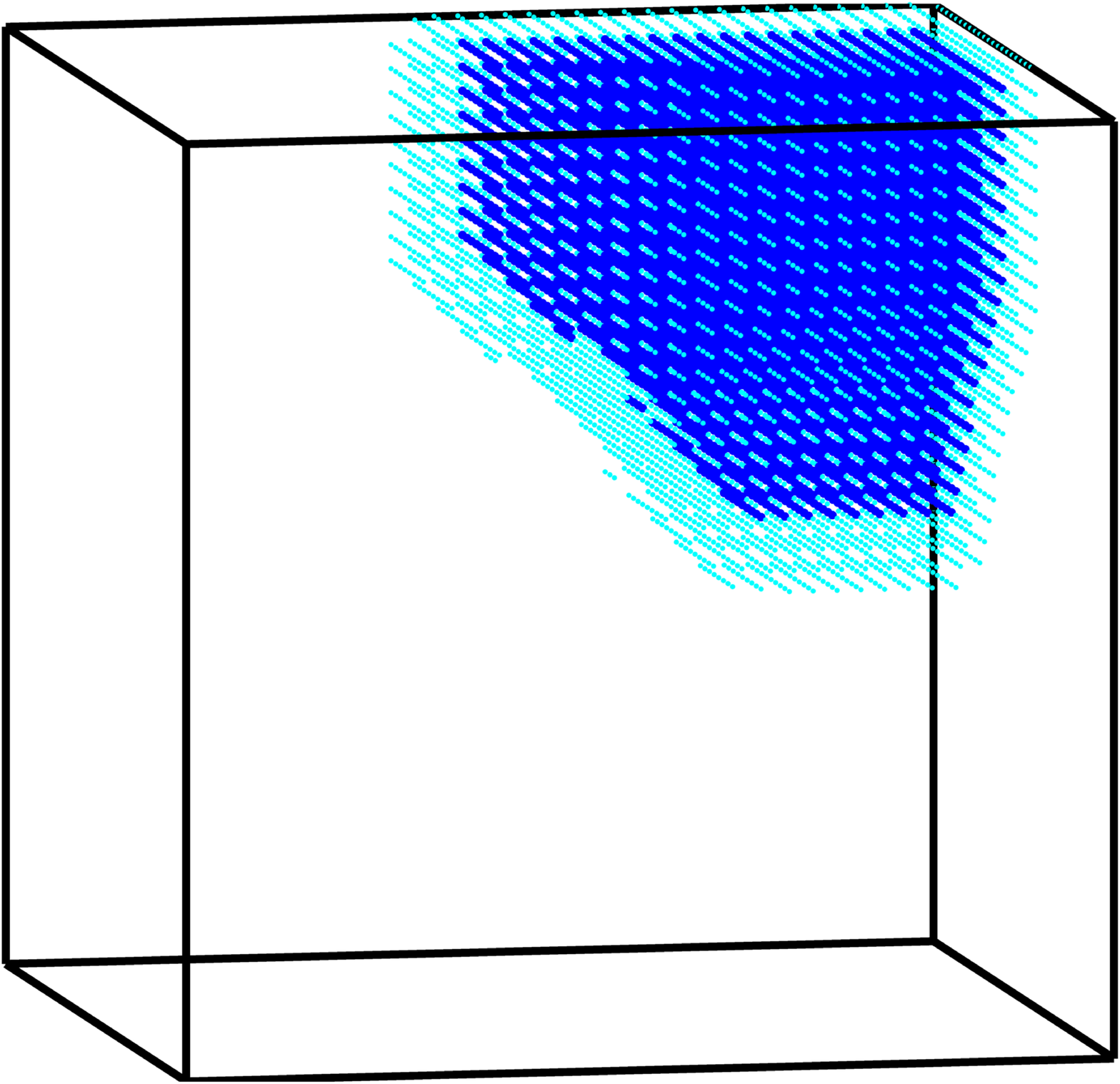}~
\includegraphics[height=0.3\textwidth]{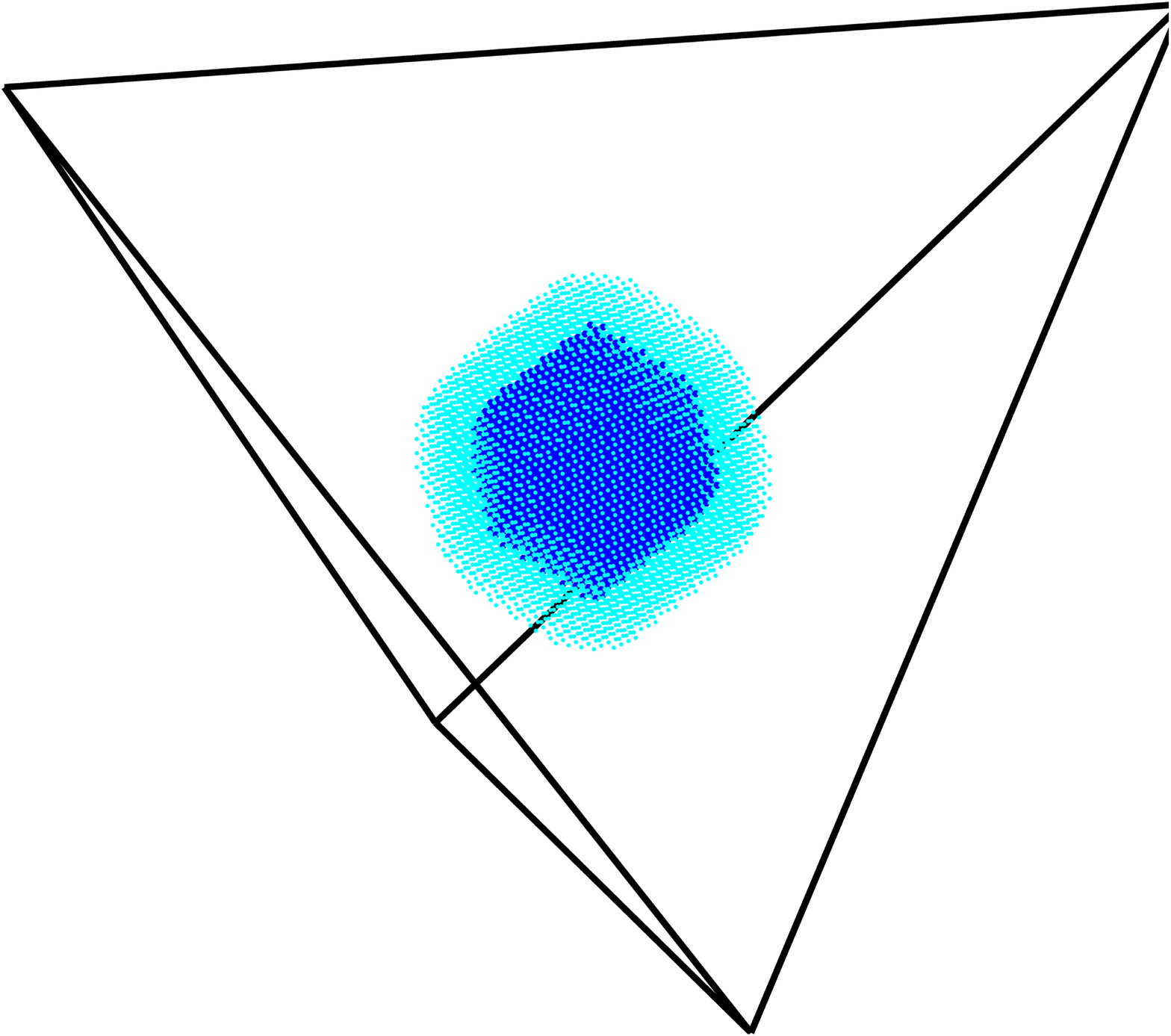}
\caption{Examples of reduced grids in 3D. Cells are defined as functions on the whole grid, but only the values corresponding to the points shown are used in the computation. Points in the cell are represented with dark blue. Points in the computational neighborhood are represented with cyan.}
\label{grid_3D}
\end{figure}

Note that until now no approach is known to work in the three dimensional case for high resolutions. If we want to use the approach in \cite{buboou} the problem quickly becomes too costly from a computational point of view since matrices involved are really large. Just to have an idea, working with finite differences without grid localization in 2D we can treat discretizations of up to $500\times 500$. This corresponds to a $60\times 60 \times 60$ discretization on the cube, with a denser matrix. Computing eigenvalues for such a grid in 3D is extremely slow. This was the main reason for which no conclusive 3D results were obtained in \cite{buboou}. On the other hand, when using a grid refinement procedure like the one described above computational costs are reduced and three dimensional computations become possible. We present a few such cases in Section \ref{simulations}. Some examples of reduced grids corresponding to three dimensional cells are presented in Figure \ref{grid_3D}.

Since we are also interested in working with partitions on surfaces we recall the finite element approach presented in \cite{bove-multiphase}. If $S$ is a surface in $\Bbb{R}^3$ with a corresponding triangulation then we may consider the mass matrix $M$ and the rigidity matrix $K$ associated to the $P1$ finite elements on this triangular grid. The equation satisfied by Laplace Beltrami eigenvalues corresponding to a density $\varphi: S \to [0,1]$ is
\[ \int_S \nabla u \nabla v + C\int_S (1-\varphi) uv = \lambda\int_S uv,\]
for every $v \in H^1(S)$. For simplicity we keep the same notation for a function defined on $S$ and its finite element approximation. In the finite element formulation this becomes
\[ v^T K u + C v^T \text{diag}(1-\varphi) M u = \lambda v^T Mu,\] for every vector $v$. Thus we are left with a generalized eigenvalue problem
\begin{equation} 
(K+C\text{diag}(1-\varphi) M) u = \lambda M u,
\label{beltrami}
\end{equation}
which can also be solved using the Matlab function \texttt{eigs}. Note that in \cite{elliott-ranner} the authors studied partitions on three surfaces (the sphere a torus and another parametric surface) for up to $32$ cells. Partitions of the sphere into cells which minimize the sum of the first Laplace-Beltrami eigenvalues were studied in \cite{beni-fsol}.  One of the drawbacks of the approach used in \cite{beni-fsol} is the fact that it was limited to the case of the sphere and that the number of discretization points was limited to about $5000$, since the generalized eigenvalue problem obtained was not sparse.

\begin{figure}
\centering
\includegraphics[width = 0.23\textwidth]{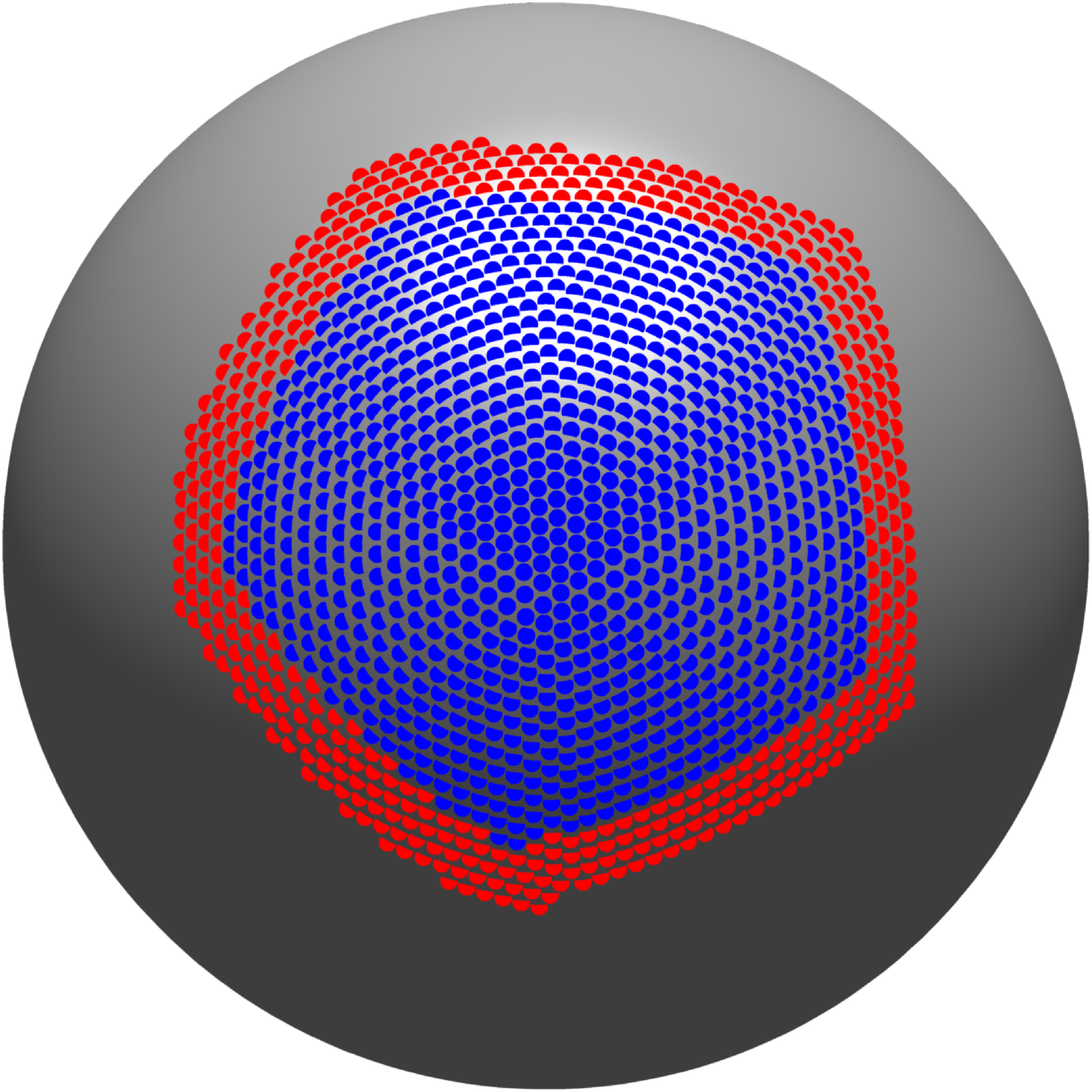}~
\includegraphics[width = 0.3\textwidth]{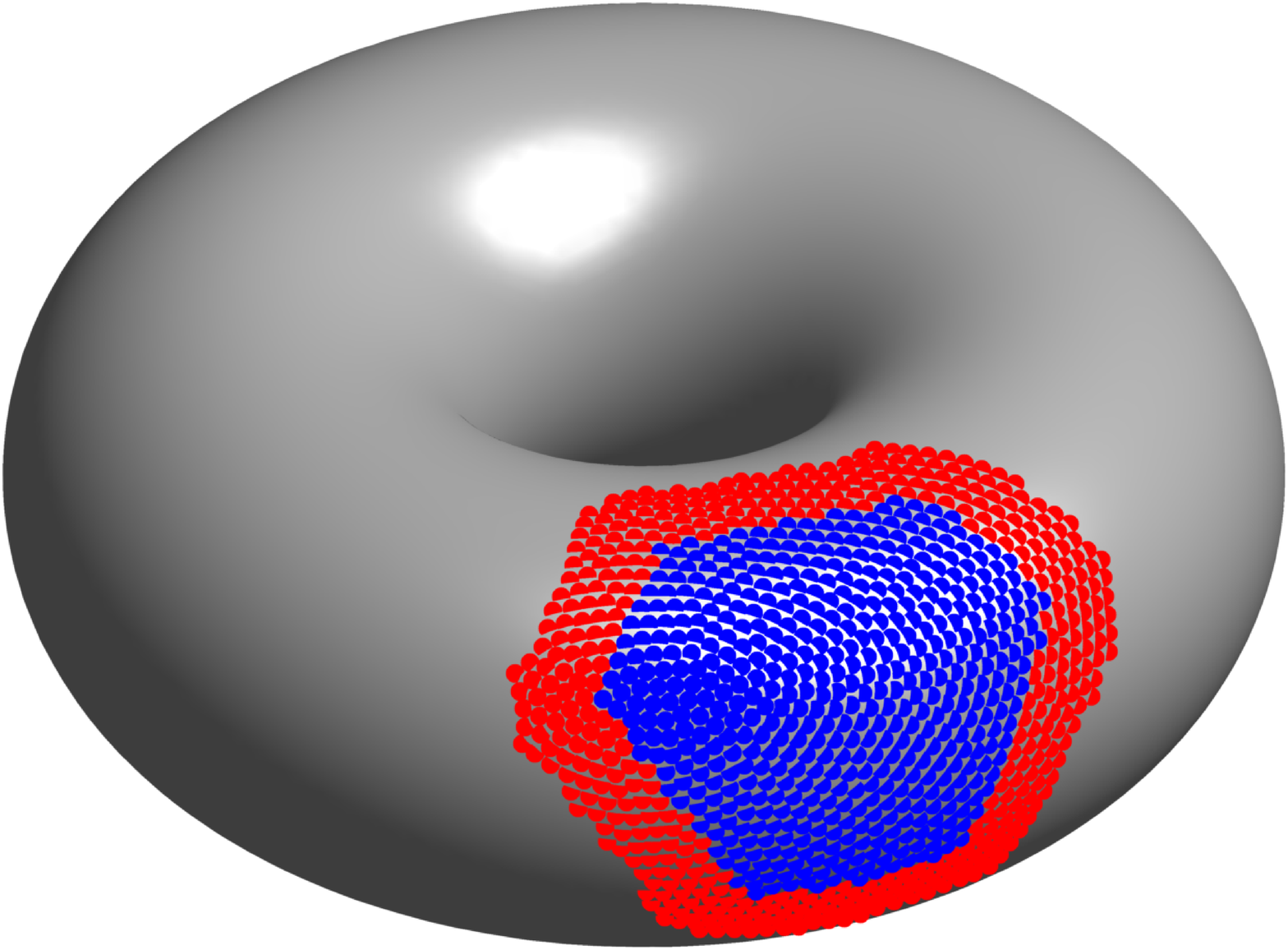}~
\includegraphics[width = 0.2\textwidth]{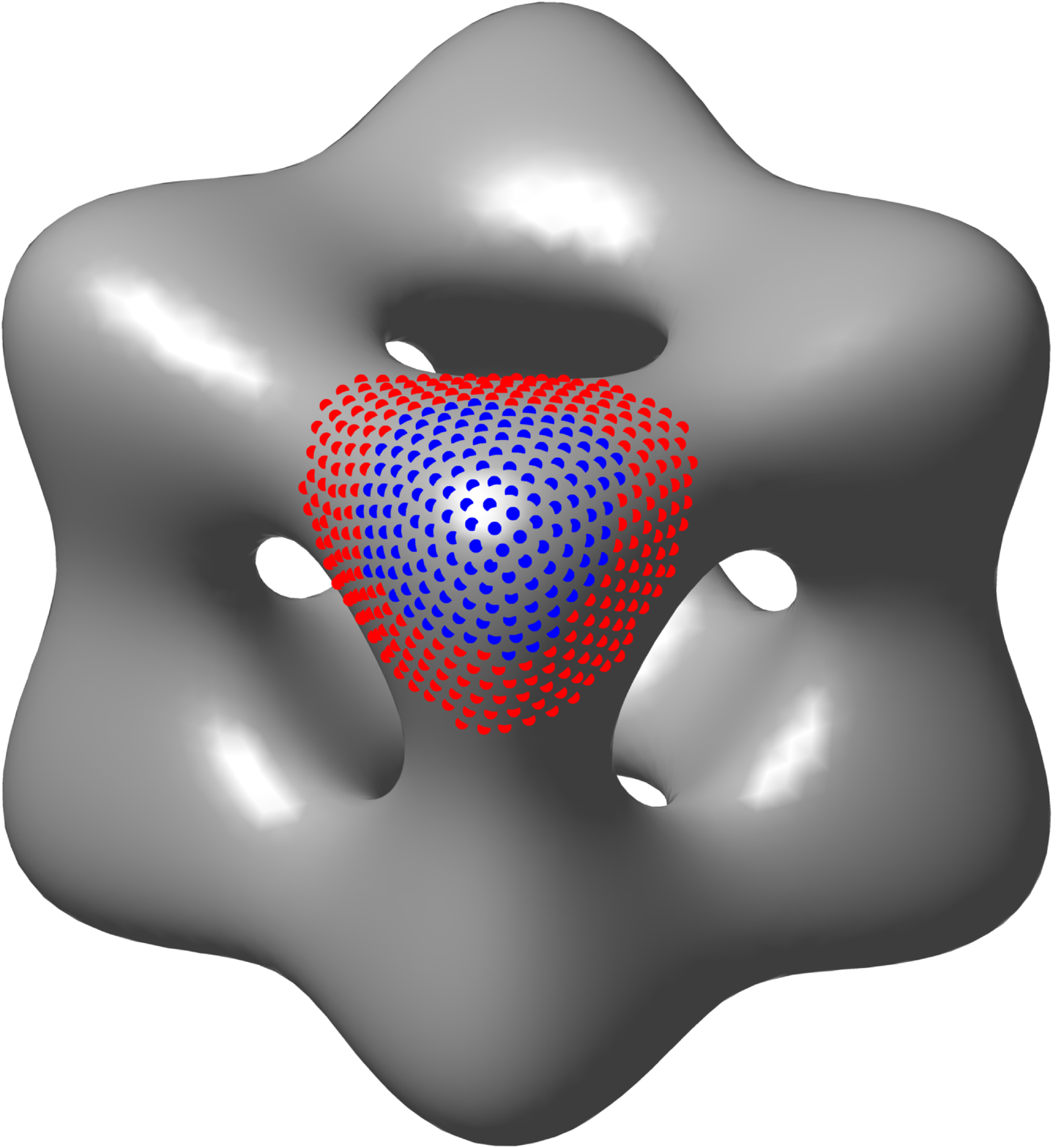}
\caption{Defining the computation grid on various surfaces. The blue points represent the and the red points are additional points taken in the computation neighborhood.}
\label{surface_grids}
\end{figure}

 In this setting we propose a similar grid restriction procedure which will allow us to work with a much larger number of discretization points while decreasing the computation time. On the sphere we preform computations on discretizations containing over $160000$ points. See section \ref{complexity} for more details. 
If $\varphi$ is the current density then we look for the points where $\varphi$ is greater than $0.01$. Once we have identified these points we search for all their neighbours which are at a distance of at most $5$ edges on the triangulation. See Figure \ref{surface_grids} for some examples. Then we compute the eigenvalue corresponding to this cell by solving a problem of the type \eqref{beltrami} only on this restricted part of the surface. 
  
In each of the situations presented above, partitions of the domain $D$ will be represented by $n$-tuples of functions  $(\varphi_1,...,\varphi_n)$ defined on $D$ which verify the partition condition
\[ \varphi_1+...+\varphi_n = 1\]
point-wise at each grid point. The initial condition will always be generated randomly and then projected on the constraint. In \cite{buboou} the authors proposed the following projection:
\[ \varphi_i \mapsto \frac{ |\varphi_i|}{\sum_{i=1}^n |\varphi_i|}.\]
Apart from the fact that this projects a family of functions on constant sum equal to 1 constraint, this also makes every function $\varphi_i$ have values in $[0,1]$. This is important in light of the fact that the application
\[ \varphi \mapsto \lambda_1(C,\varphi)\]
is concave in $\varphi$. Thus when minimizing $\sum_{i=1}^n\lambda_1(\varphi_i,C)$ each $\varphi_i$ will become an extremal point in the family of functions with values in $[0,1]$, and these extremal points are exactly the characteristic functions, taking values in $\{0,1\}$.

At each iteration of the optimization algorithm we have the following:
\renewcommand{\labelitemi}{$\cdot$}
\setlist{nolistsep}
\setlist{noitemsep}
\begin{itemize}[leftmargin=0.3cm]
\item compute the computational neighborhood of $\varphi_i$
\item compute the associated eigenvalue and eigenfunction by solving \eqref{discretev}, \eqref{beltrami} on the restricted grid
\item using the eigenfunction, compute the associated gradient of $\lambda_1$ with respect to each point of the grid. According to \cite{buboou} this is given by
\[ \partial_i \lambda_1(C,\varphi) = -C u_i^2 \]
where $u$ is the solution of \eqref{discretev}, \eqref{beltrami} associated to $\lambda_1(C,\varphi)$.
\item evolve each function in the direction of the gradient
\[ \varphi_i = \varphi_i+\alpha C u_i^2\]
and test if the new partition has a lower energy. If yes, then we continue, if not then we decrease the value of the step size $\alpha$.
\item since the evolution in the direction of the gradients will break the partition condition we project back on the constraint using the algorithm in \cite{buboou}, described above.
\item continue repeating these steps until we reach a preset maximal number of iterations, or the step $\alpha$ is smaller than $10^{-6}$.
\end{itemize}

The considerations above allow us not only to work on rectangular domains, but also in general domains in the following way. Suppose $D$ is a generic domain. Then we consider a rectangular domain $D'\supset D$ which contains it. On $D'$ we construct the finite difference framework with the difference that any point in $D'\setminus D$ will always be penalized in \eqref{discretev}, \eqref{beltrami}. This is equivalent to imposing $\varphi=0$ outside $D$. 

All figures were created using Matlab standard commands: \texttt{plot}, \texttt{plot3}, \texttt{surf}, \texttt{patch}.

\subsection{Automatic classification of cells} 

Since we are interested in studying partitions with many cells, we would like a way to automatically classify these cells, which would be more rigorous than visual inspection. In 3D it is even more difficult to visualize clearly what happens in the interior of the computational domain, hence the necessity of the tools presented below. 

For two dimensional domains and surfaces we have an immediate way of classifying cells by counting neighbours. We will loosely use the terms pentagons, hexagons, heptagons to denote cells with five, six or seven neighbours. When looking at patches of cells with six neighbours these seem to be made out of polygonal hexagons, in the classical sense. However, when pentagons or heptagons are present, boundaries seem to be slightly curved in certain cases.

Counting the neighbours of each cell is not a difficult task. At the end of the optimization procedure described previously, we have each cell $\omega_i$ represented by a density function \mbox{$\varphi_i: D \to [0,1]$} evaluated at the nodes of the discretized domain. As in the construction of the computational cell, we make use of the matrix of neighbors $N$. We recall that $N$ is symmetric and line $\ell$ of $N$ contains ones on positions corresponding to nodes neighbors to node $\ell$ and zero on the other positions. Given two density functions $\varphi_i,\varphi_j$ we test if they are adjacent in the partition in the following way. Find indices $I_i$ and $I_j$ for positions where $\varphi_i,\varphi_j$ are greater than $0.5$. Then compute the submatrix of $N$ corresponding to lines $I_i$ and columns $I_j$. If this submatrix contains non-zero elements then we conclude that cells $i$ and $j$ are neighbors in the partition. This procedure is described in Algorithm \ref{detect_top2D}.

\begin{algorithm}
\caption{Counting neighbors for each cell}
\label{detect_top2D}
\begin{algorithmic}[1]
\Require 
\begin{itemize}
\item $n$: number of cells
\item $\varphi_i$, $i=1,...,n$ the density functions
 corresponding to the numerical optimizer
\item adjacency matrix $N$ of neighbors on the grid
\end{itemize}
\State For $1\leq i\leq n$ find indexes $I_i$ of the points in the grid where $\varphi_i>0.5$.
\State For each $1\leq i<j \leq n$ decide that cells $i$ and $j$ are neighbors if the submatrix of $N$ corresponding to lines $I_i$ and columns $I_j$ contains non-zero elements.
\end{algorithmic}
\end{algorithm}

In 3D, things are more complicated. Counting neighbors is not difficult, but it does not bring such significant information as in 2D, since polyhedral domains with the same number of faces may have quite different structure. In order to detect cells which are similar we use a heuristic method based on the spectra of the Laplace-Beltrami operators on closed surfaces. It is known that the spectrum of the Laplace operator does not completely characterize a shape. There are known isospectral domains which are not isometric. However, experimental works in \cite{reuter06} and \cite{laispectral} show that spectral analysis of the Laplace-Beltrami spectrum may help in the automatic classification of surfaces. Since this is not difficult to implement, in view of the fact that we study spectral problems in this work, we applied these ideas in order to classify the cells obtained in the 3D computations. The method manages to provide remarkably accurate classifications for the test cases we tried.

As in the 2D classification, we are given the output of the optimization algorithm, consisting of density functions $\varphi_i:D \to [0,1]$. Knowing the values of $\varphi_i$ on a structured grid, we may use the \texttt{isosurface} command in Matlab in order to construct a couple $(p,T)$ defining a triangulation of the level set $\{\varphi_i = 0.5\}$. Here $p$ contains the coordinates of the points and $T$ contains on each line indexes of points belonging to the same triangle. Given $(p,T)$ we may use $P1$ finite elements in order to approximate the first few Laplace-Beltrami eigenvalues of the surface. Once the mass matrix $M$ and the rigidity matrix $K$ are constructed, one simply needs to solve the generalized eigenvalue problem $Kv=\lambda Mv$. For this we use the command \texttt{eigs} in Matlab.
\begin{algorithm}
\caption{Automatic classification of surfaces using Laplace-Beltrami spectrum}
\label{detect_top3D}
\begin{algorithmic}[1]
\Require 
\begin{itemize}
\item $n$: number of cells
\item $\varphi_i$, $i=1,...,n$ the density functions
 corresponding to the numerical optimizer
 \item $\varepsilon$: treshold for spectrum comparison
 \item $k$ the number of Laplace-Beltrami eigenvalues computed for each surface
\end{itemize}
\State for each $\varphi_i$,\ $i=1,...,n$ compute the triangulation $(p_i,T_i)$ of the level set $\{\varphi_i=0.5\}$.
\State for each $i=1,...,n$ compute the mass matrix $M_i$ and rigidity matrix $K_i$ associated to the triangulation $(p_i,T_i)$.
\State for each $i=1,...,n$ compute the first $k$ smallest eigenvalues of the generalized eigenvalue problem $K_iu = \lambda M_iu$ and store them in a vector $v_i$. 
\State normalize each $v_i$: $v_i \gets v_i/\|v_i\|$. 
\State compute the matrix $\mathcal{M} = (\|v_i-v_j\|)_{1\leq i,j\leq n}$.
\State decide that two surfaces are similar if $\mathcal{M}_{ij} < \varepsilon$.
\end{algorithmic}
\end{algorithm} 

Once we have computed the first $k$ Laplace-Beltrami eigenvalues $\lambda_1,...,\lambda_k$ for each surface corresponding to $\{\varphi_i=0.5\}$ we normalize the vector $v_i = (\lambda_1,...,\lambda_k)$. We use normalized vectors since eigenvalues change proportionally when we dilate the domains. Next we want to see for which cells the normalized vectors of eigenvalues are correlated. In order to to this we simply compute the norm of the differences between $v_i$ and $v_j$. We assume that cells $i$ and $j$ are similar if $\|v_i-v_j\|<\varepsilon$, where $\varepsilon$ is a chosen threshold, which in our tests was $\varepsilon=0.01$. The procedure is described in Algorithm \ref{detect_top3D}.

\section{Numerical Simulations and Observations}
\label{simulations}

\subsection{Dimension two}
For planar domains we observe the emergence of hexagonal patches, as soon as the number of cells is large enough so that the boundary does not have much influence on the interior cells. Thus, it seems that numerical simulations, in accordance with previous simulations \cite{buboou}, \cite{CBH05}, confirm the Caffarelli-Lin conjecture.

\begin{figure}[!ht]
\centering
\includegraphics[width = 0.4\textwidth]{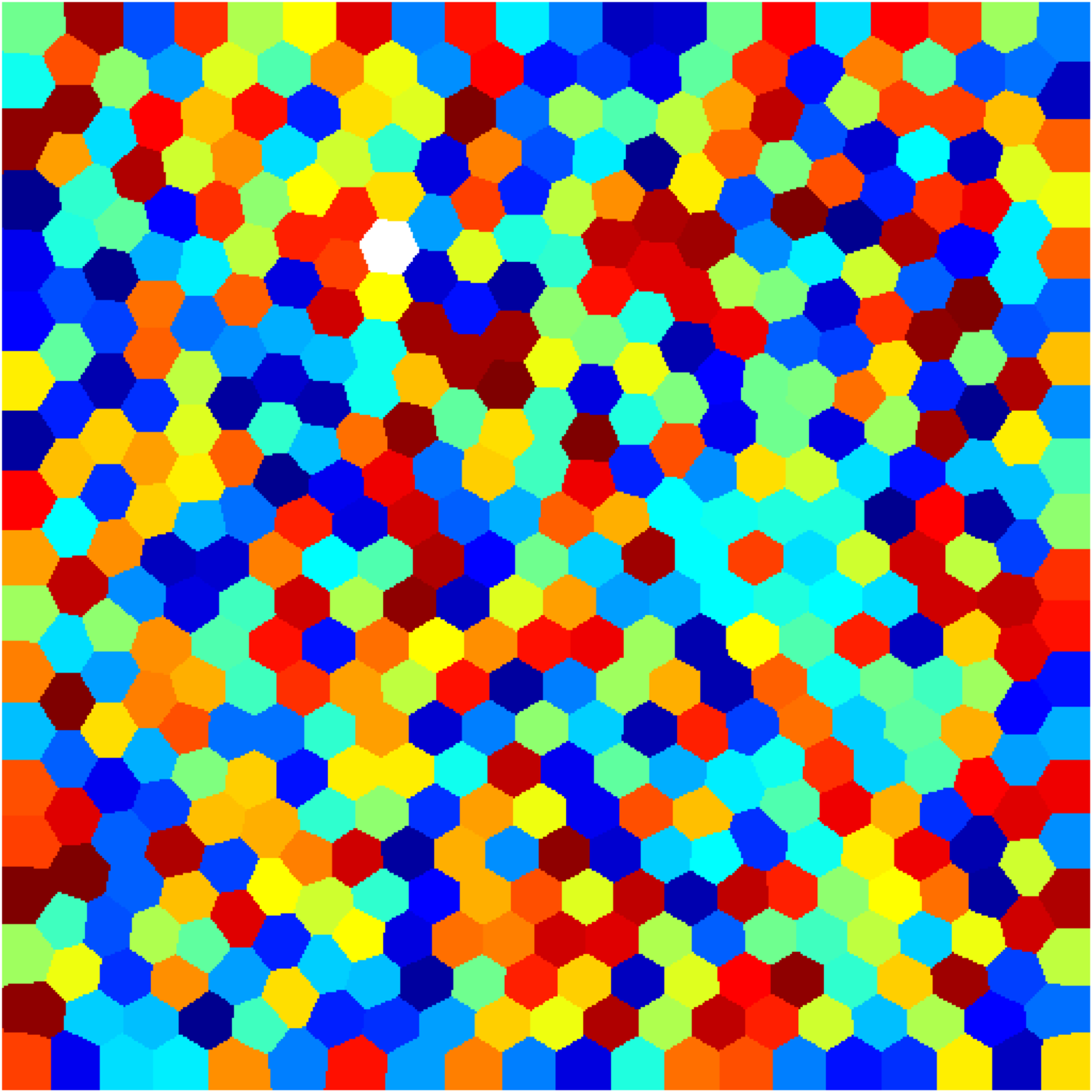}\quad
\includegraphics[width = 0.4\textwidth]{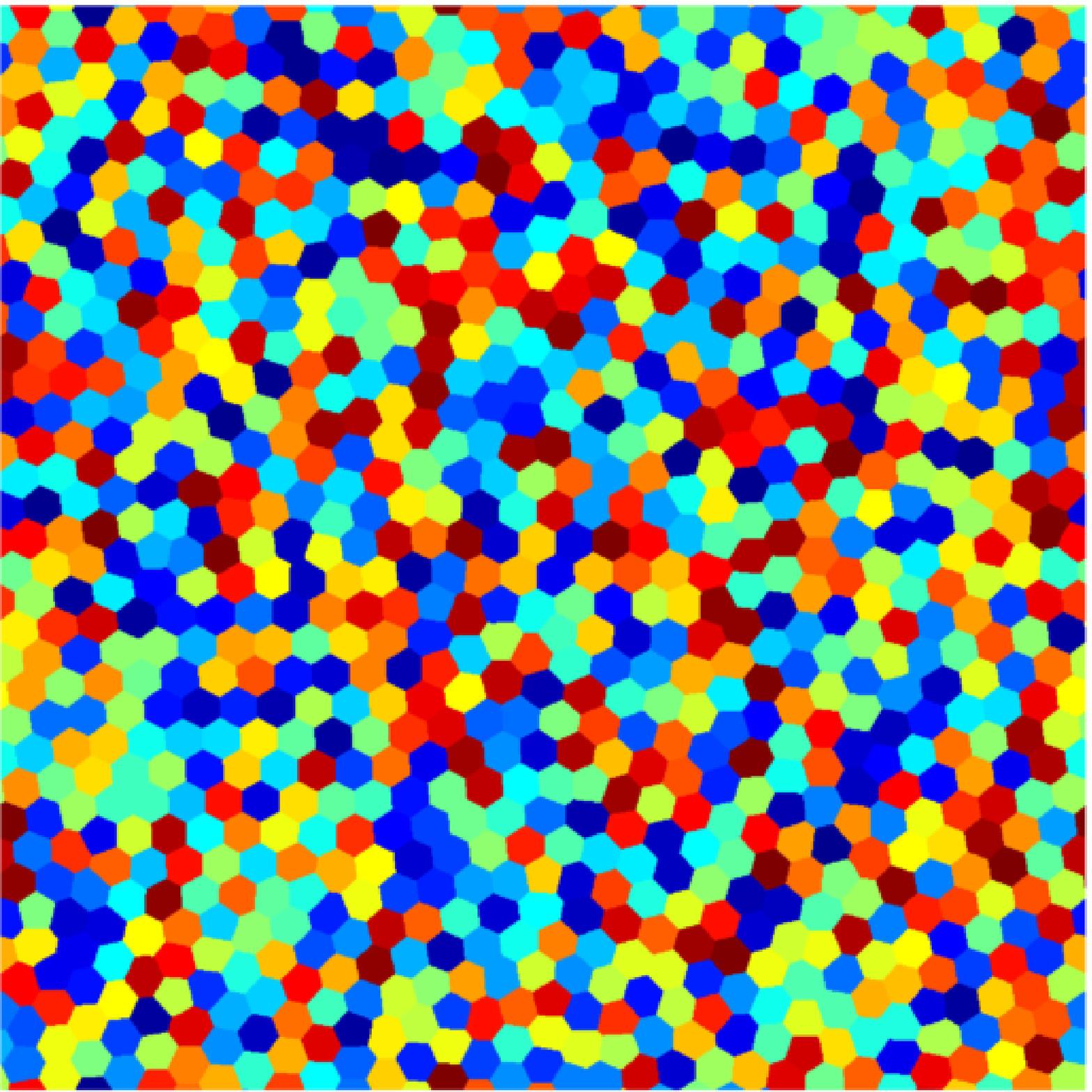}
\caption{$512$ cells on a grid of size $512\times 512$ and $1000$ cells on a grid of size $1000\times 1000$.}
\label{BBO}
\end{figure}

\begin{figure}
\centering
\includegraphics[width = 0.2 \textwidth]{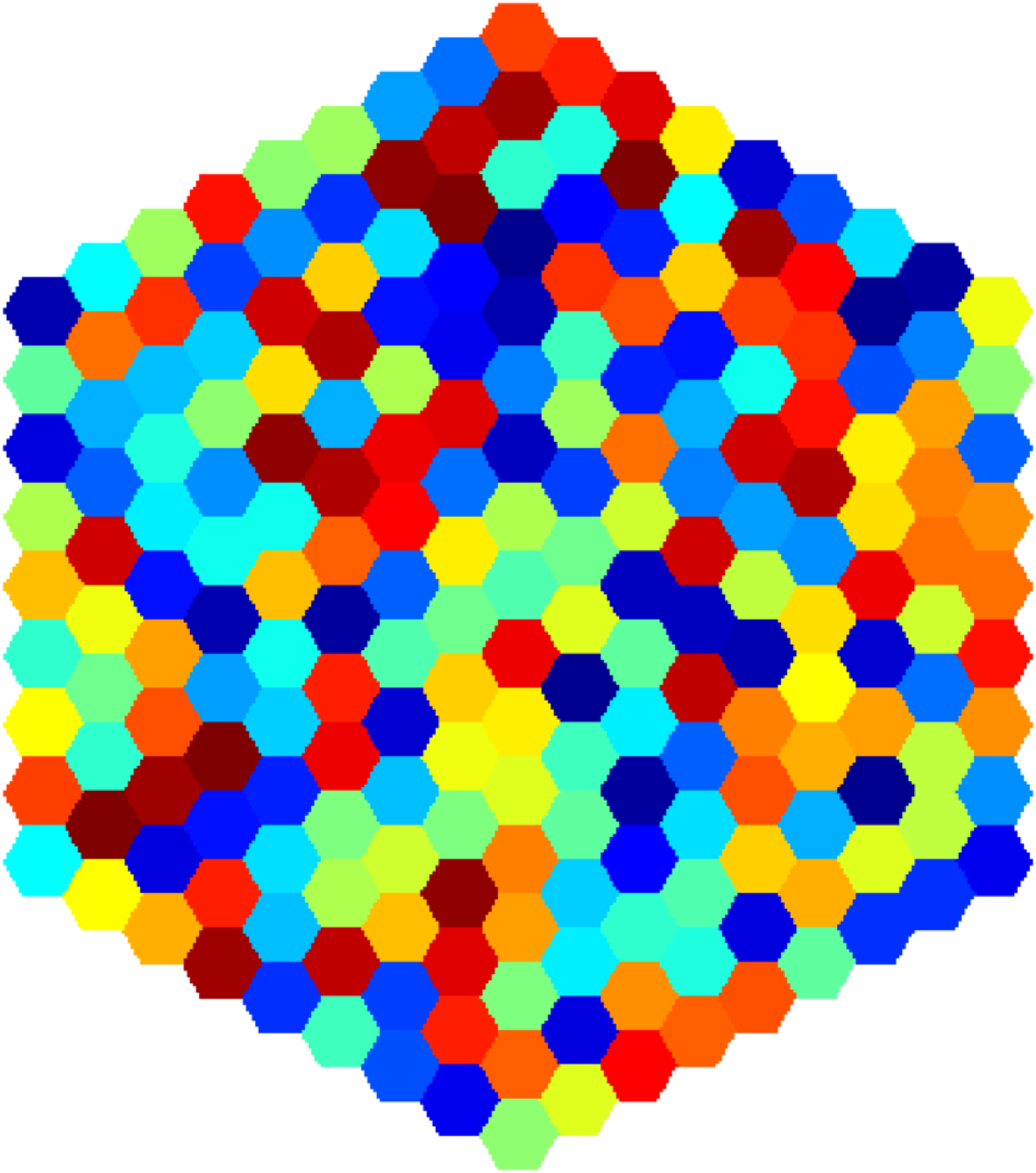}
\includegraphics[width = 0.2 \textwidth]{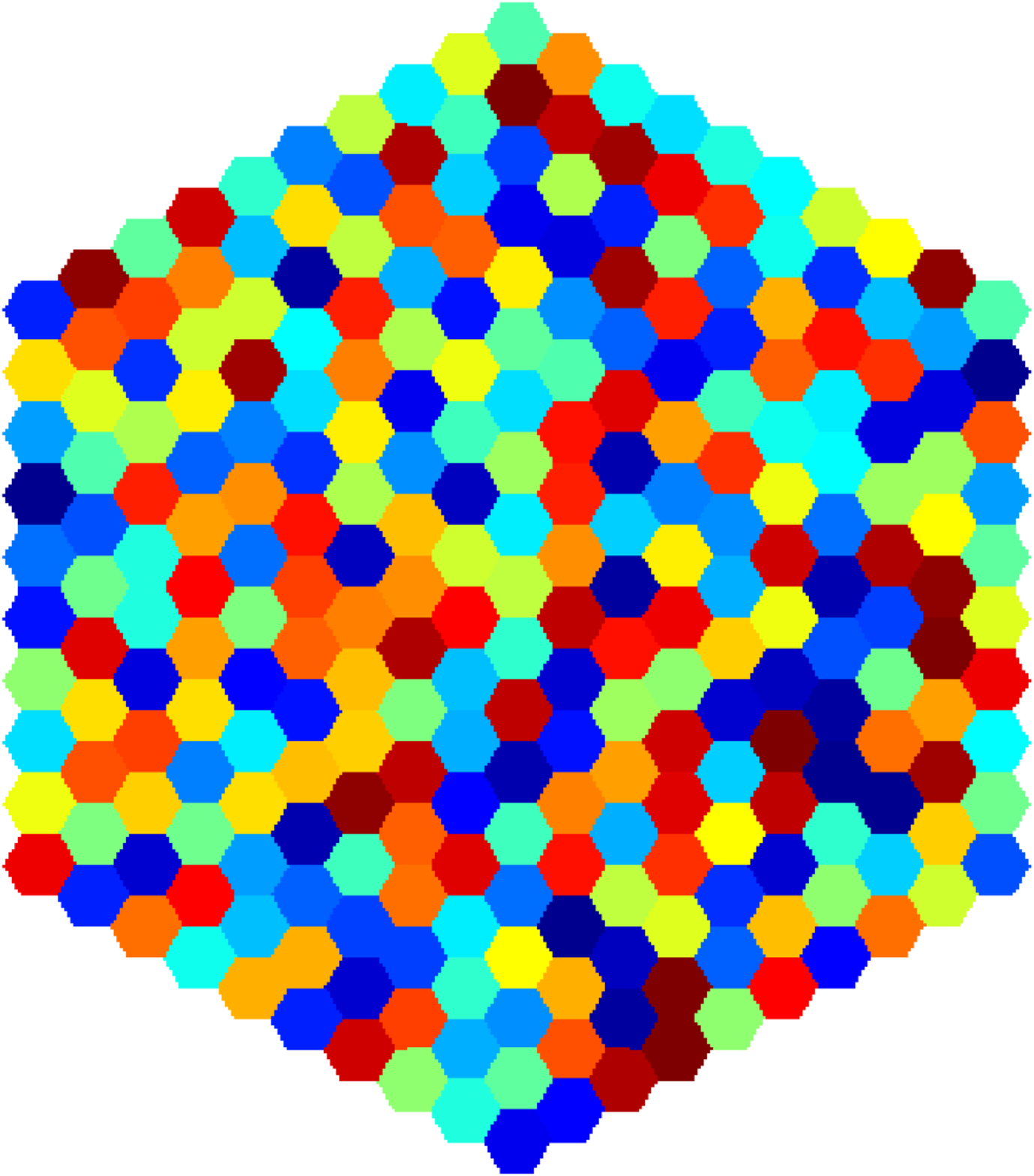}
\includegraphics[width = 0.2 \textwidth]{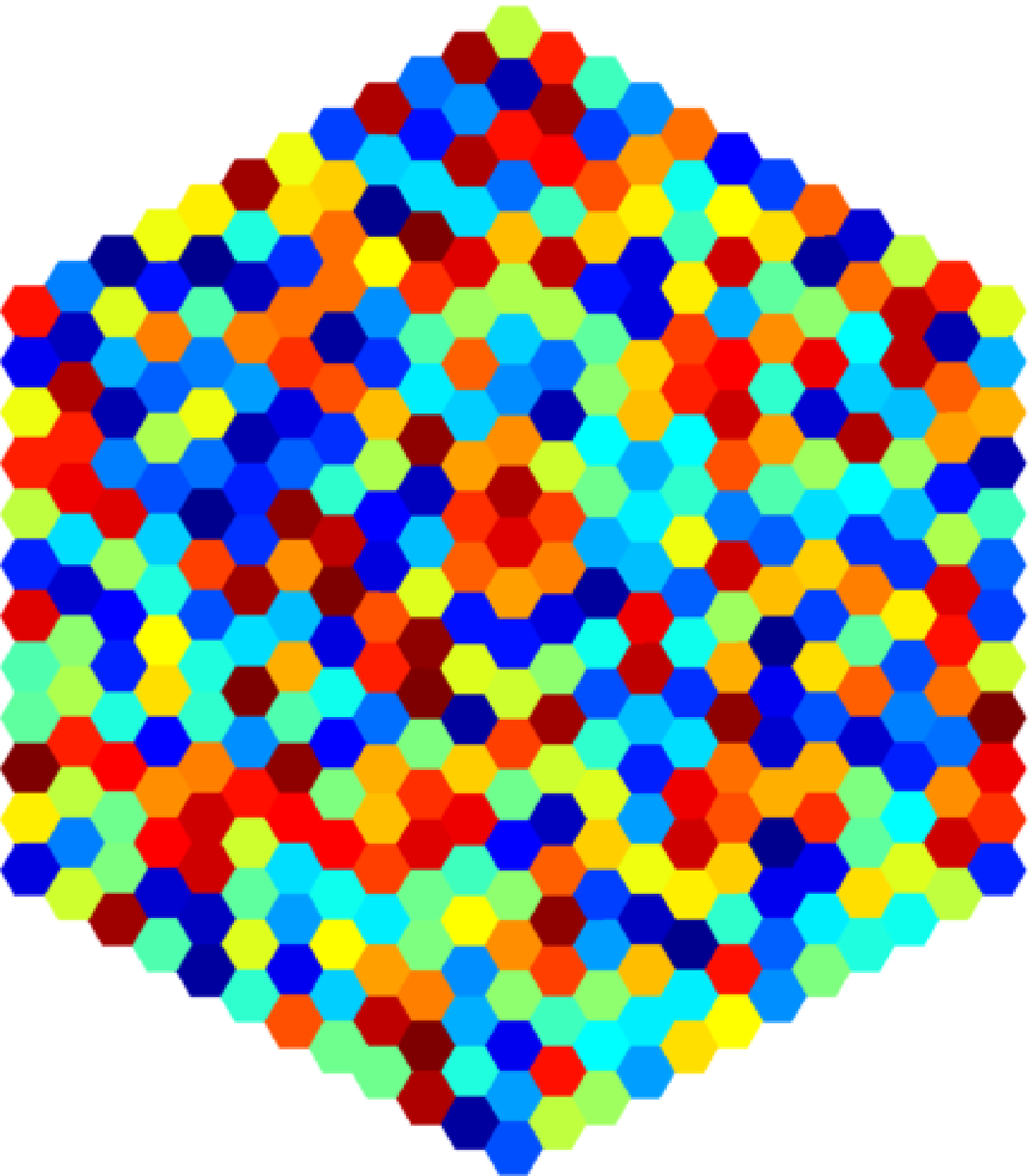}
\includegraphics[width = 0.2 \textwidth]{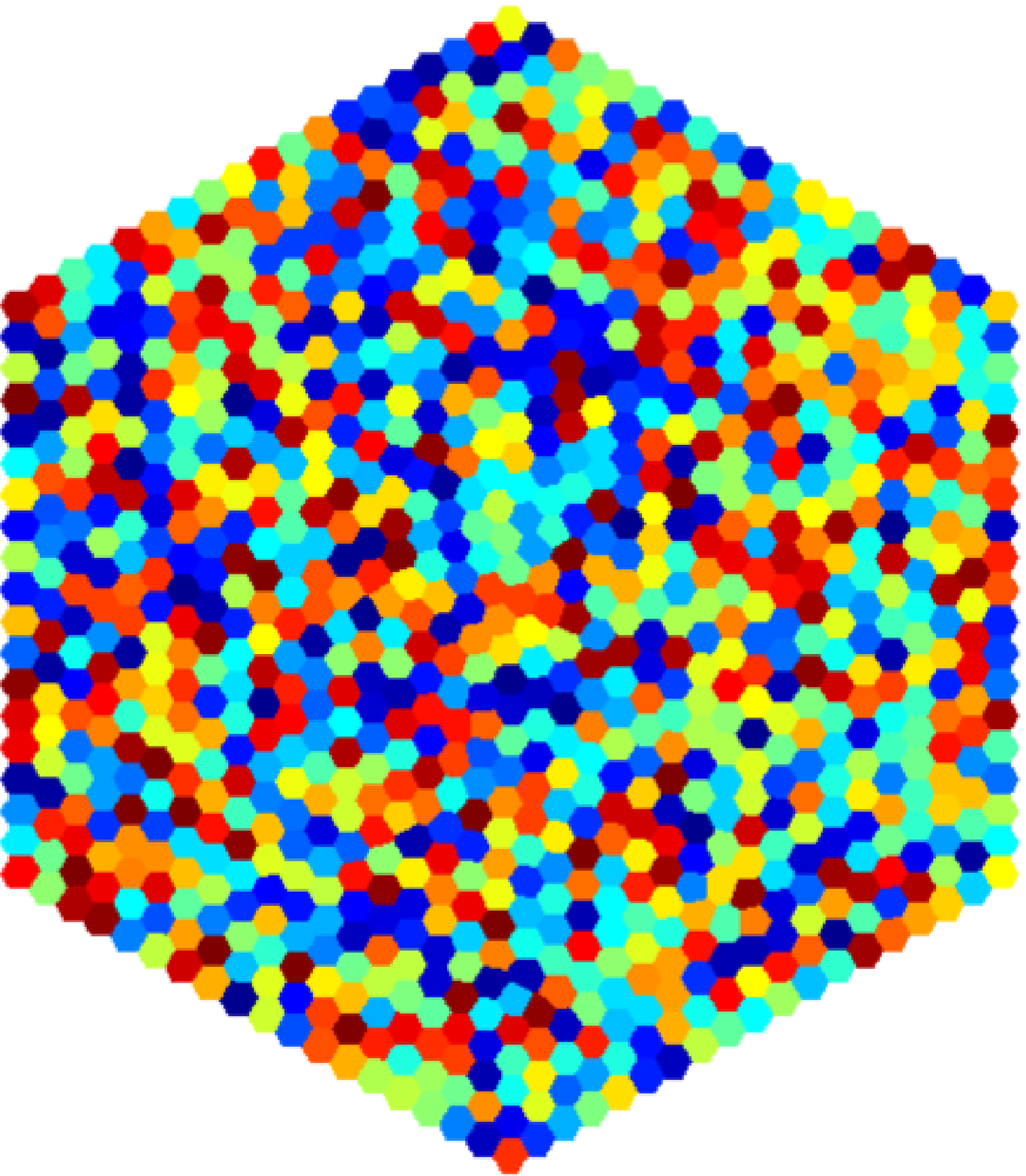}
\caption{Optimal partitions on unions of hexagons for $217, 271, 397$ and $1027$ cells. The partition of $397$ cells is computed on a grid of size $1024\times 1024$ and the partition into $1027$ cells on a grid of size $510\times 510$.}
\label{union_hex}
\end{figure}

We start by presenting the example of the square with $512$ cells on a $512\times 512$ grid. We take random densities on a $32\times 32$ grid and we perform an optimization on this grid in order to get an initial condition. We double the size of the grid points along each axis direction and we interpolate the previous result on the new grid. We optimize again the result and we double the grid size until we reach the desired resolution. In Figure \ref{BBO} we present the result obtained using the approach given above. We observe the presence of patches of regular hexagons as expected. We underline here that this test was made on a laptop with a 3.5GHz quad-core processor and $16$GB of RAM in a few hours of computation time. Having more cells and better resolution is possible, as it can be seen in Figures \ref{BBO}, \ref{union_hex}. Working with sparse matrices we optimize the memory cost used for storing the computational structure. For example in Figure \ref{BBO}, we present a case with $1000$ cells on a grid of size $1000\times 1000$ the computation takes about $12$ hours and RAM consumption is less than $16$GB when using sparse matrices to represent the cells. The $1000\times1000$ computation was made on a machine with an $8$-core Xeon processor with $32$ GB of RAM. However, the code was not parallelized, so eigenvalues were not computed in parallel. The RAM usage was of about $4$GB in general, with spikes up to $12$GB when performing the iteration step, probably due to the fact that several variables containing information of the size of the whole partition structure were present: among these one contains the current densities and one contains the gradient. Also, additional variables were created during the projection step. This example is presented to show the advantage of the grid restriction method with respect to the algorithm presented in \cite{buboou}, where for the same computation a supercomputer was needed\footnote{The computation in \cite{buboou} was made at Texas Advanced Computing Center}. With a fully parallelized code and a super computer the improved algorithm could go even further.

It is possible to observe exact hexagonal patterns if we look at partitions of 2D shapes which are exact union of regular hexagons. These types of domains have also been investigated by Bonnaillie-No\"el, Helffer and Vial in \cite{BNHV} by adding cracks with Dirichlet conditions. In the computations presented here, we always start with an initial condition consisting of random densities and we arrive at an exact partition made of regular hexagons. This is further evidence that the Caffarelli-Lin conjecture seems to be true. The partitions obtained are given in Figure \ref{union_hex}.

\begin{figure}
\centering
\includegraphics[width=0.18\textwidth]{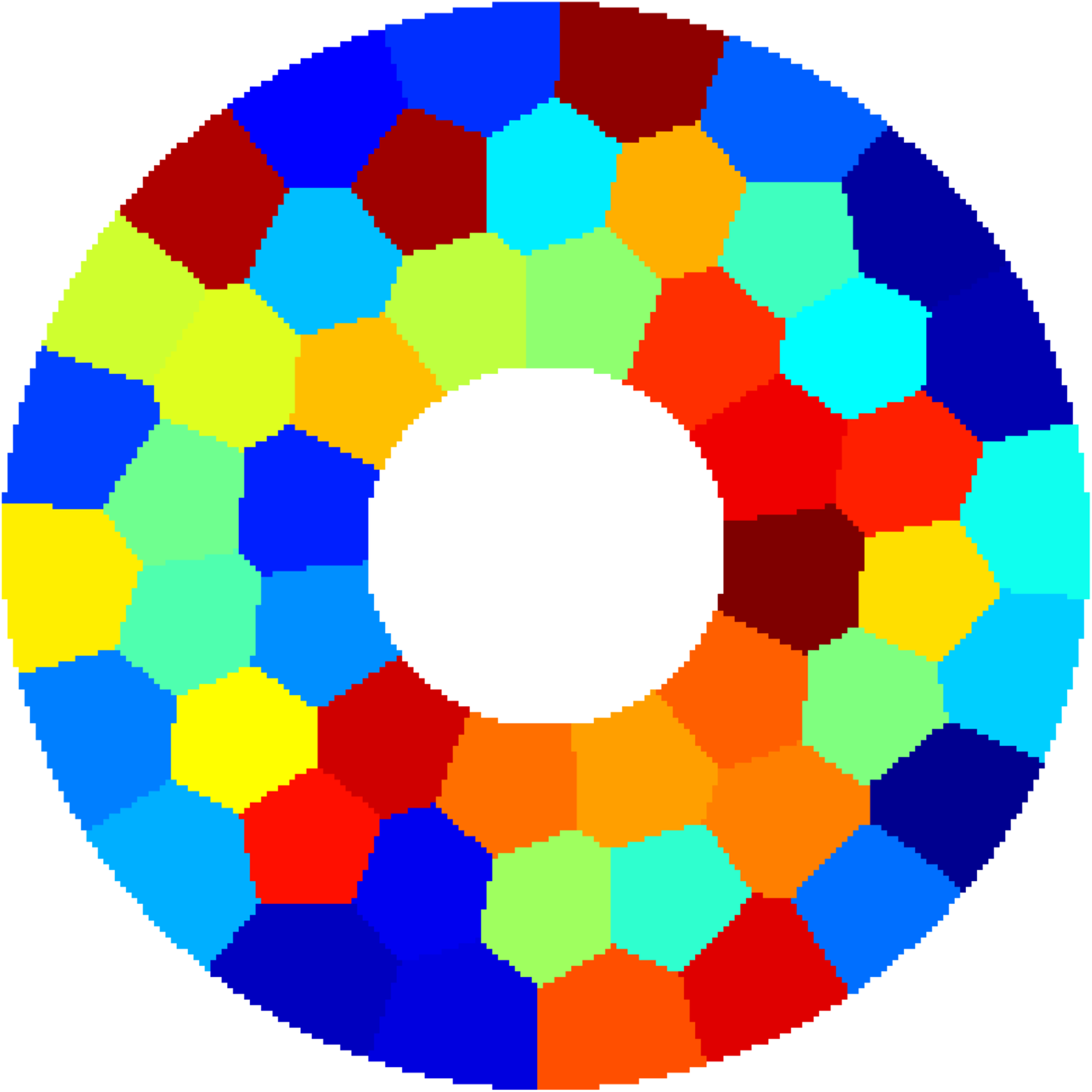}~
\includegraphics[width=0.19\textwidth]{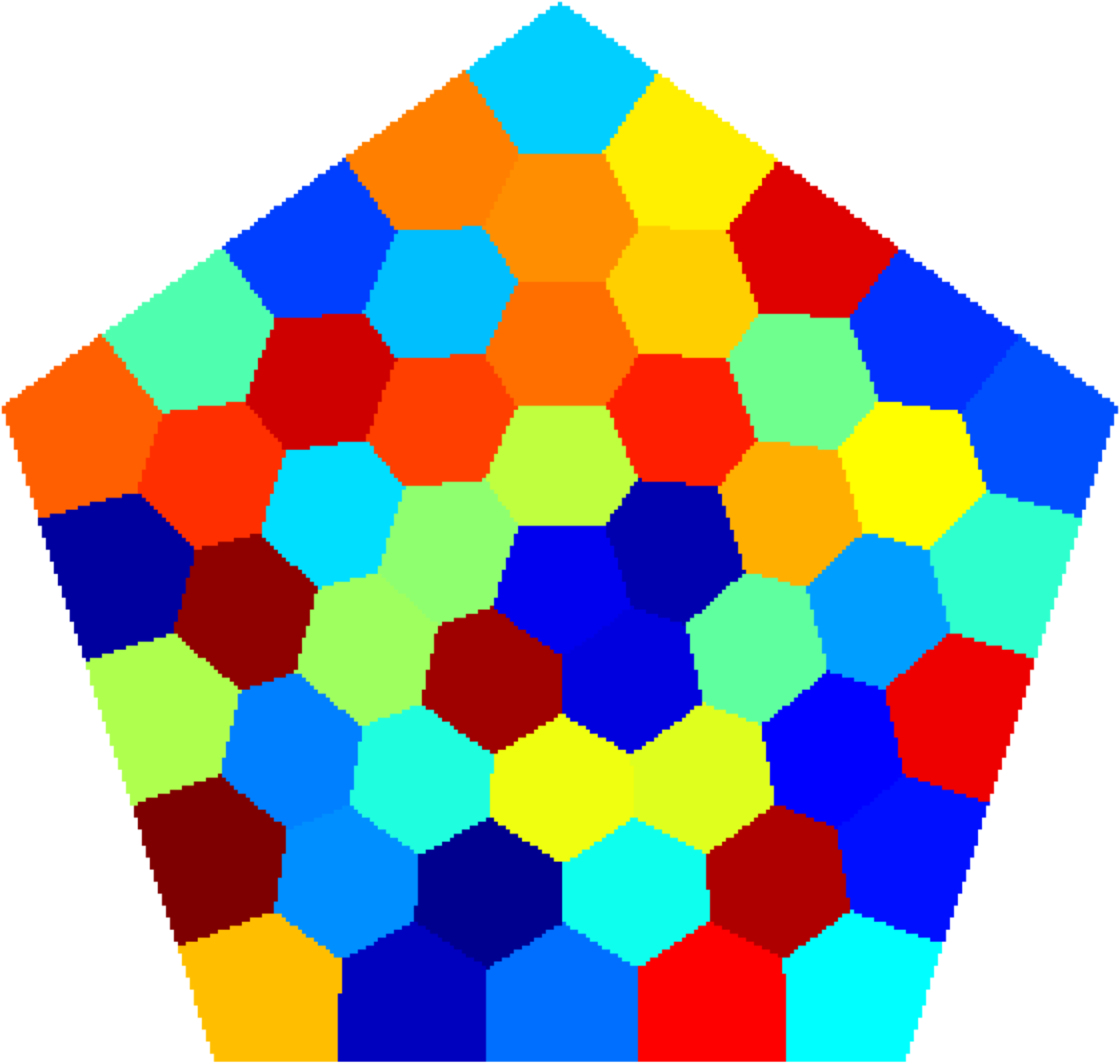}~
\includegraphics[height=0.17\textwidth]{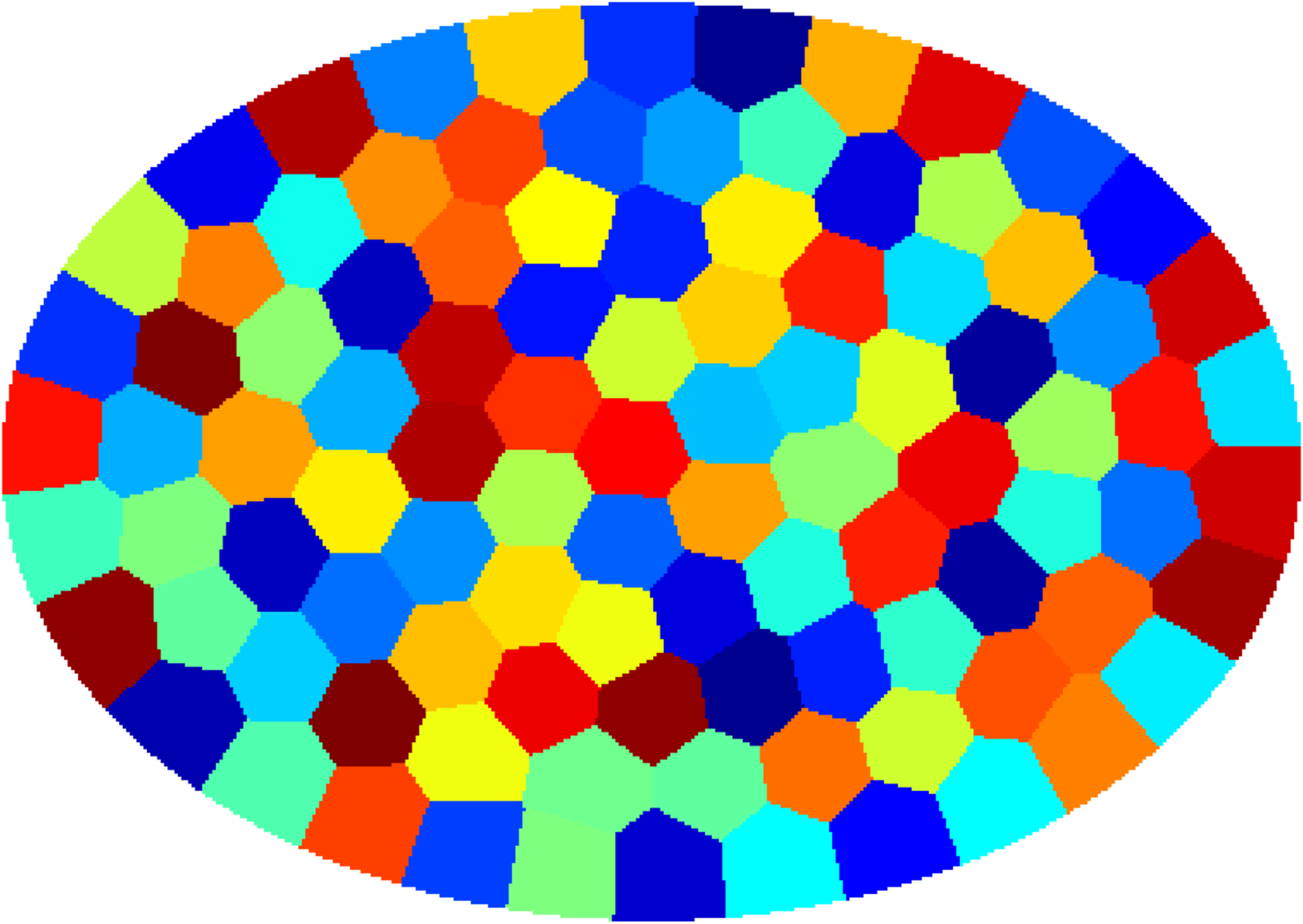}~
\includegraphics[width=0.19\textwidth]{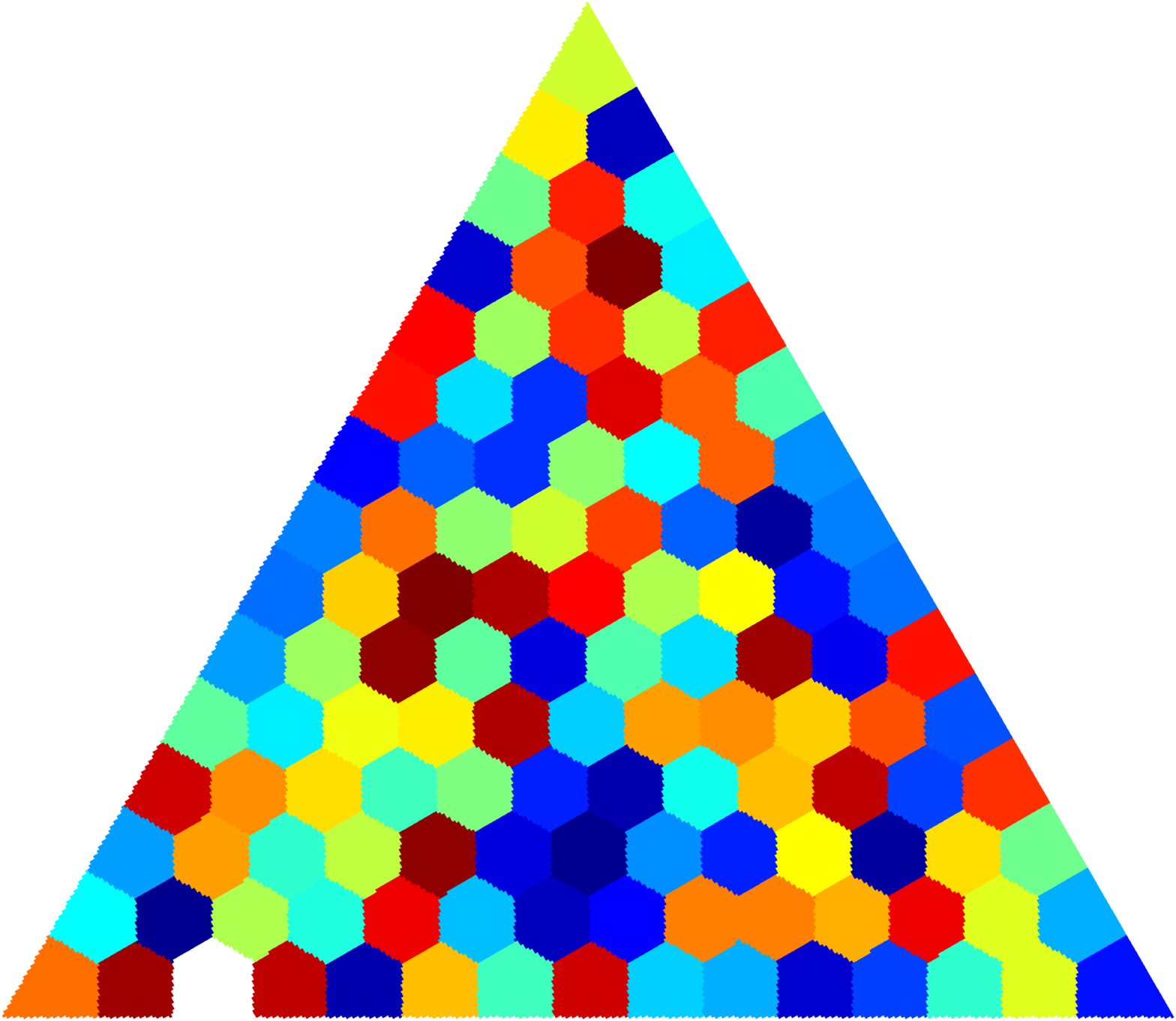}
\caption{Partitions on various domains}
\label{various}
\end{figure}

We may also work on general geometries like the disk, the equilateral triangle, ellipses, polygons, etc. Note that even if we need to consider a computational region which is rectangular and contains our domain, the computational complexity remains the same when we use the grid restriction procedure. Some results on various geometries are available in Figure \ref{various}. 

As was already noted in \cite{BoBN16}, for the equilateral triangle we observe that partitions into $k(k+1)/2$ elements seem to be made out of polygons of three types. This is not very surprising if we note that the edges of the equilateral triangle have the right orientations to make it compatible with hexagonal patterns. More precisely, for numbers of cells $n$ of the form $n = k(k+1)/2$, i.e. triangular numbers, we obtain cells which are polygons of three types: three quadrilaterals corresponding to the corners, $3(k-2)$ pentagons on the edges and a hexagonal cells in the interior.

In \cite{buve} the authors introduced the multiphase problem
\begin{equation}
 \min \sum_{i=1}^n \lambda_1(\omega_i) +\alpha|\omega_i|
 \label{multiphase}
\end{equation}
where the minimization is made in the class of families $(\omega_i)_{i=1}^n$ of non-overlapping subsets $\omega_i$ of the domain $D$. Note that in this case the $\omega_i$ do not form a partition of $D$. In fact in \cite{buve} the authors show that as soon as $\alpha>0$ there are no triple points and there are points of $D$ which are not covered by any of the sets $\omega_i$. In \cite{bove-multiphase} the authors proved that if $D$ is Lipschitz, then there are no singular points on $\partial D$ and they proposed a numerical algorithm for the study of such multiphase configurations. The numerical algorithm should take into account the non-overlapping condition and this was achieved by adding a supplementary phase which represents the part of $D$ not covered by any of the sets $\omega_i$. This allows the use of a similar algorithm, since we turn the non-overlapping condition in to a partition condition. The grid reduction procedure presented before can also be applied for the multiphase problem, leading to an acceleration of the optimization algorithm for problem \eqref{multiphase} presented in \cite{bove-multiphase}.

The reason we mention problem \eqref{multiphase} here is to underline the potential connection between the spectral partitioning problem, corresponding to $\alpha=0$ in \eqref{multiphase}, and the circle packing problem. The circle packing problem seeks the maximal radius of $n$ disks which can be placed in $D$ without overlap. Since for given $\alpha>0$ the one phase problem 
\[ \min \lambda_1(\omega)+\alpha|\omega|\]
is solved by a disk of radius $r_\alpha = \left(\frac{\lambda_1(B_1)}{\alpha \pi}\right)^{1/4}$, we see that for $\alpha$ large enough, $n$ disks of radii $r_\alpha$ fit into the domain $D$ (which is assumed open and Lipschitz). Moreover, as $\alpha$ decreases, the optimal radii $r_\alpha$ increase. There is a critical value $\alpha_0$ for which $r_{\alpha_0}$ is equal to the maximal radius of a family of $n$ disks which fits into $D$. We see that in this way, the spectral optimization problem \eqref{multiphase} includes the circle packing problem for a particular value of $\alpha$.
\begin{figure}
\centering
\includegraphics[width=0.3\textwidth]{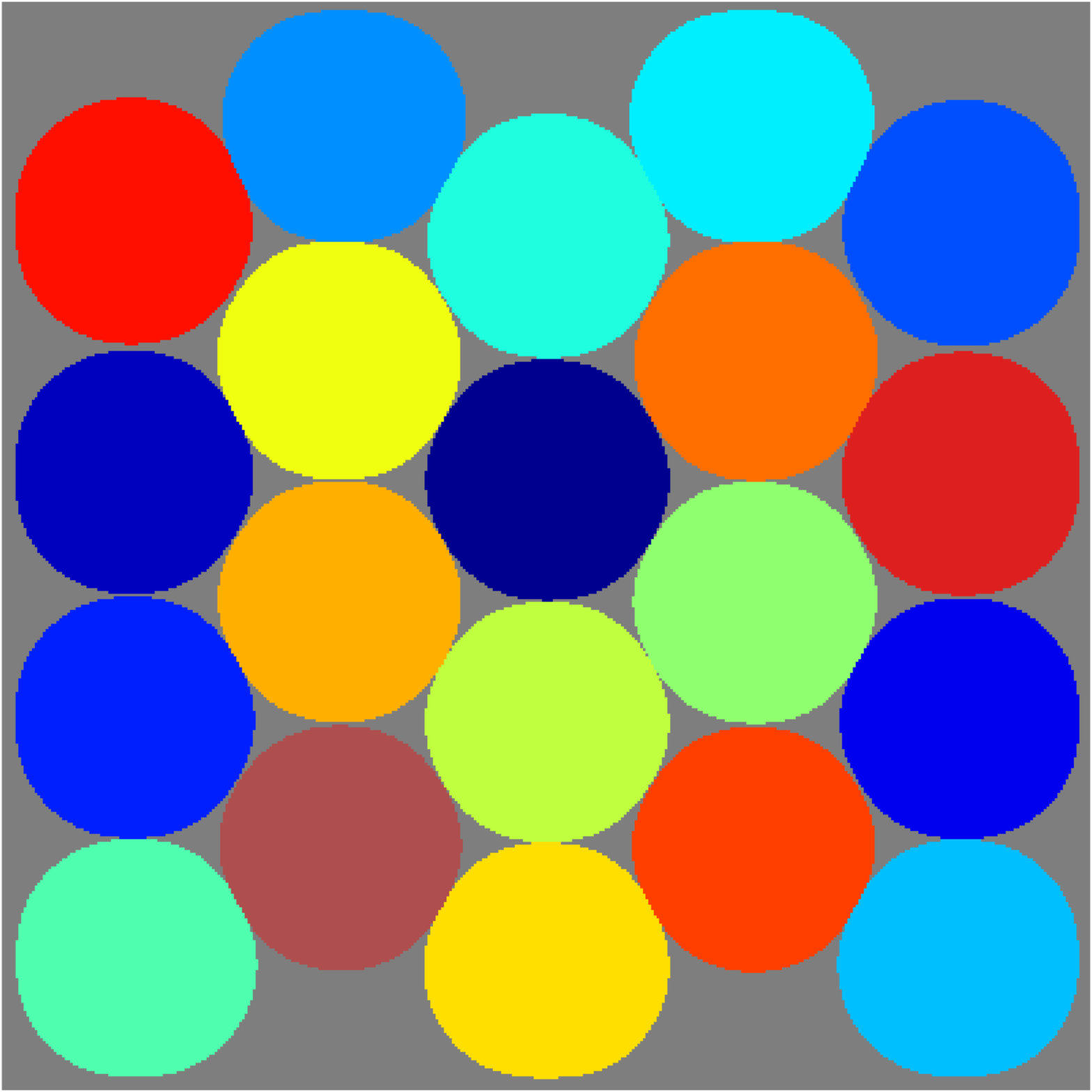}~
\includegraphics[width=0.34\textwidth]{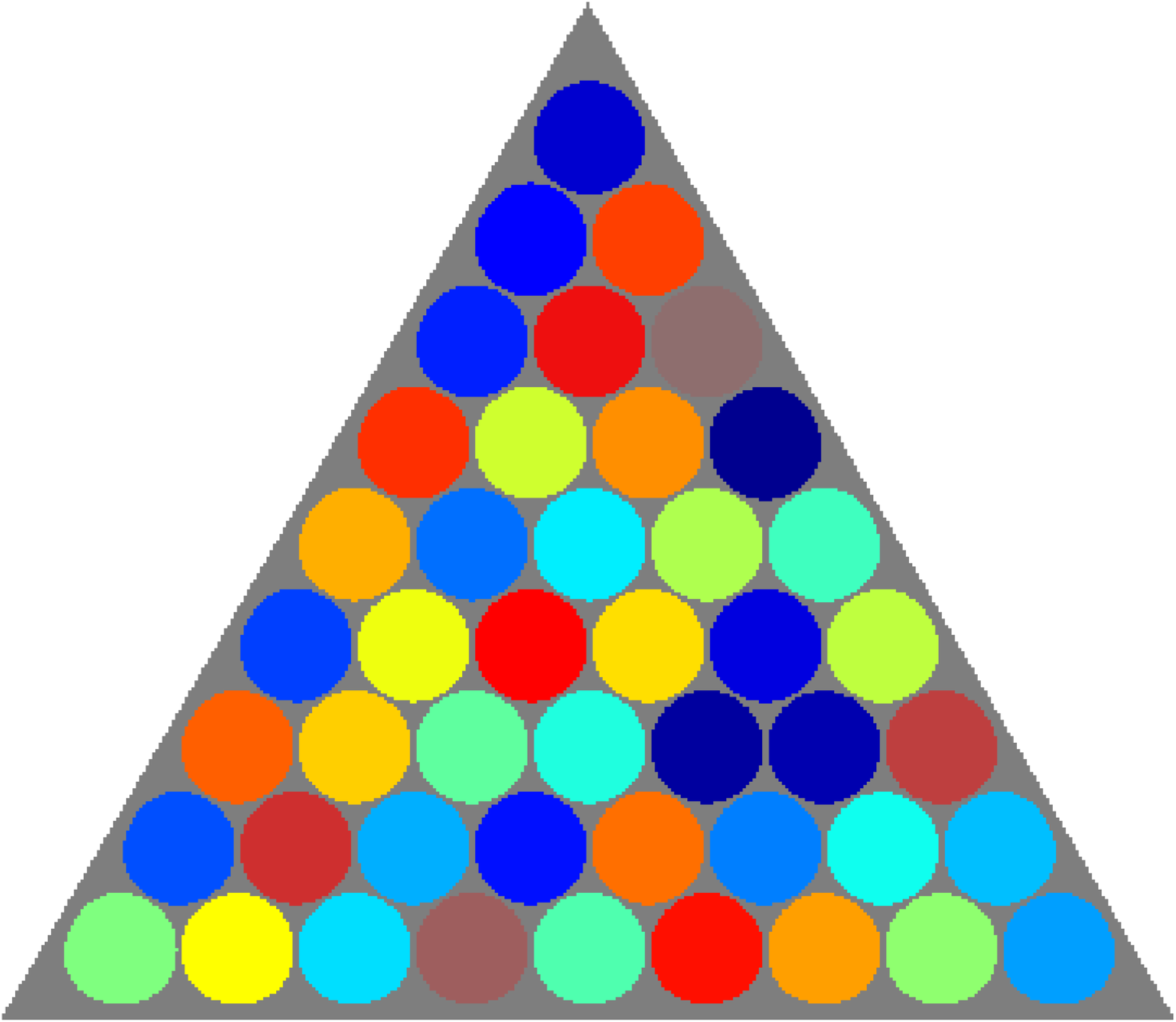}
\caption{Examples of computations of circle packings}
\label{comp_cpack}
\end{figure}

\begin{figure}
\centering
\includegraphics[width=0.15\textwidth]{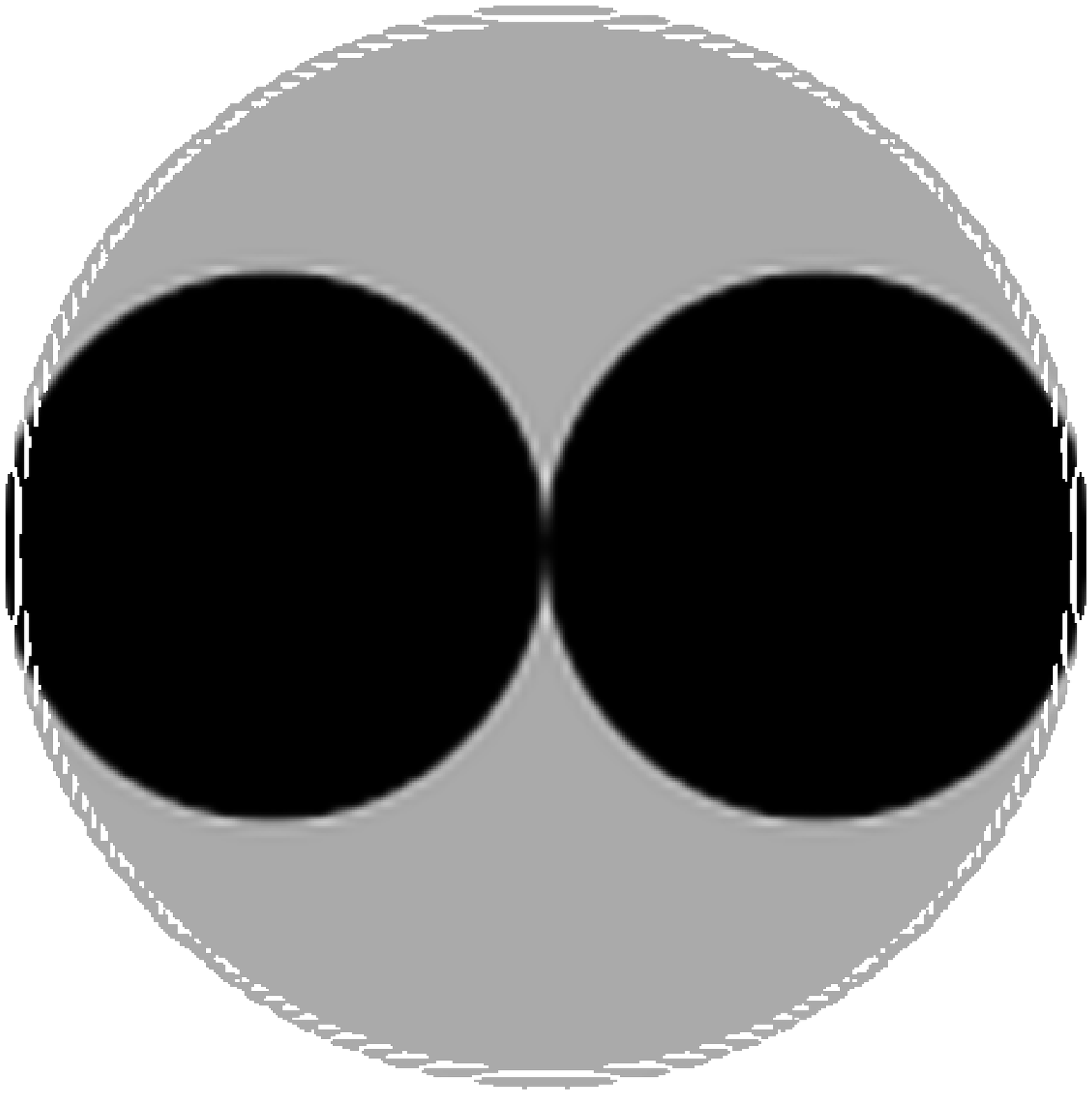}~
\includegraphics[width=0.15\textwidth]{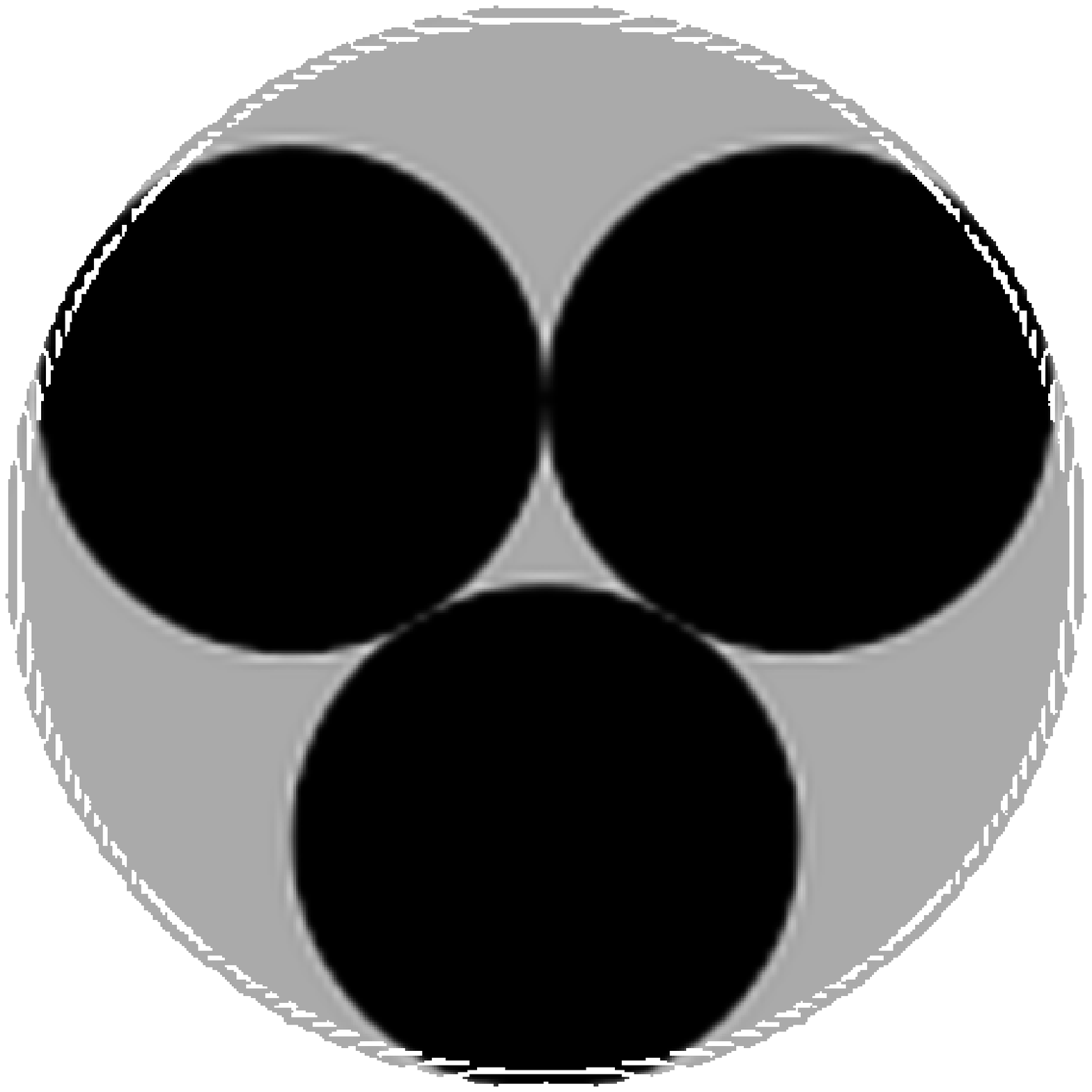}~
\includegraphics[width=0.15\textwidth]{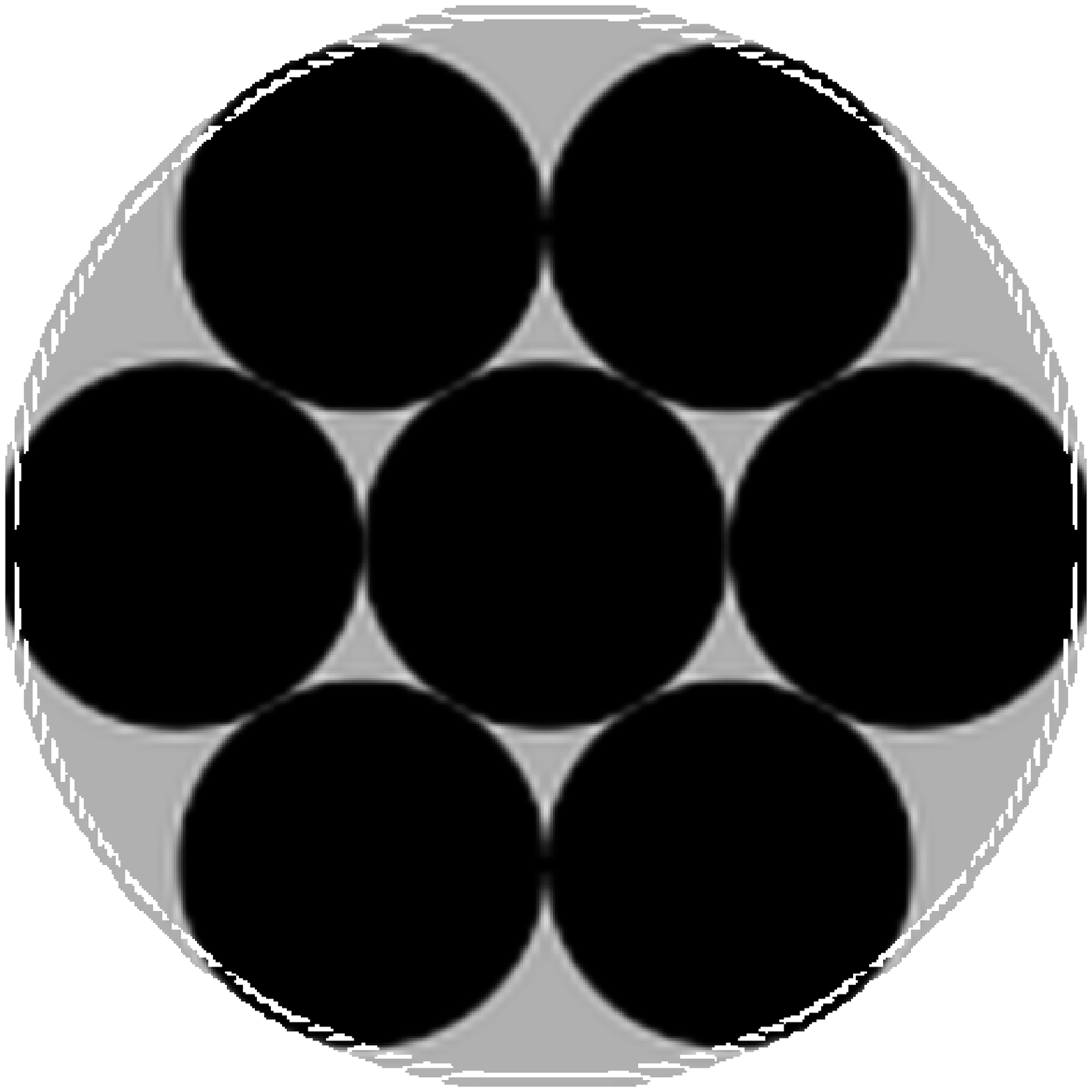}~
\includegraphics[width=0.15\textwidth]{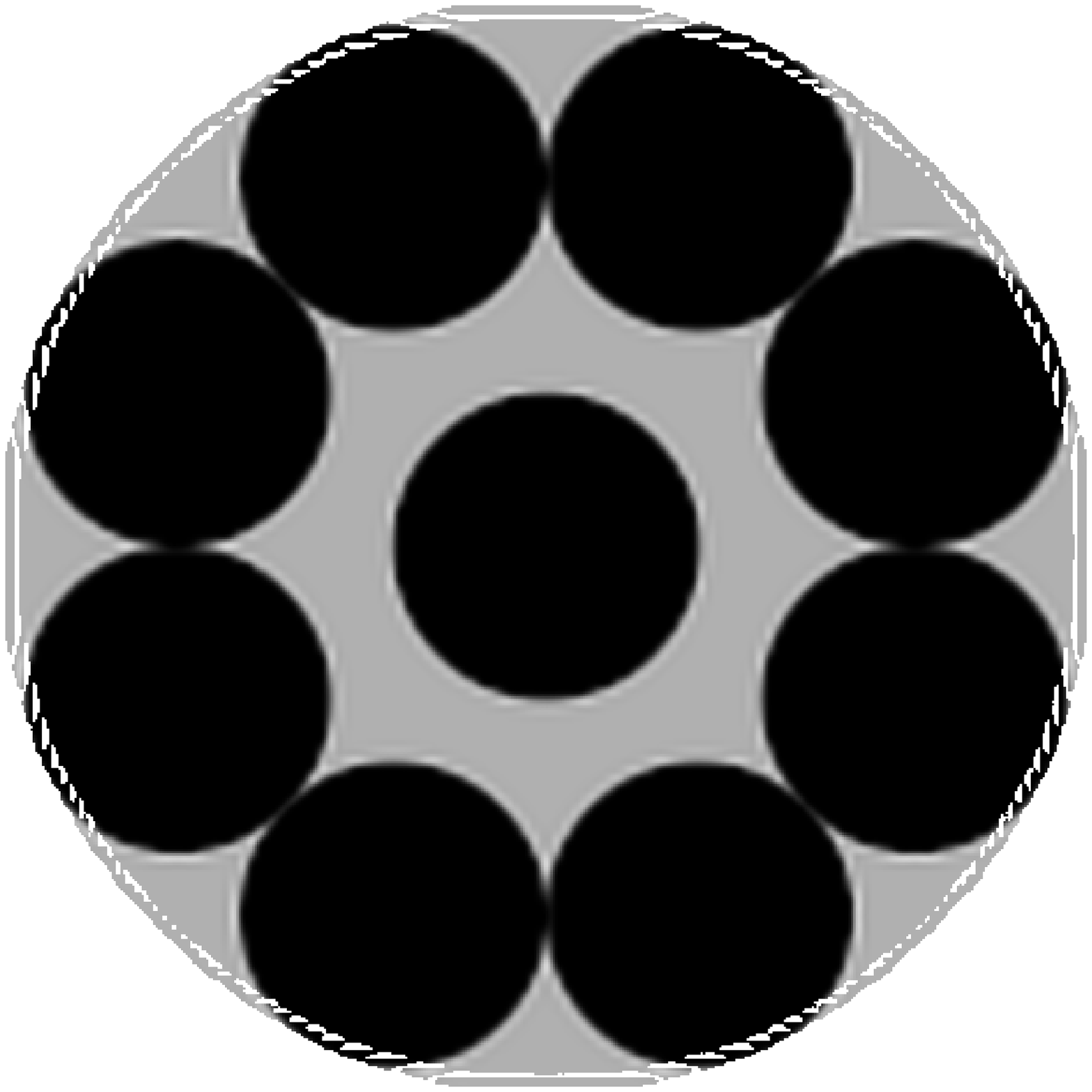}~
\includegraphics[width=0.15\textwidth]{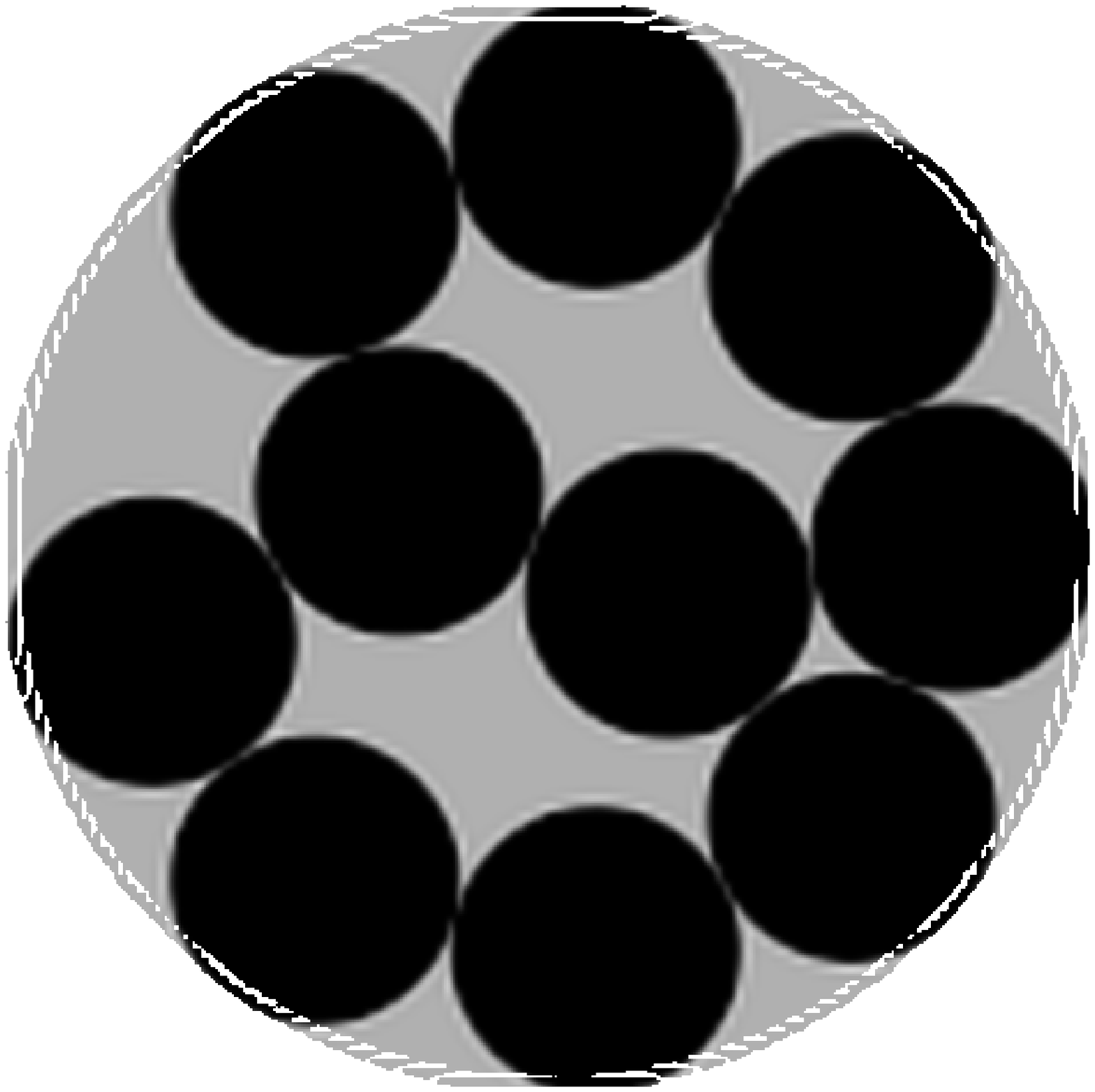}~

\includegraphics[width=0.15\textwidth]{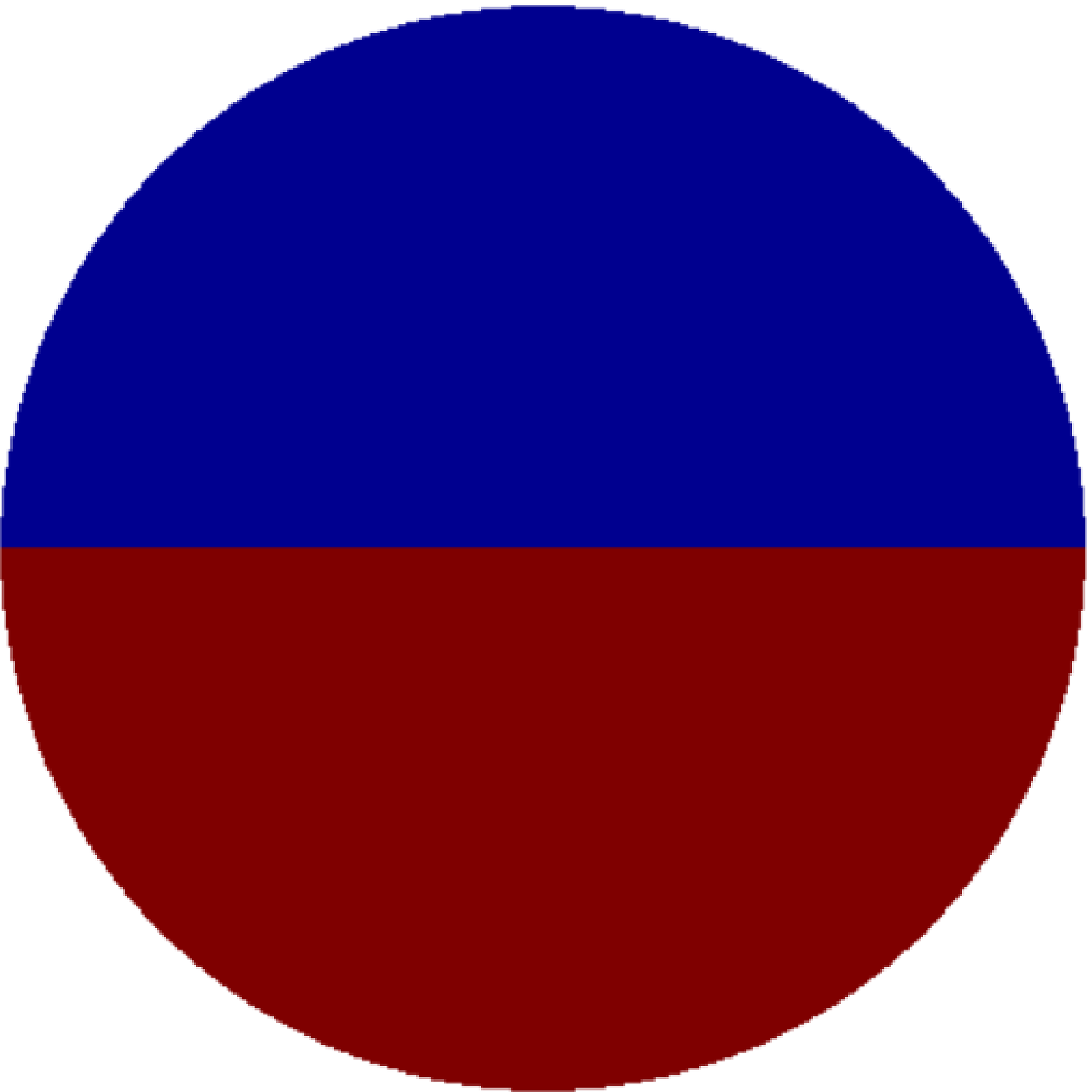}~
\includegraphics[width=0.15\textwidth]{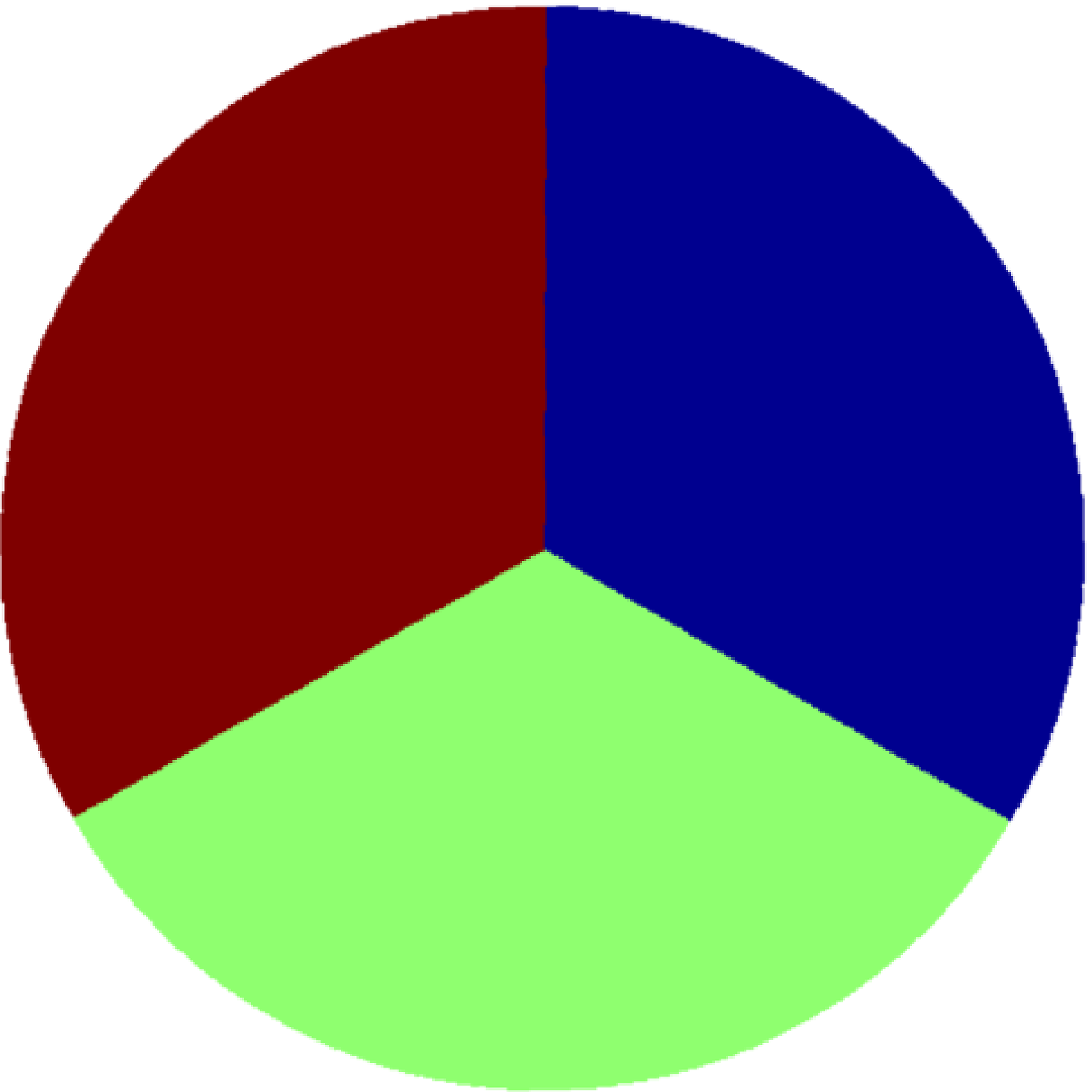}~
\includegraphics[width=0.15\textwidth]{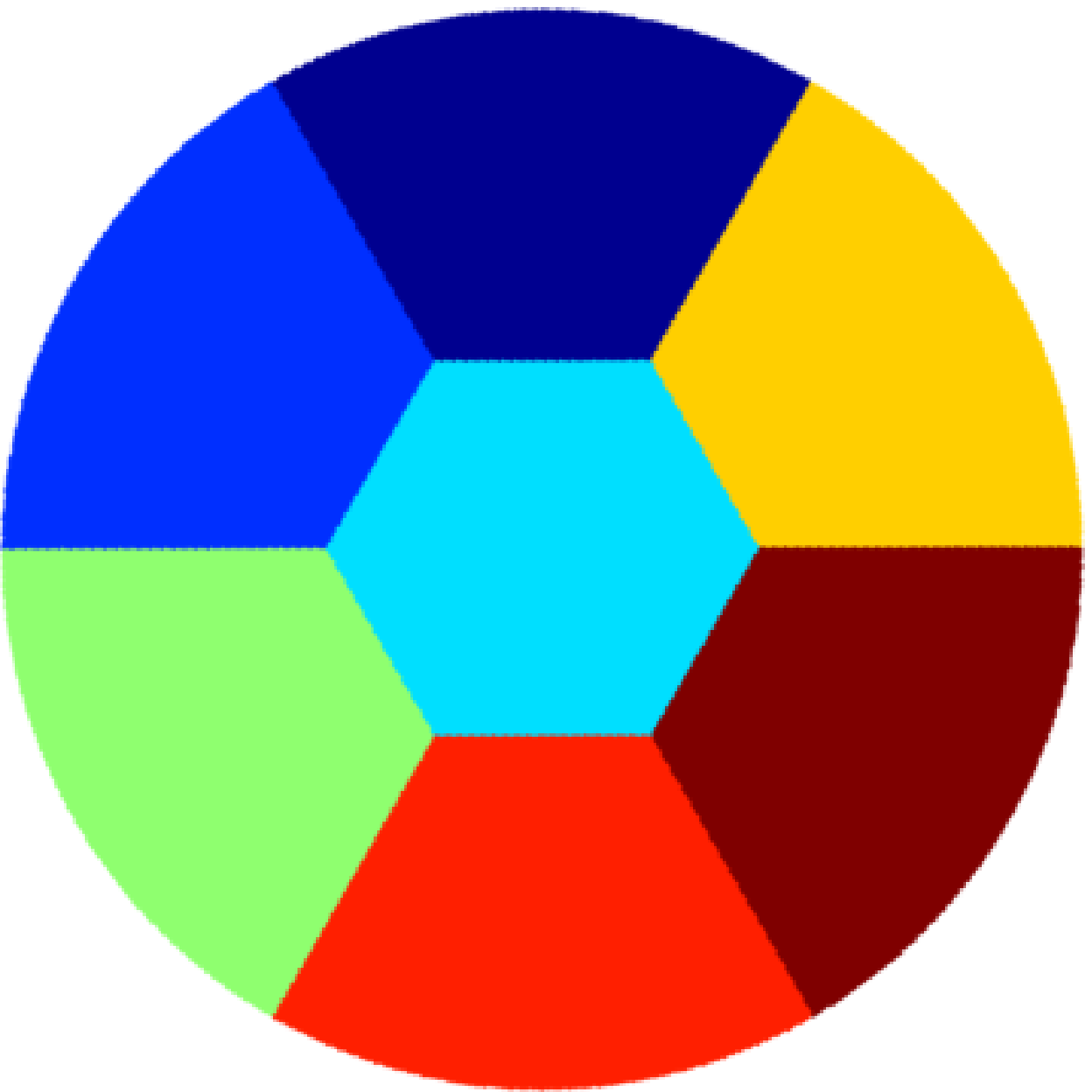}~
\includegraphics[width=0.15\textwidth]{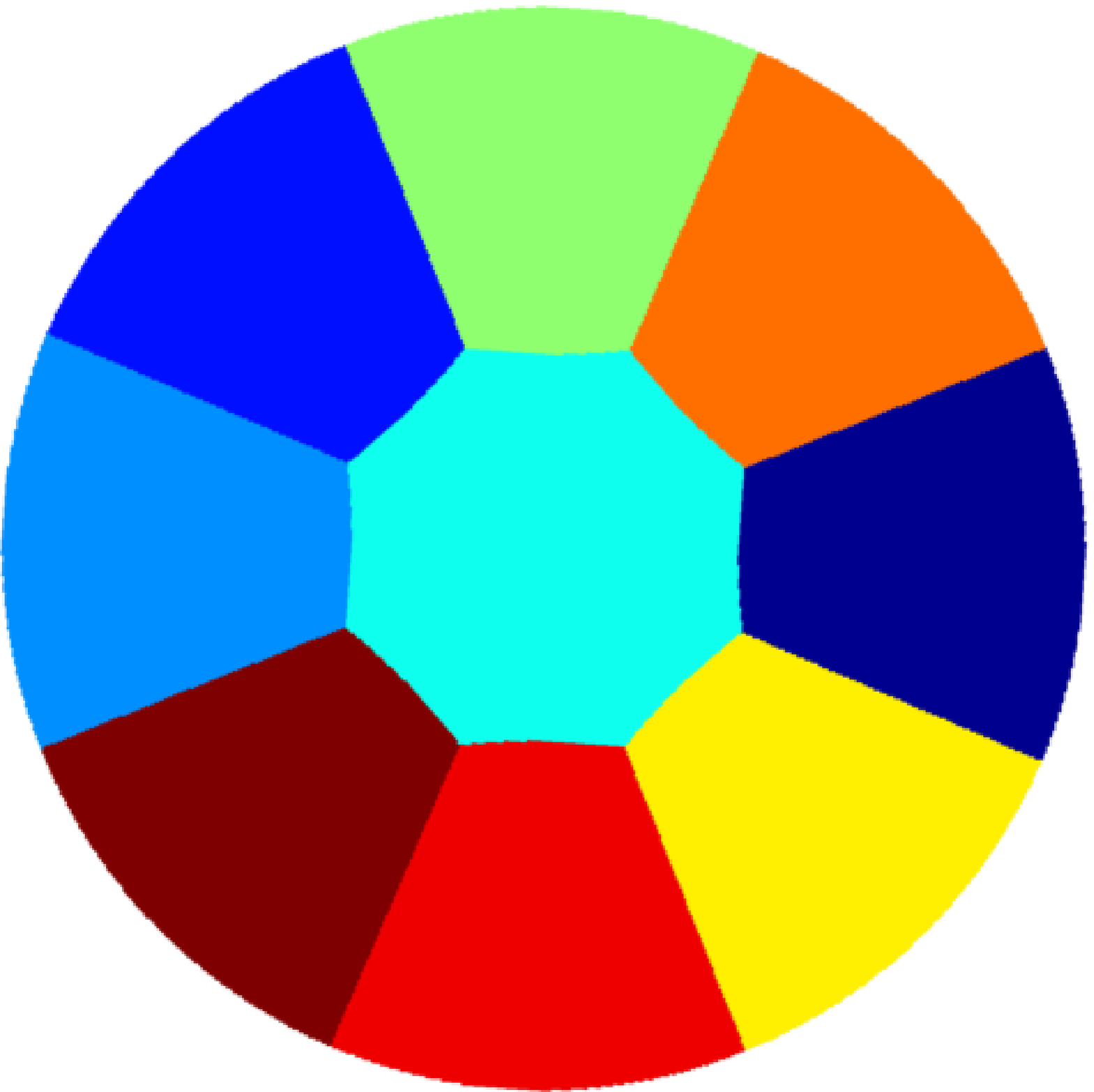}~
\includegraphics[width=0.15\textwidth]{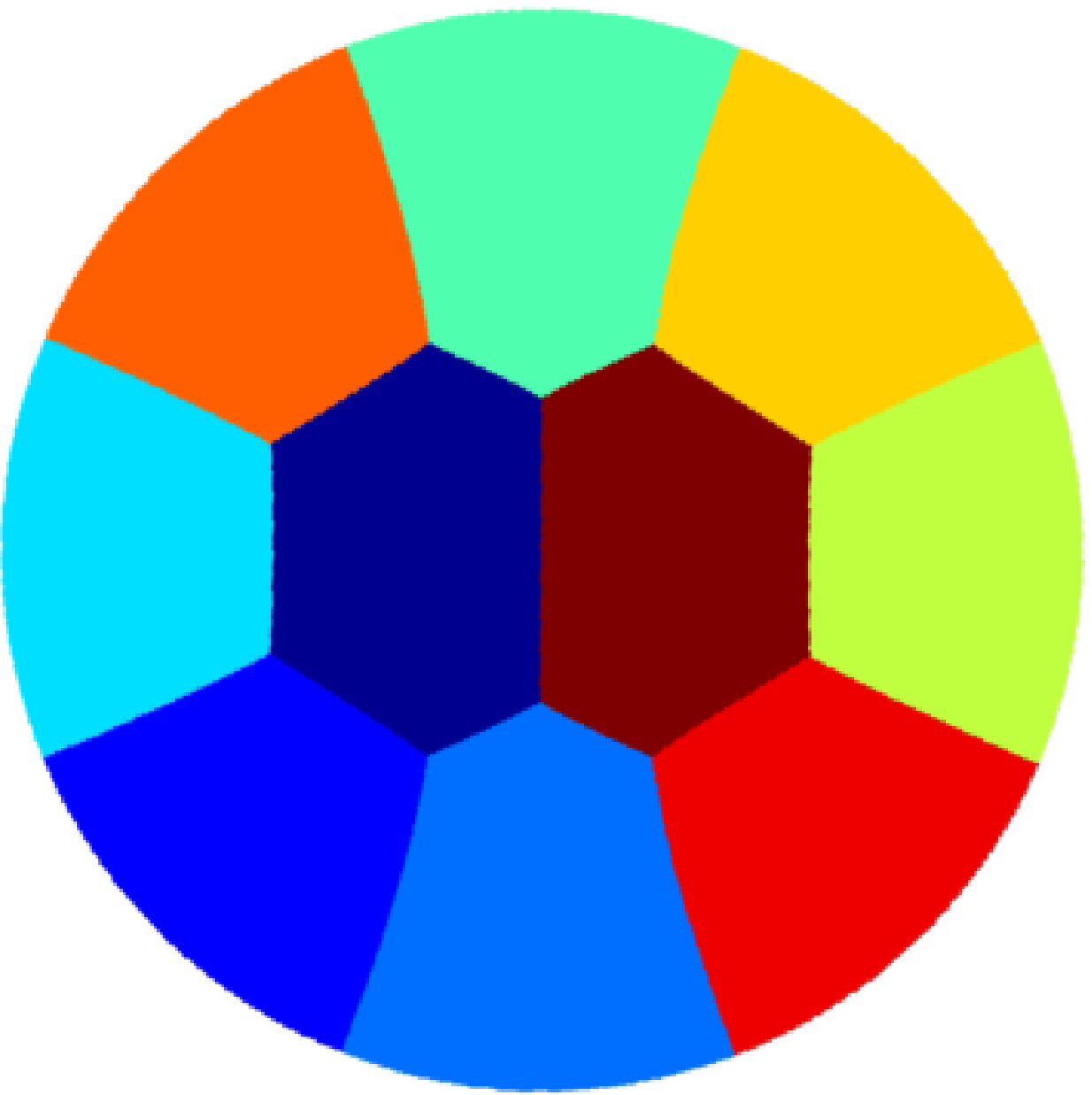}~
\caption{Similarities between optimal circle packing configurations and minimal spectral partitions}
\label{similarity_cpack}
\end{figure}

At a first sight it may seem that using a spectral formulation for the computation of the optimal circle packing is not feasible numerically. However, using the optimized algorithm presented in this article, we may approach successfully the optimal circle packing problem in certain cases. Some examples are shown in Figure \ref{comp_cpack}. One may check these results against best known ones presented at \href{http://www.packomania.com}{packomania.com} Another reason why problem \eqref{multiphase} was recalled here is to underline the similarity between some circle packing configurations and the corresponding optimal partition minimizing the sum of eigenvalues. It turns out that in the cases where the symmetries of the circle packing configurations are also symmetries of $D$, the distribution of the cells in the spectral partition is very similar to the one of the disks in an optimal circle packing. See Figure \ref{similarity_cpack} for a few examples. We may wonder if there is a general connection between the circle packing problem and the minimization of the sum of the fundamental eigenvalues. There is a resemblance even in the asymptotic case: it is conjectured that the hexagonal tiling of the plane minimizes the sum of the fundamental eigenvalues. On the other hand it is known that the densest circle packing in the plane is the one corresponding to a hexagonal grid \cite{toth72}.

\subsection{Surfaces}
Initial computations of spectral partitions on surfaces were motivated by Bishop's conjecture \cite{bishop}. It is known that if we consider the partition of the sphere into two cells $(\omega_1,\omega_2)$ then the sum of the first Laplace-Beltrami eigenvalues is minimized in the case of two hemispheres. This is in close relation to the symmetrization techniques used in \cite{FH76} and the monotonicity formulas proved by Alt, Caffarelli and Friedman in \cite{ACFm}. The case of partitions with three cells is still open for the sum, although Bishop \cite{bishop} conjectured for some time now that the $Y$ partition consisting of three slices of the sphere of angle $2\pi/3$ is the optimal one.  The computations in \cite{elliott-ranner} and \cite{beni-fsol} confirmed numerically that it is probable that Bishop's conjecture is valid.  A number of other conjectures arise from these numerical computations, notably those related to regular partitions of the sphere: $4,6$ and $12$ cells. We recall that if instead of the sum we wish to minimize the largest eigenvalue then the optimal partition is known for $n=3$ on the sphere. It is proved in \cite{helffer-3sphere} that the $Y$ partition into three slices of angle $2\pi/3$ is optimal.

\begin{figure}
\centering
\includegraphics[width = 0.3\textwidth]{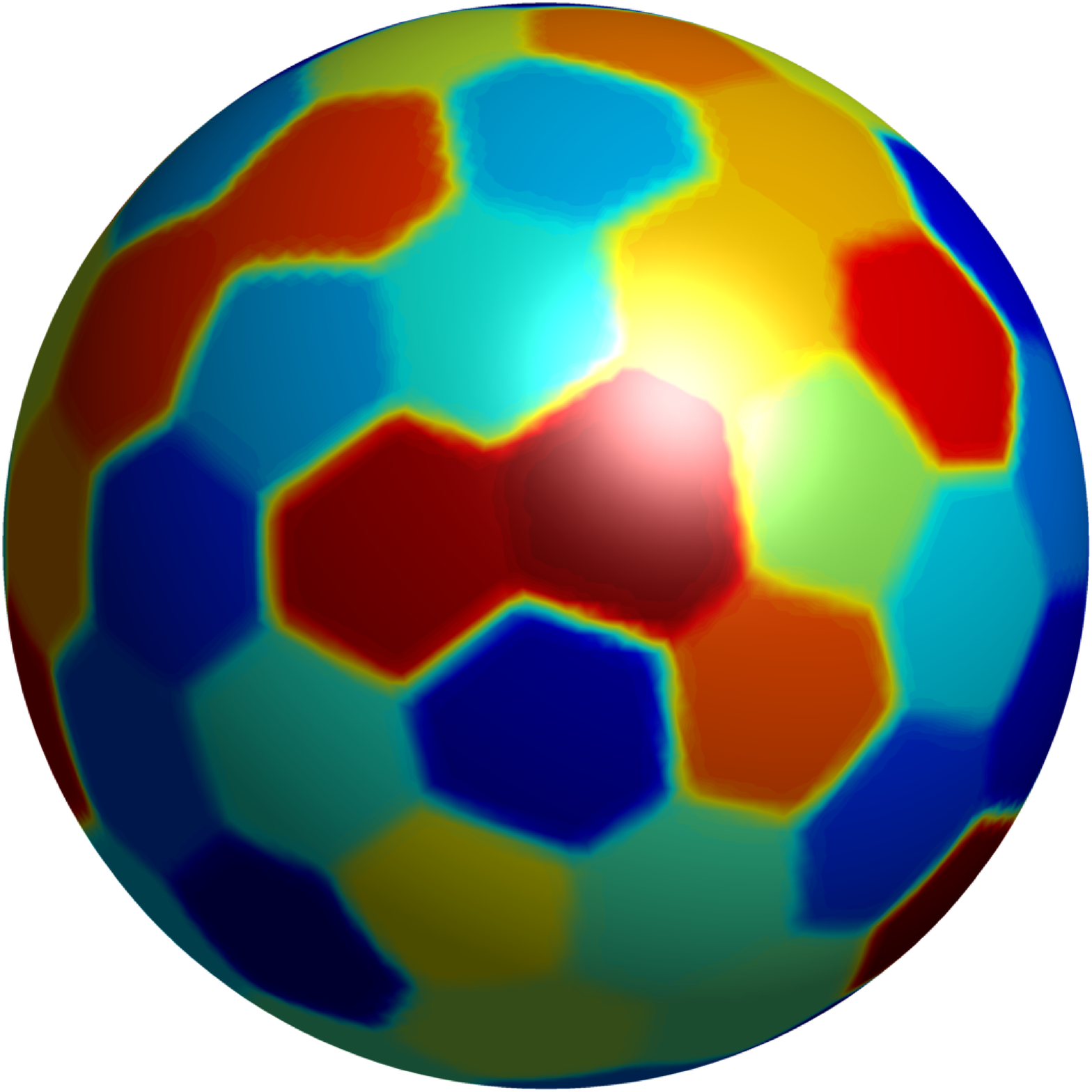}~
\includegraphics[width = 0.3\textwidth]{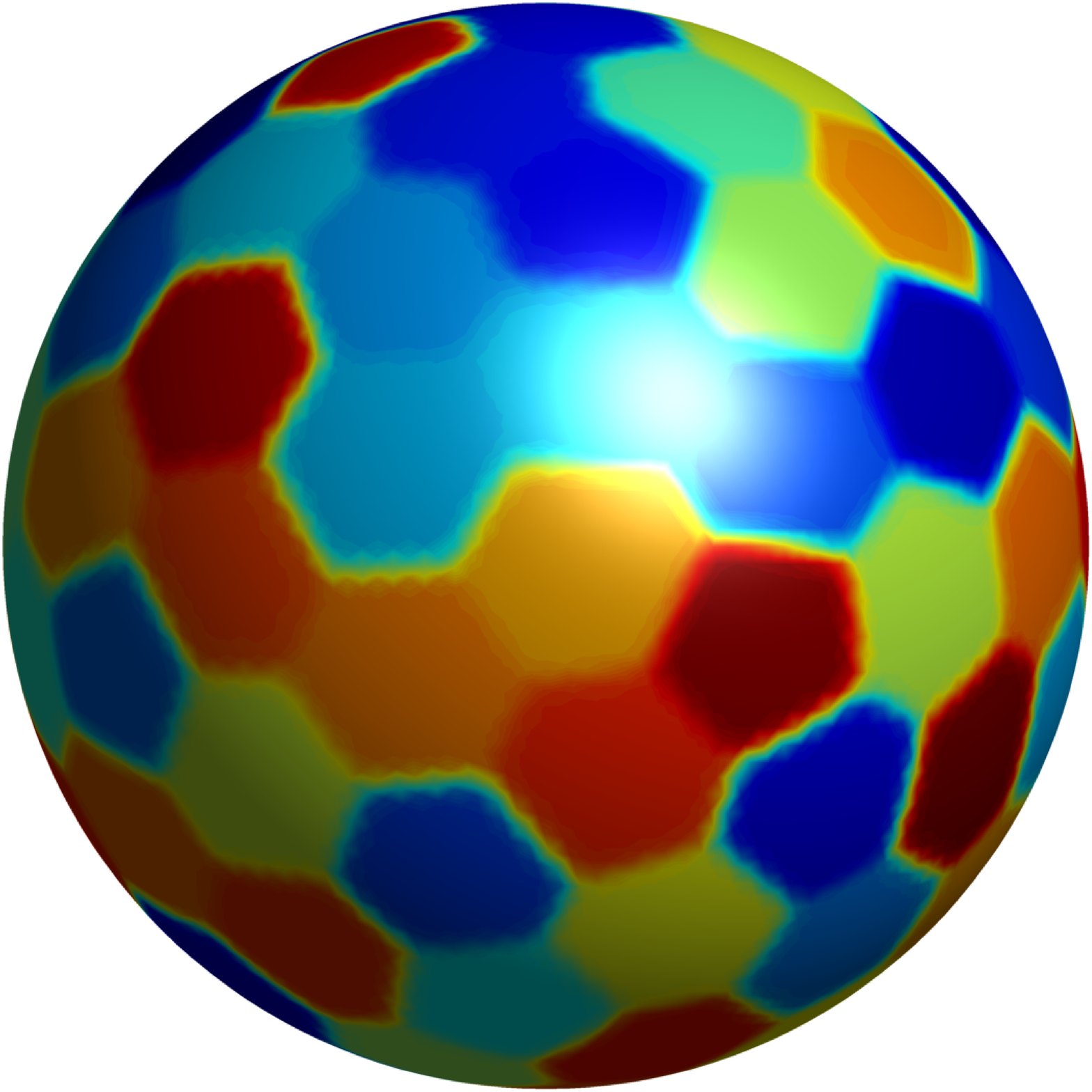}~
\includegraphics[width = 0.3\textwidth]{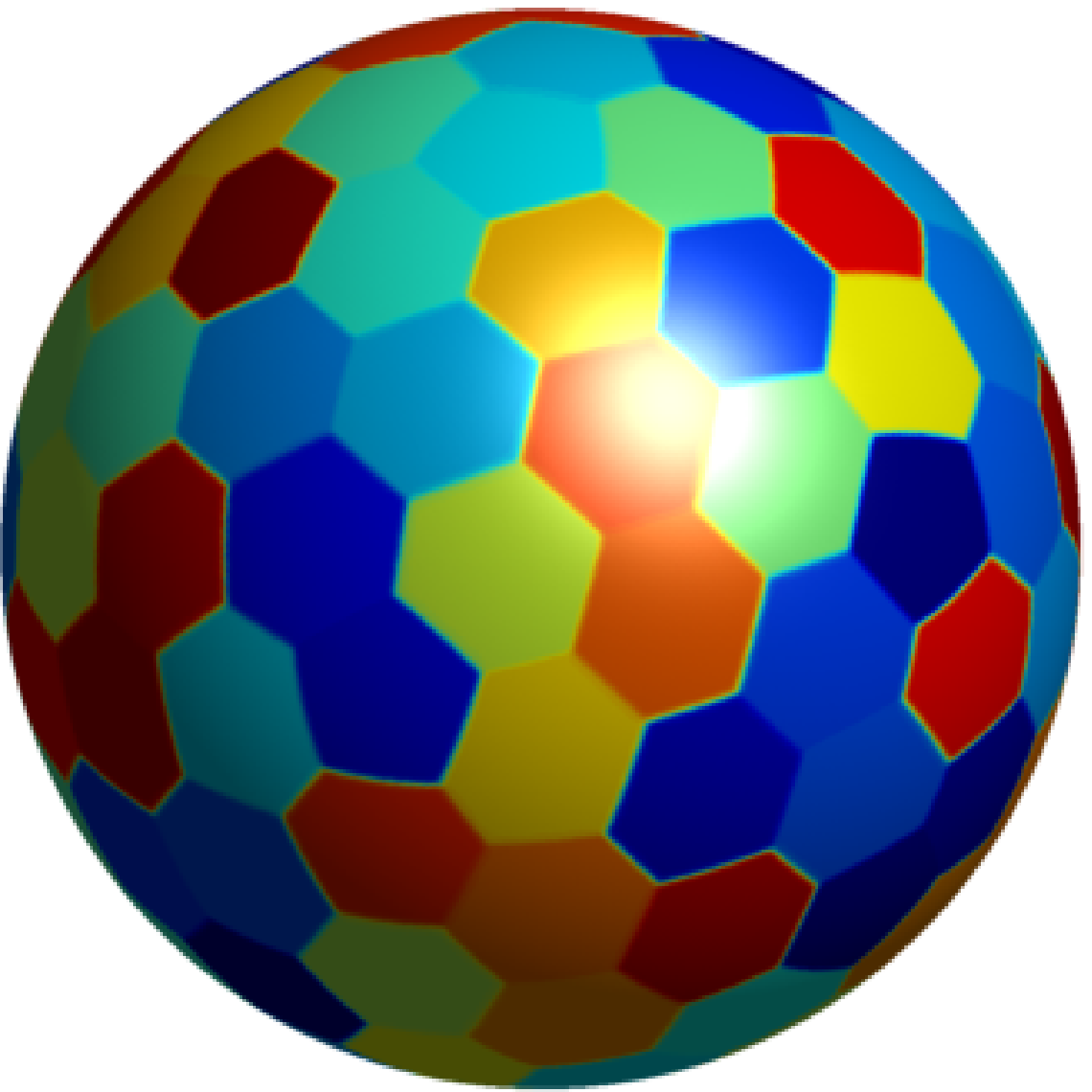}
\caption{Optimal candidates in the case of the sphere. $n \in \{72,90,120\}$}
\label{sphere_bel}
\end{figure}

\begin{figure}
\centering
\includegraphics[width = 0.3\textwidth]{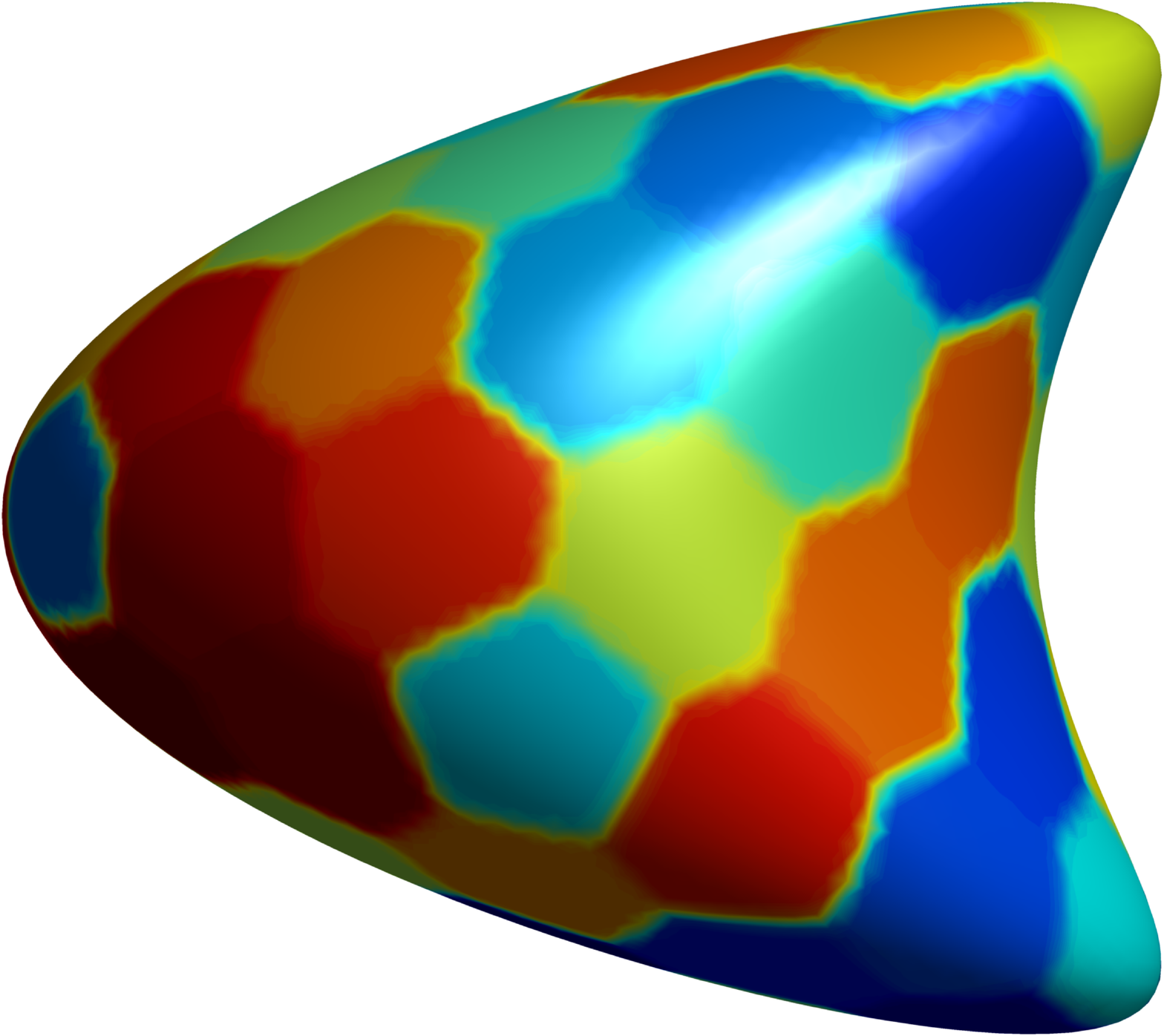}
\includegraphics[width = 0.3\textwidth]{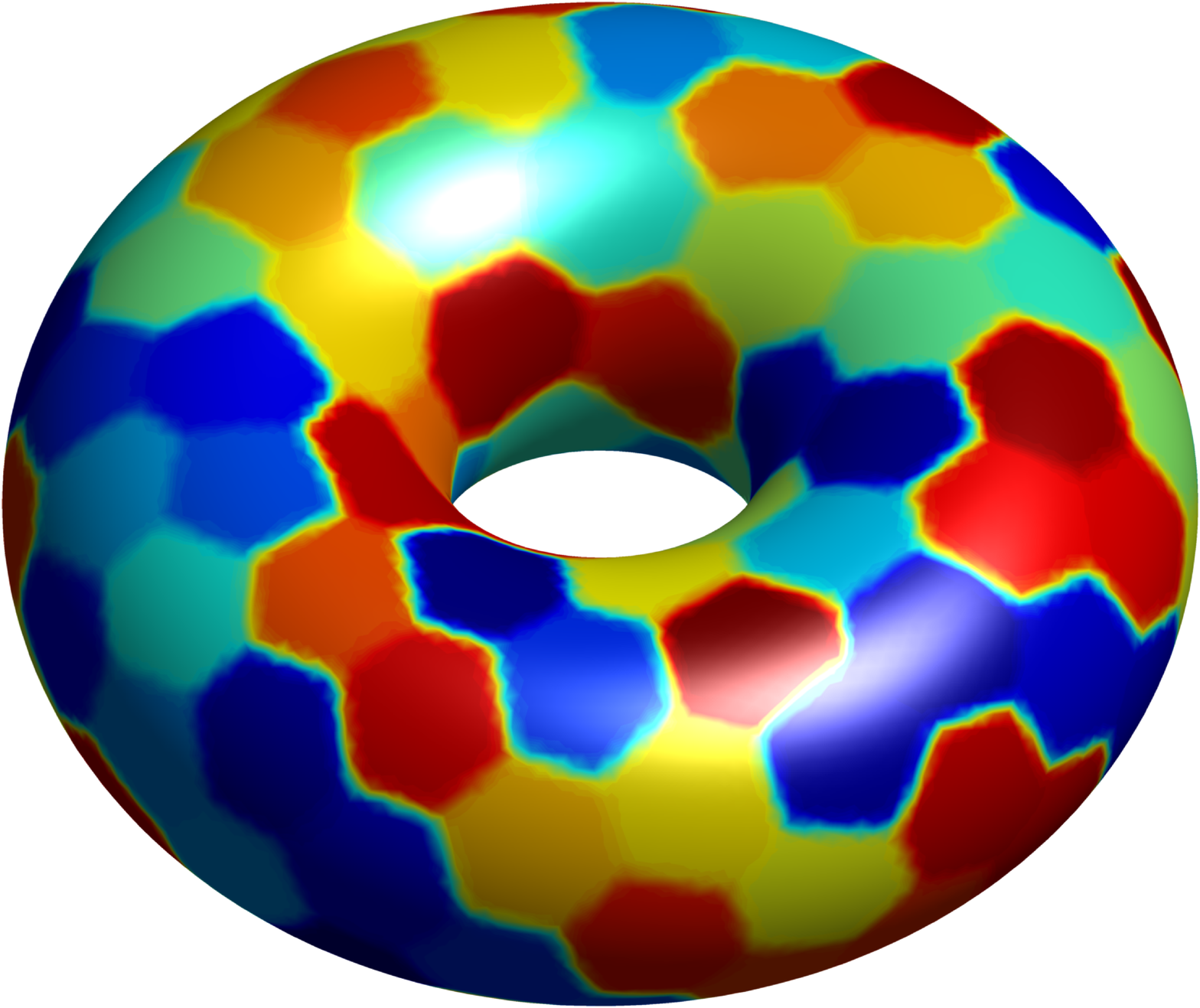}
\includegraphics[width = 0.28\textwidth]{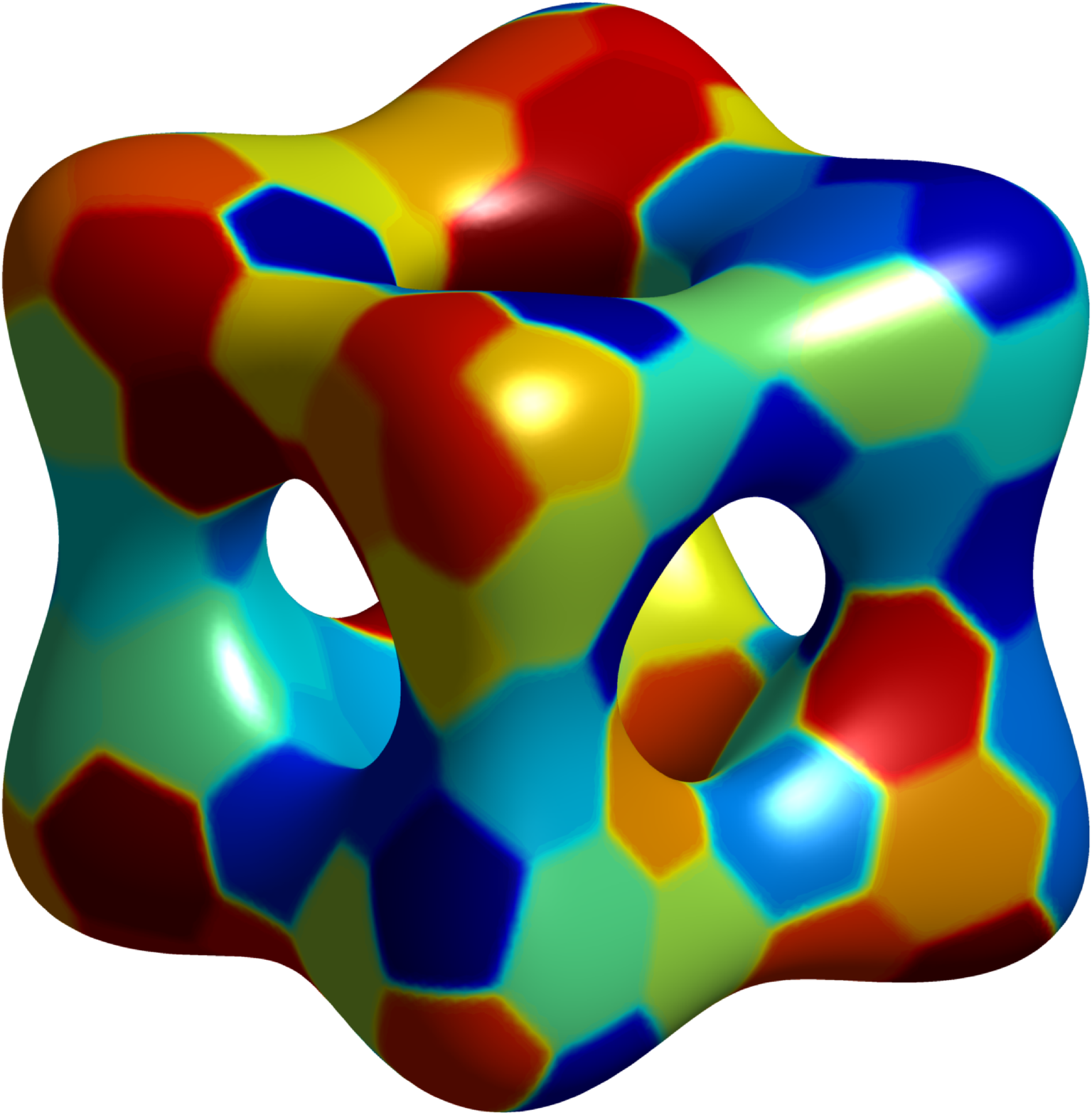}
\caption{Partitions on other types of surfaces}
\label{general_bel}
\end{figure}

On the other hand, we may ask what happens if we consider many cells cells in a partition of a surface. We may wonder if the analogue of the Caffarelli-Lin conjecture is also valid in the case of certain surfaces. The algorithm described in Section \ref{sect.num} allows us to study such partitions. In the case of the sphere computations were made for $n \leq 150$. Using Algorithm \ref{detect_top2D} we observe that for each $n\geq 14$ the numerical minimizer is composed of $n-12$ hexagons and $12$ pentagons. This is to be expected, due to Euler's formula $E+2=V+F$, where $E$ denotes the number of edges, $V$ the number of vertices and $F$ the number of components of the partition. This formula implies that there is no hexagonal partition of the sphere. Moreover, if only hexagons and pentagons are present, then we must have $n-12$ hexagons and $12$ pentagons. Some of the results in the case of the sphere are presented in Figure \ref{sphere_bel}.
Studying other types of surfaces reveals a similar behavior, as can be seen in Figure \ref{general_bel} for a torus of radii $R=1, r=0.6$, a Banchoff-Chmutov surface\footnote{defined as the zero level set of the function $T_4(x)+T_4(y)+T_4(z)$, where $T_4(x) = 8x^4-8x^2+1$} and a surface diffeomorphic to a ball\footnote{taken from \cite{elliott-ranner}, attributed to Dziuk, defined as the zero level set of $\Phi(x,y,z) = (x-z^2)^2+y^2+z^2-1$}. Hexagonal patterns tend to emerge with additional pentagons or heptagons. As already noted in \cite{elliott-ranner}, pentagons tend to appear in areas of high curvature and heptagons in areas with small curvature.

\subsection{Dimension Three}

As recalled in the introduction, previous works mostly relate to two dimensional partitions, since working in 3D requires heavy computational costs if we work on the whole grid. The grid restriction procedure described in Section \ref{sect.num} allows us to diminish the computational cost and work with resolutions up to $100\times 100\times 100$. We present below some numerical simulations concerning partitions of various three dimensional objects. 

\begin{figure}
\centering
\includegraphics[width = 0.19\textwidth]{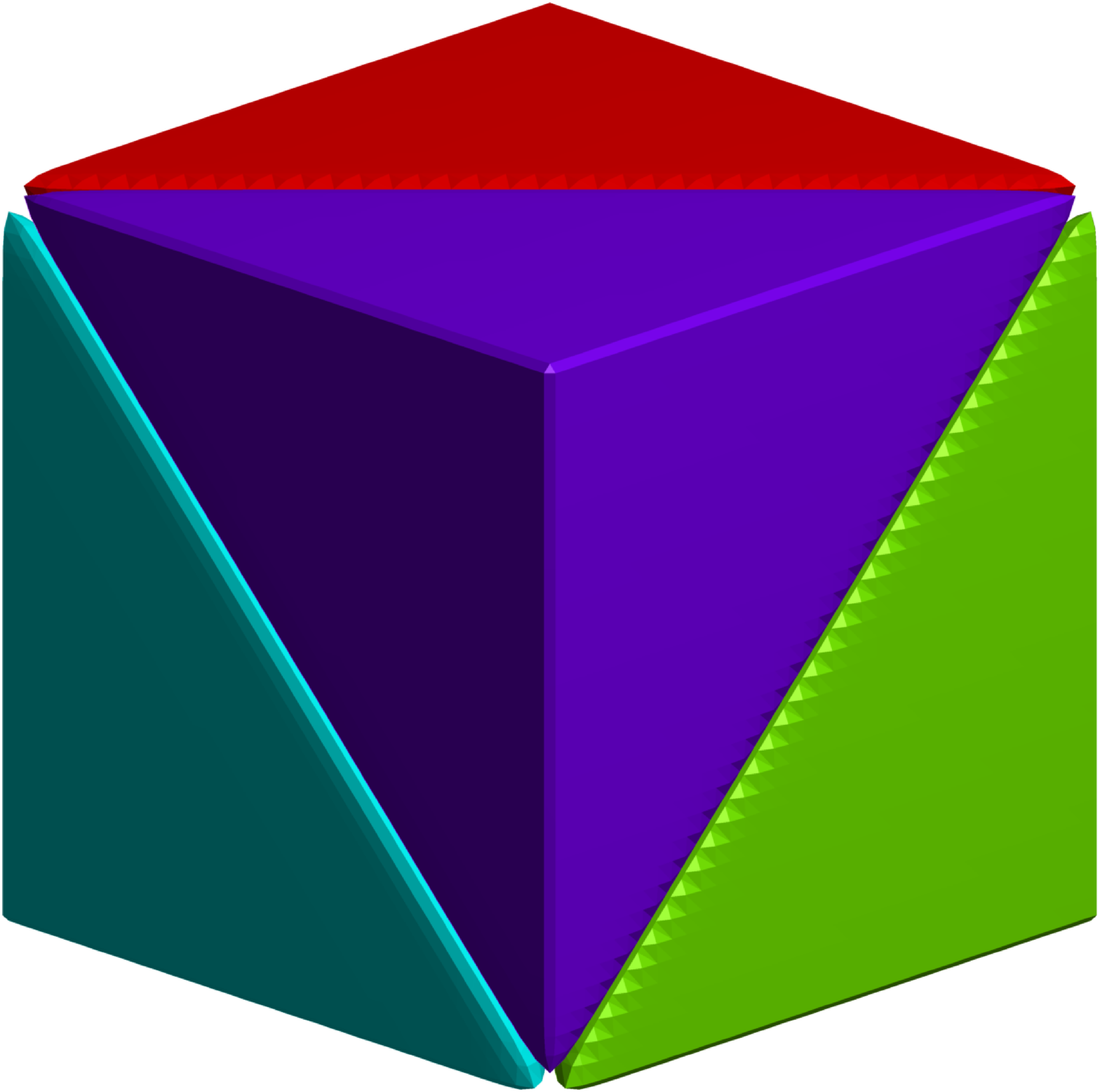}~
\includegraphics[width = 0.19\textwidth]{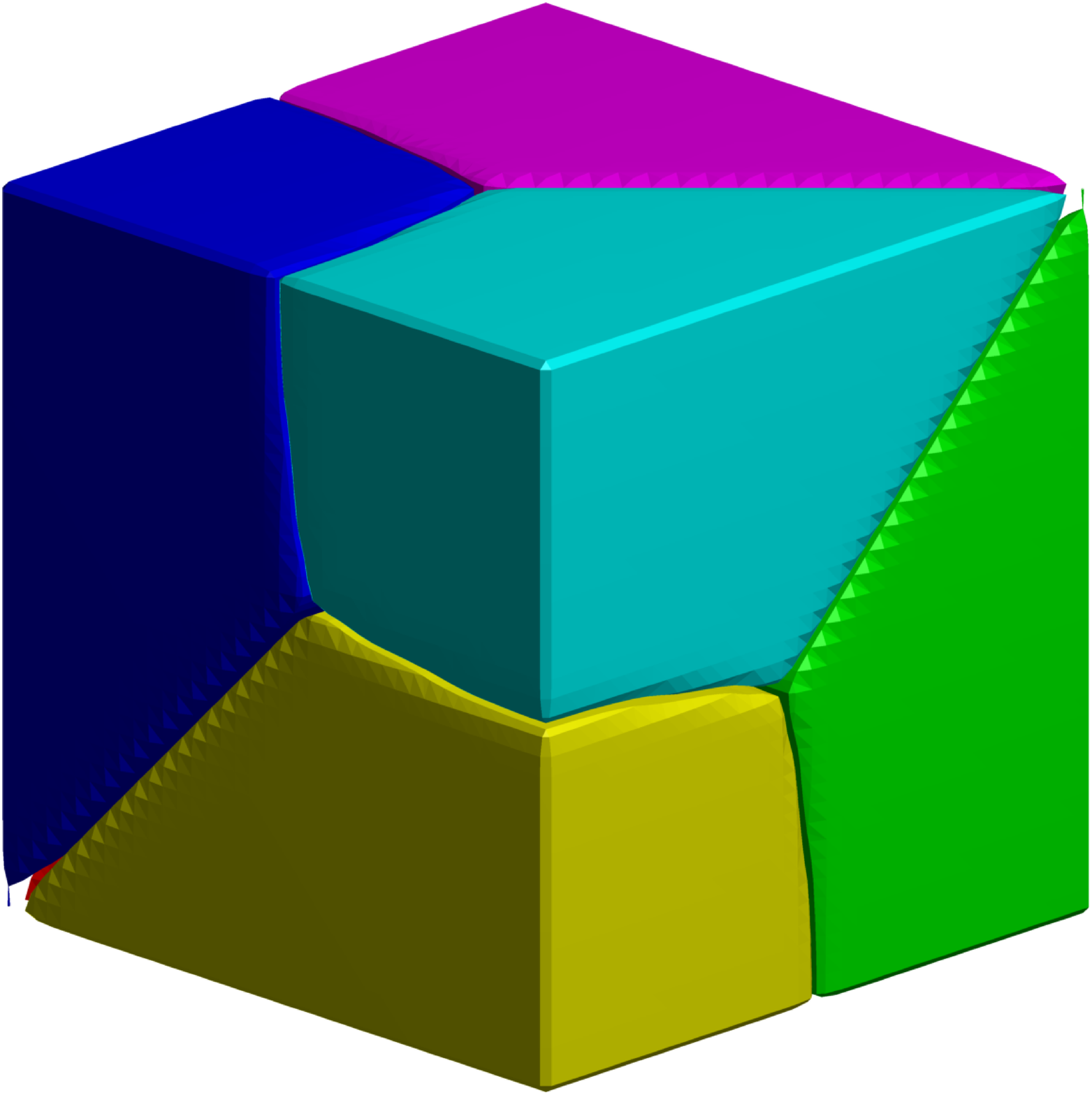}
\includegraphics[width = 0.19\textwidth]{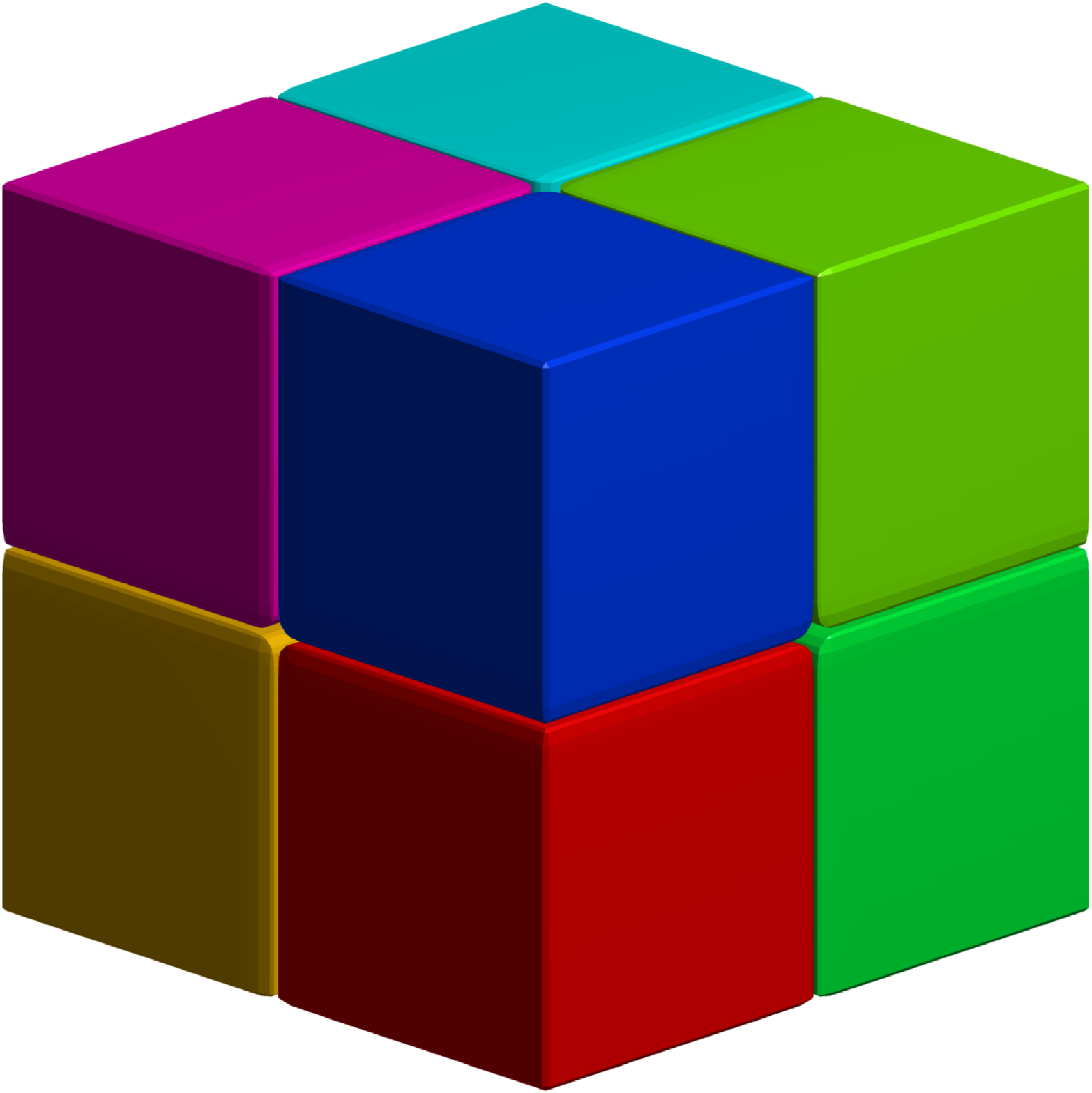}~
\includegraphics[width = 0.19\textwidth]{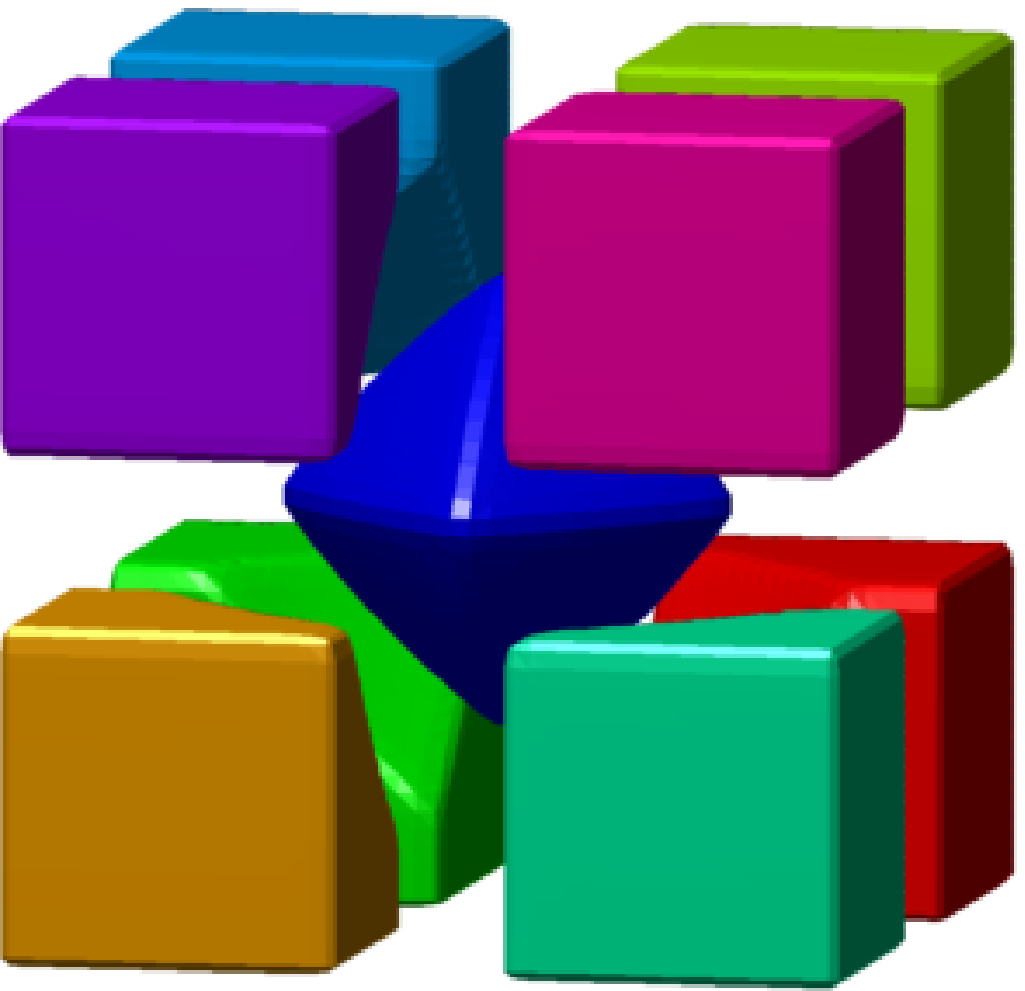}~
\includegraphics[width = 0.2\textwidth]{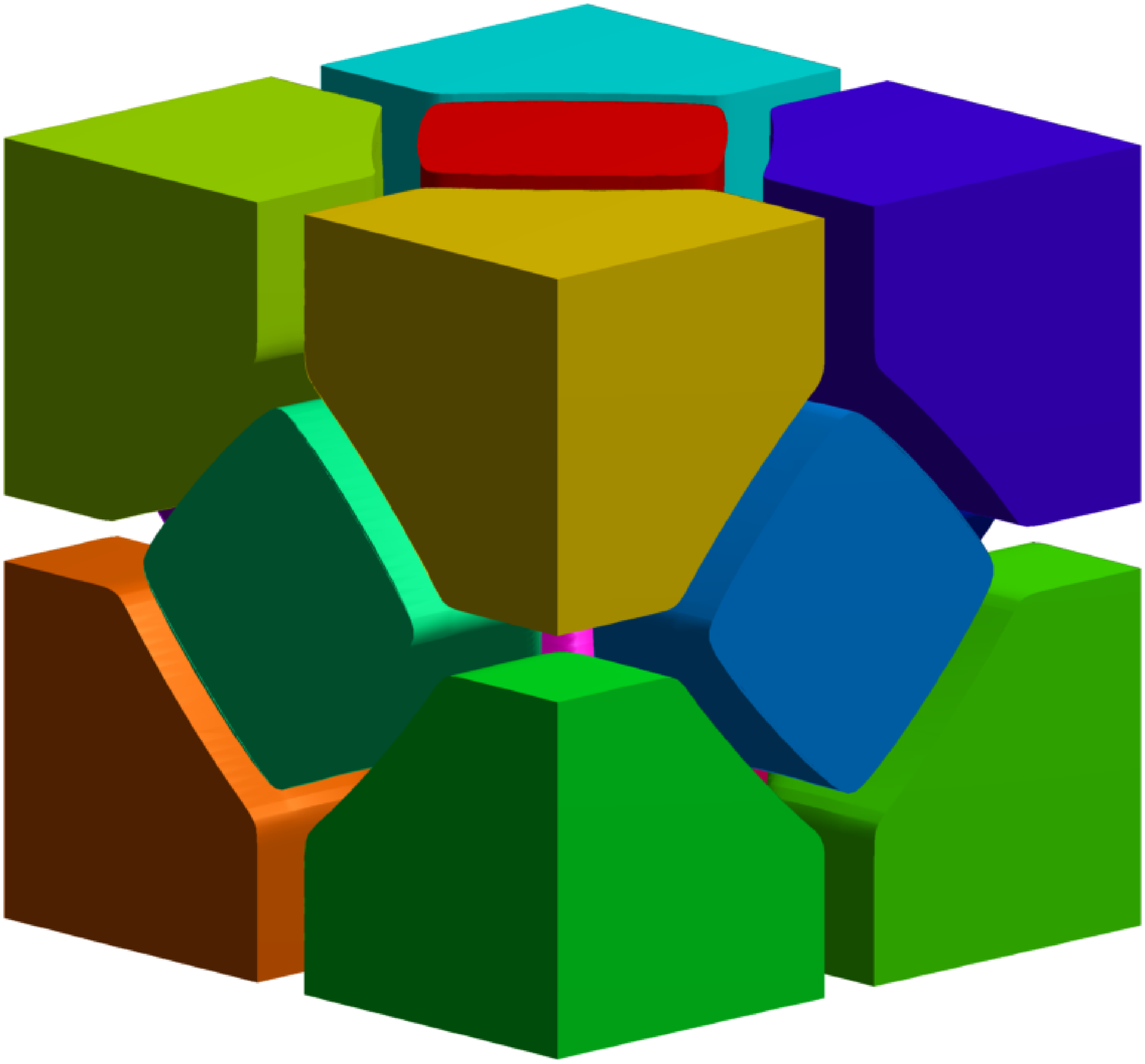}

\caption{Some results for the case of the 3D partition of the cube for $n=4,6,8,9,14$}
\label{3D_cube}
\end{figure}
\begin{figure}
\centering
\includegraphics[width = 0.19\textwidth]{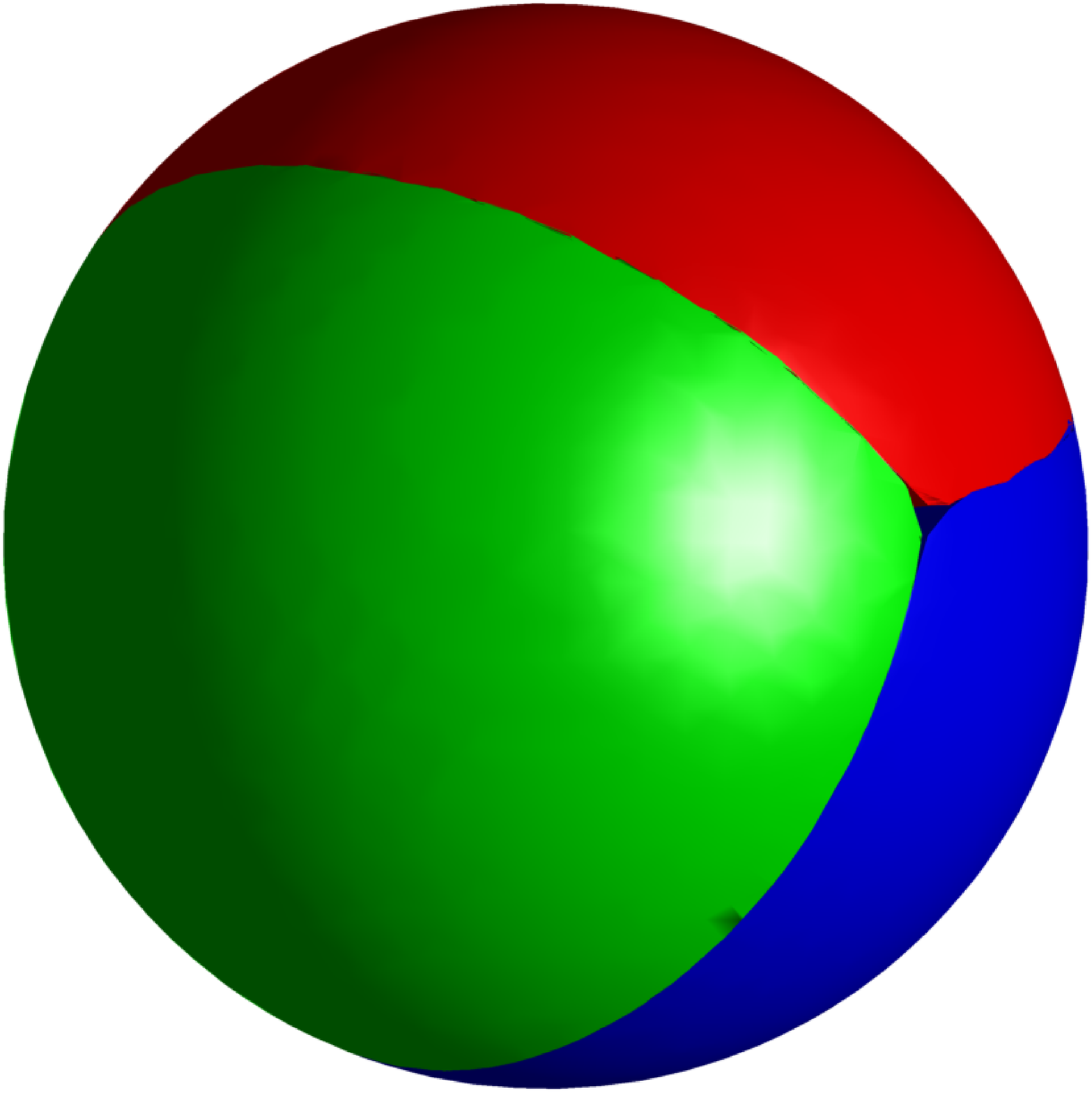}~
\includegraphics[width = 0.19\textwidth]{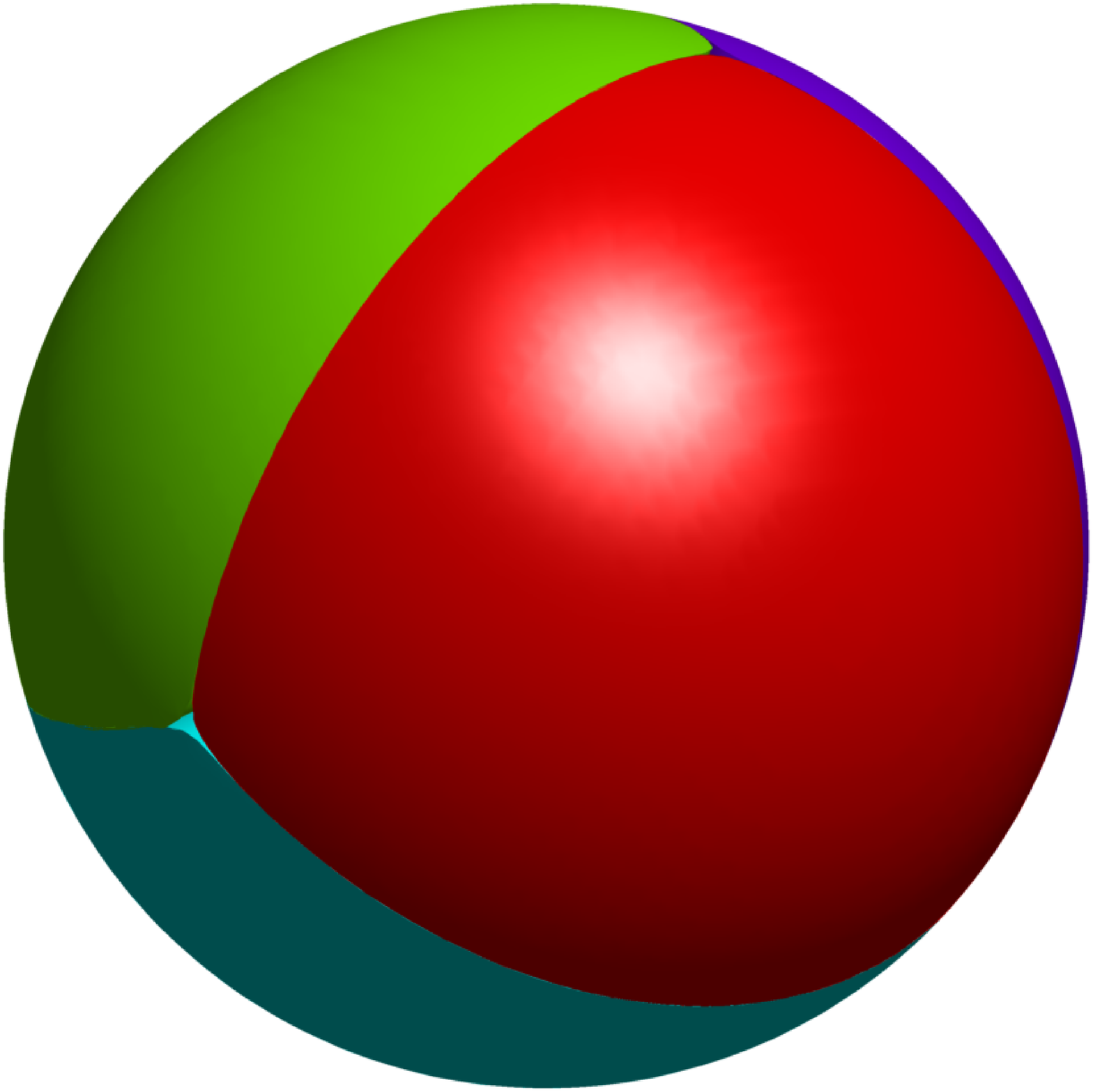}~
\includegraphics[width = 0.19\textwidth]{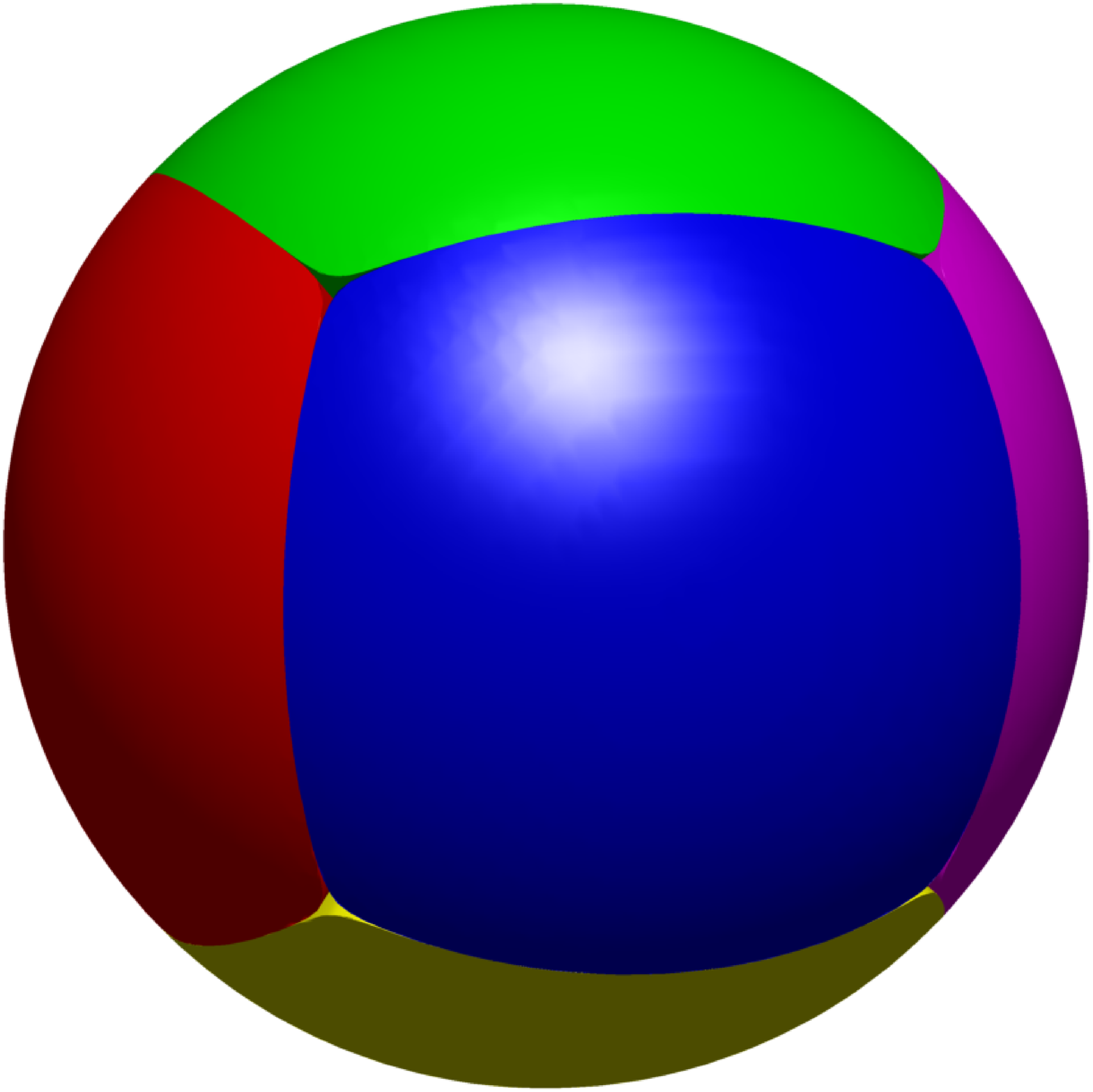}~
\includegraphics[width = 0.19\textwidth]{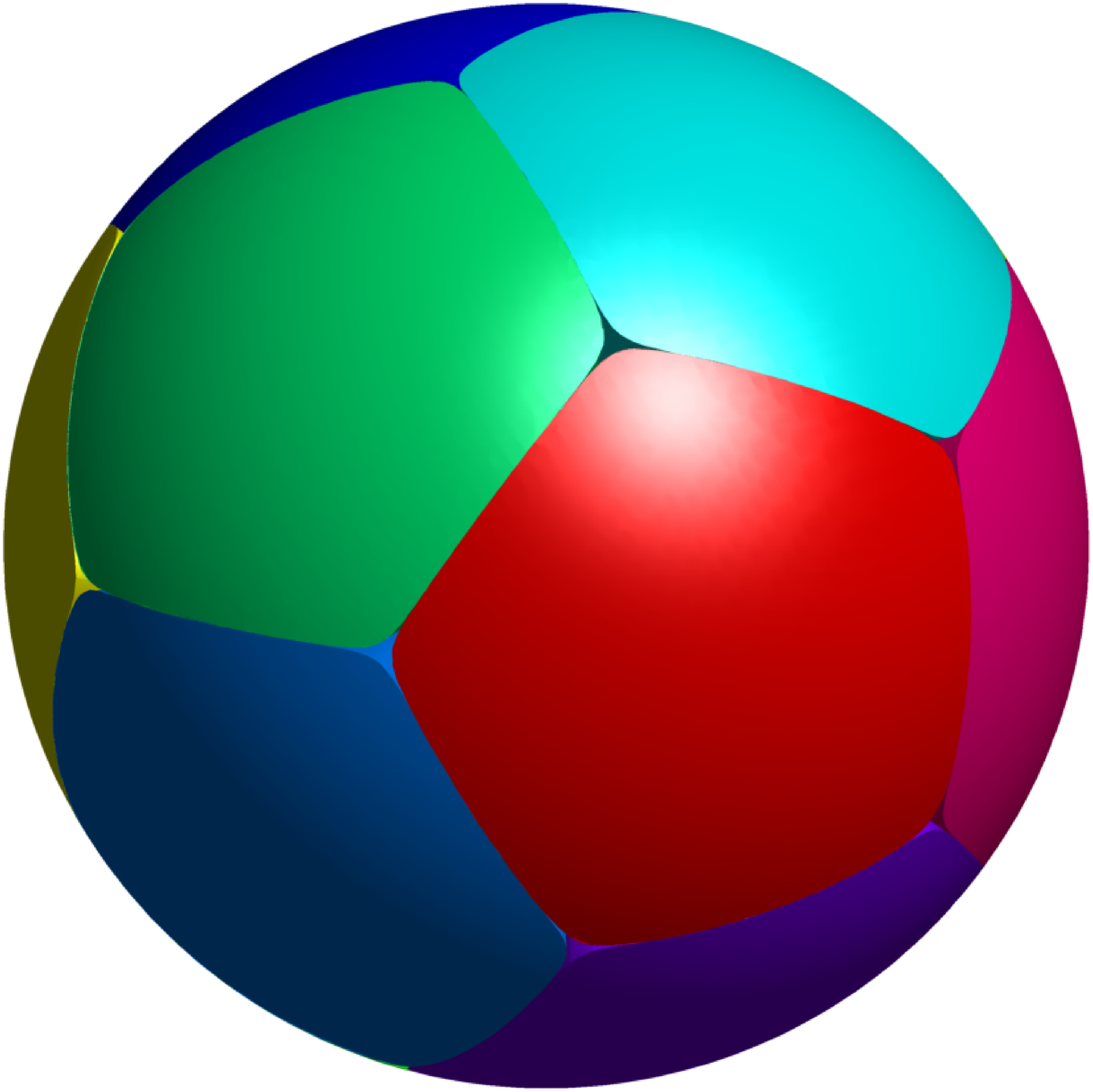}~
\includegraphics[width = 0.19\textwidth]{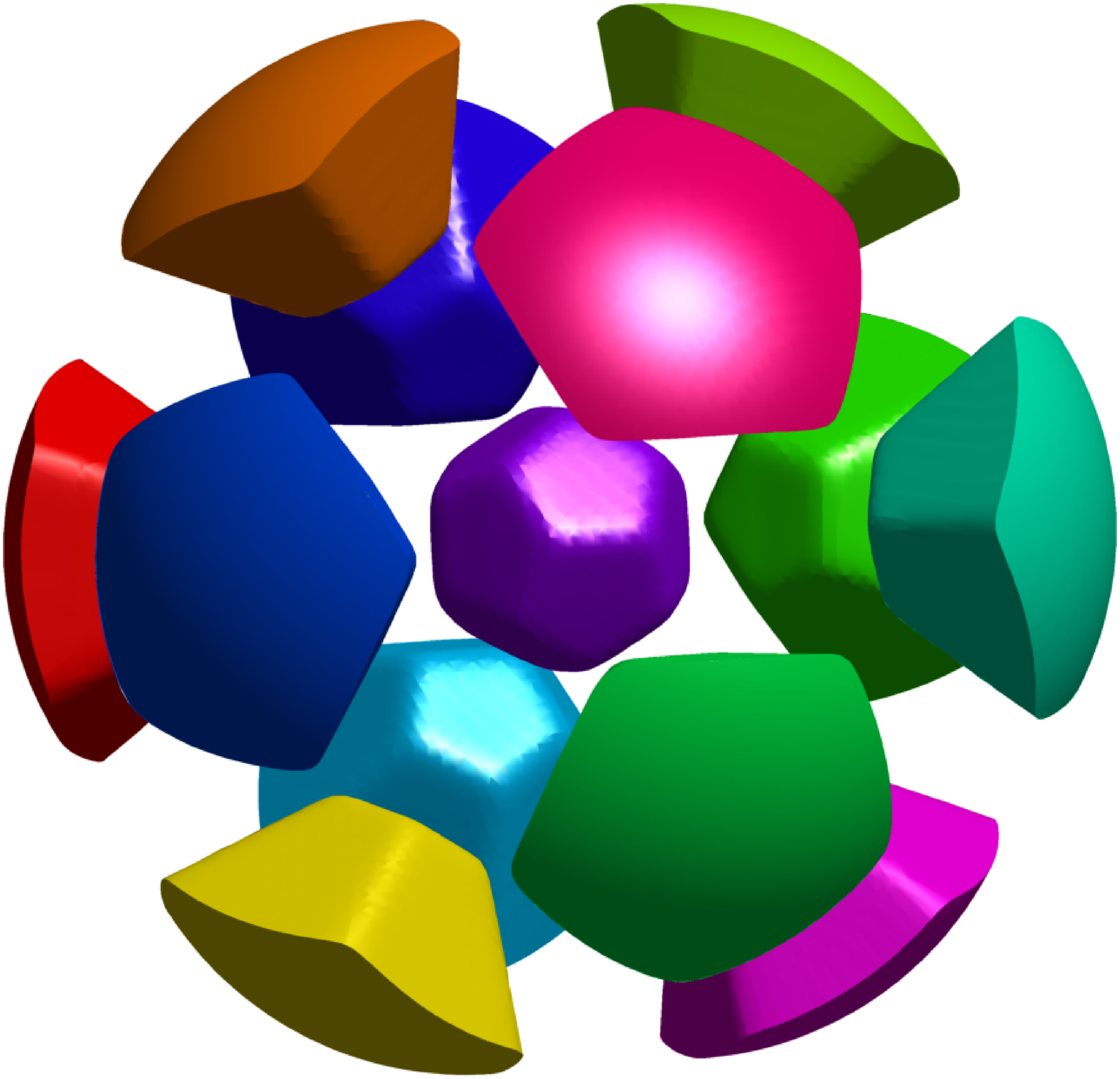}
\caption{Some numerical optimal partitions for the ball.  The regular cases: $n=3,4,6,12$. The last image is for $n=13$, the first instance when we have a cell in the interior of the partition.}
\label{3D_ball}
\end{figure}
\begin{figure}
\centering
\includegraphics[width = 0.19\textwidth]{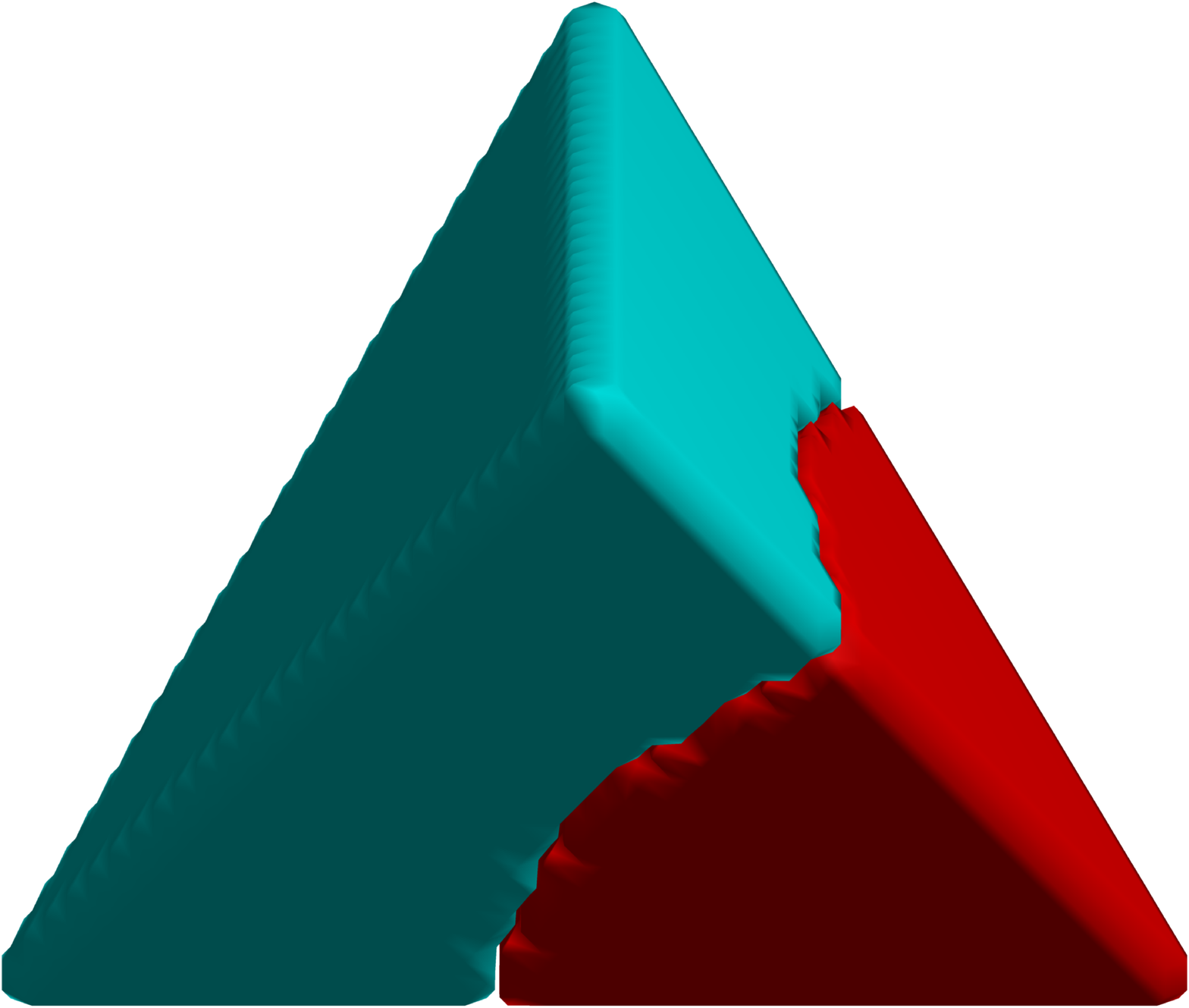}~
\includegraphics[width = 0.19\textwidth]{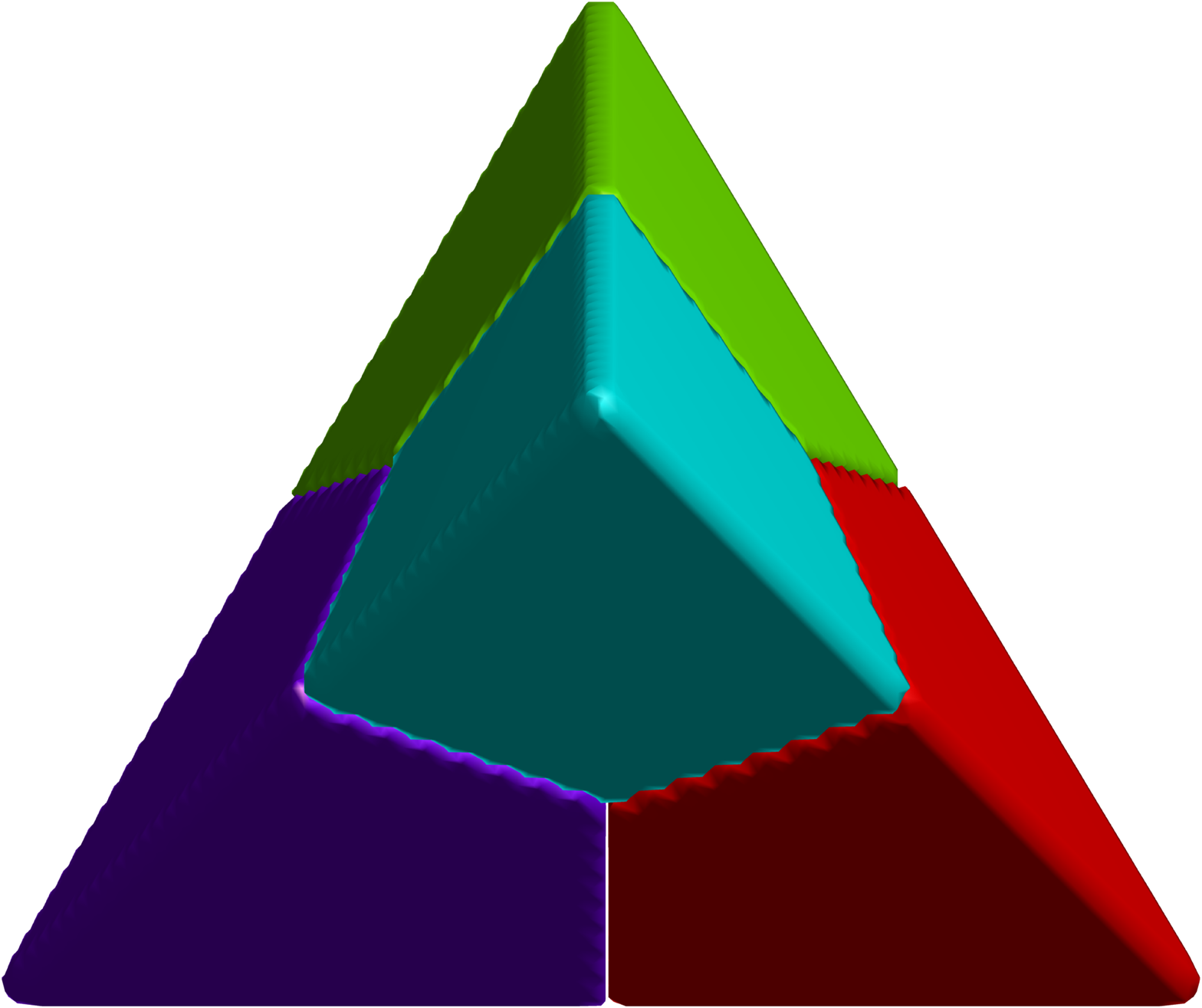}~
\includegraphics[width = 0.19\textwidth]{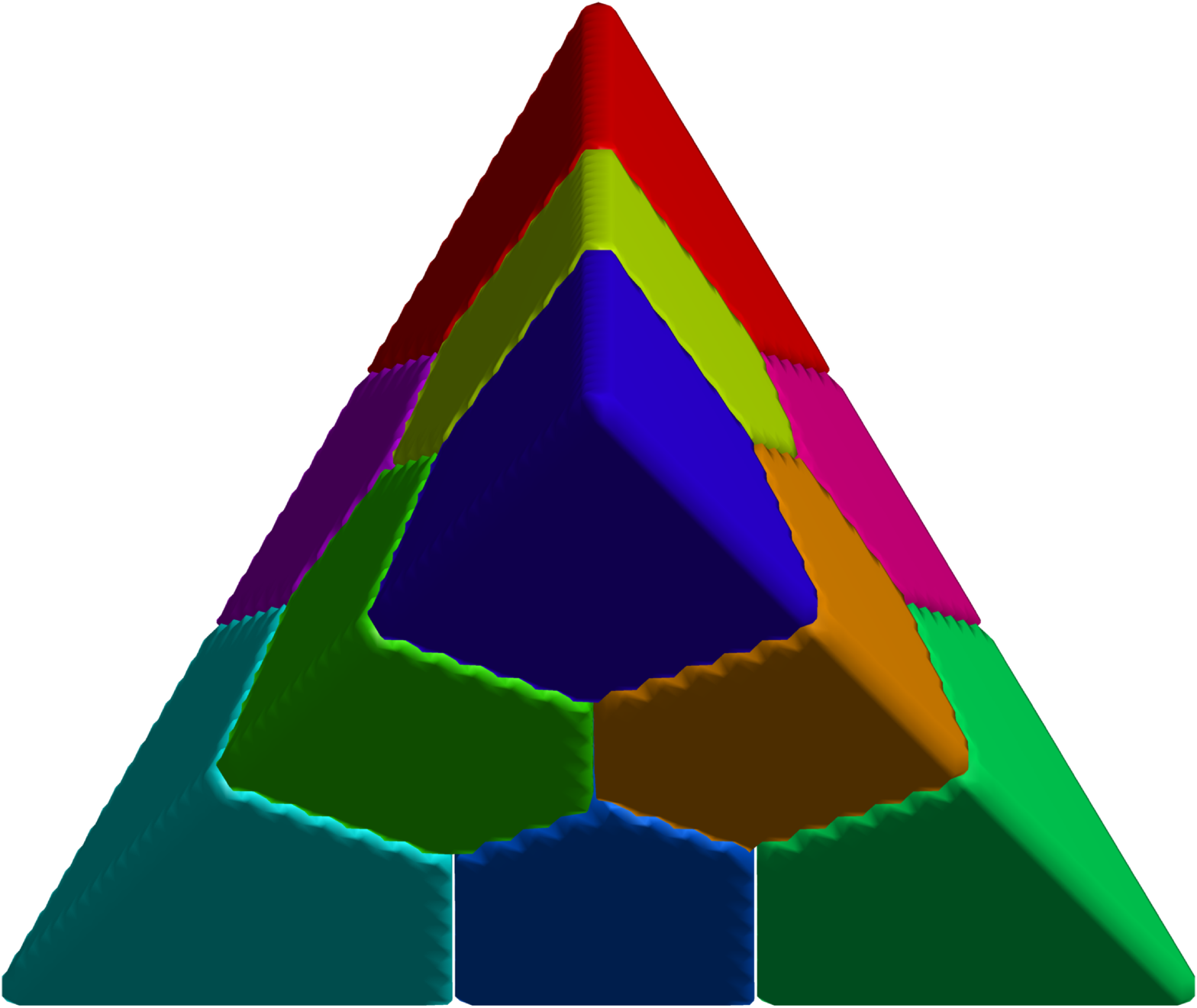}~
\includegraphics[width = 0.19\textwidth]{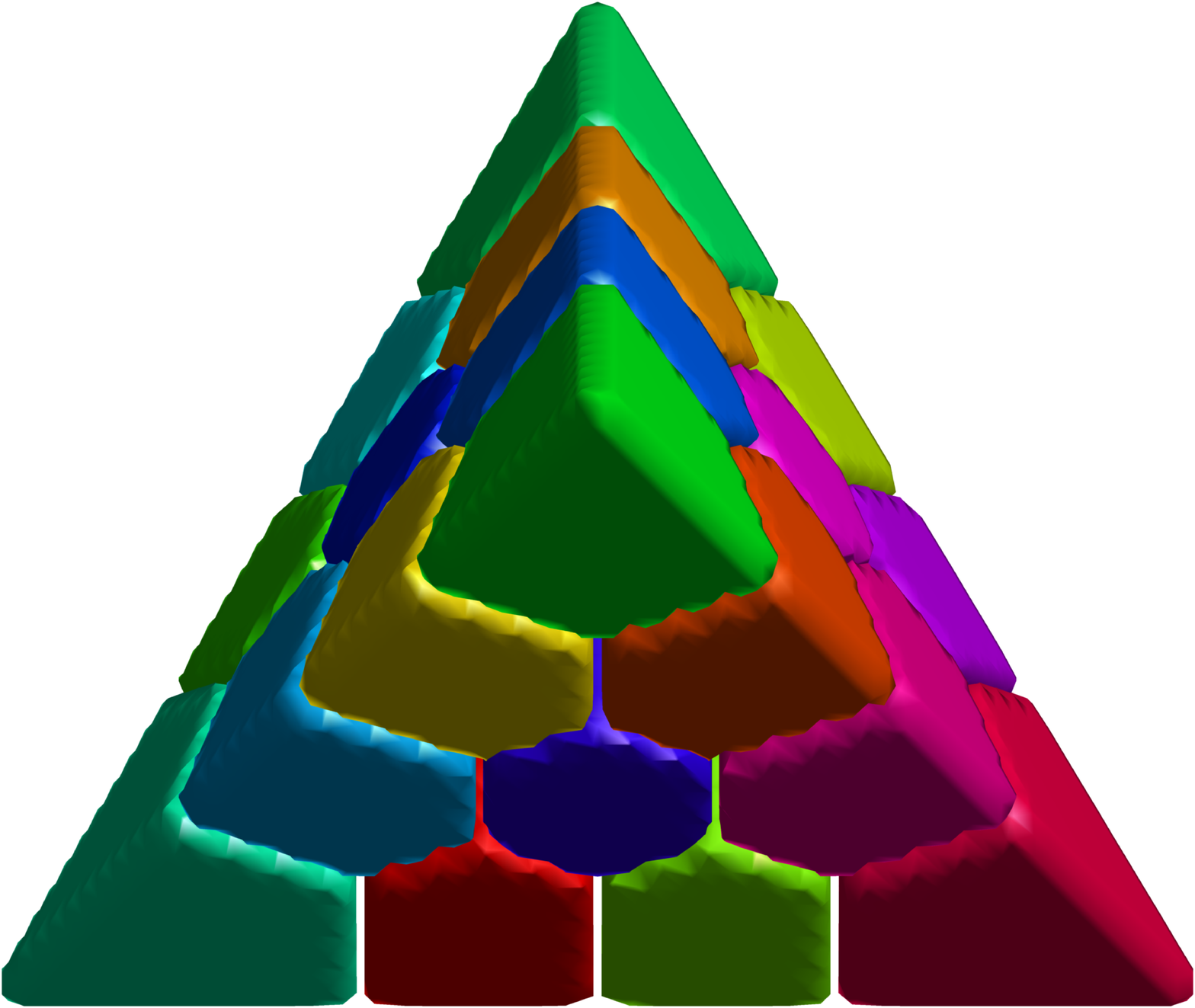}~
\includegraphics[width=0.19\textwidth]{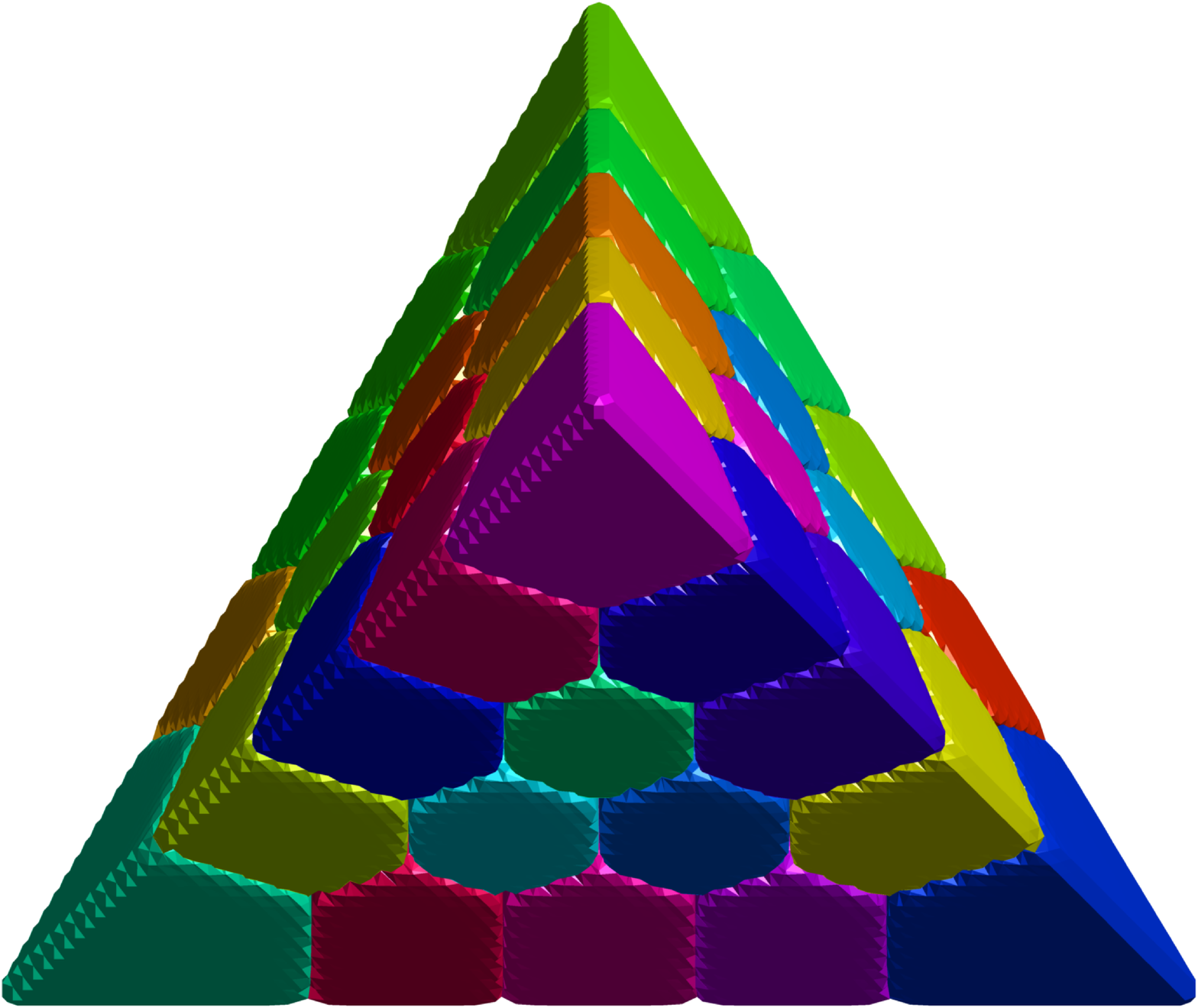}
\caption{Some numerical optimal partitions for the regular tetrahedron for $n=2,4,10,20,35$.}
\label{3D_tetra}
\end{figure}
\begin{figure}
\centering
\includegraphics[width = 0.3\textwidth]{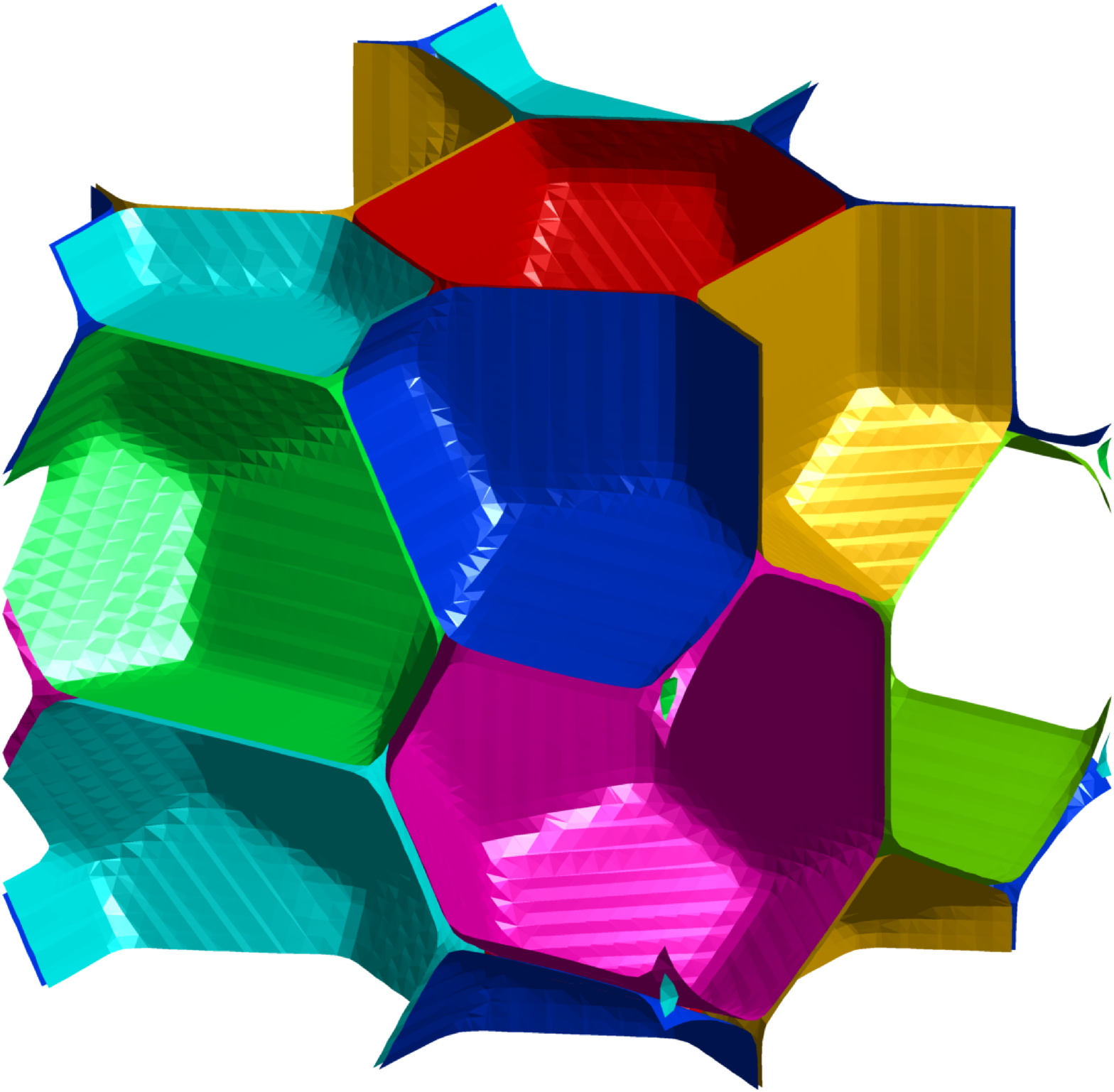}
\includegraphics[width=0.3\textwidth]{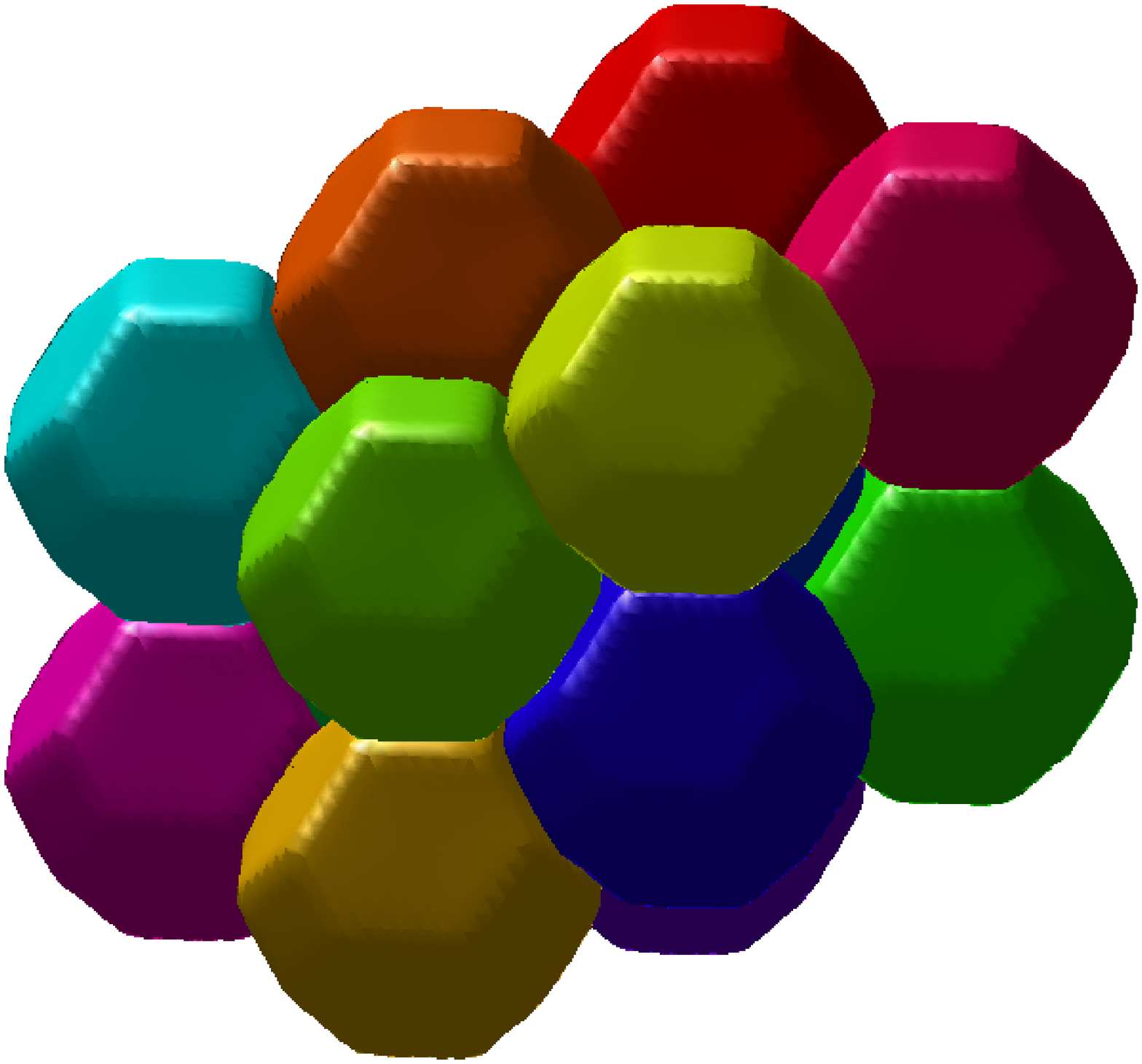}~
\includegraphics[width = 0.3\textwidth]{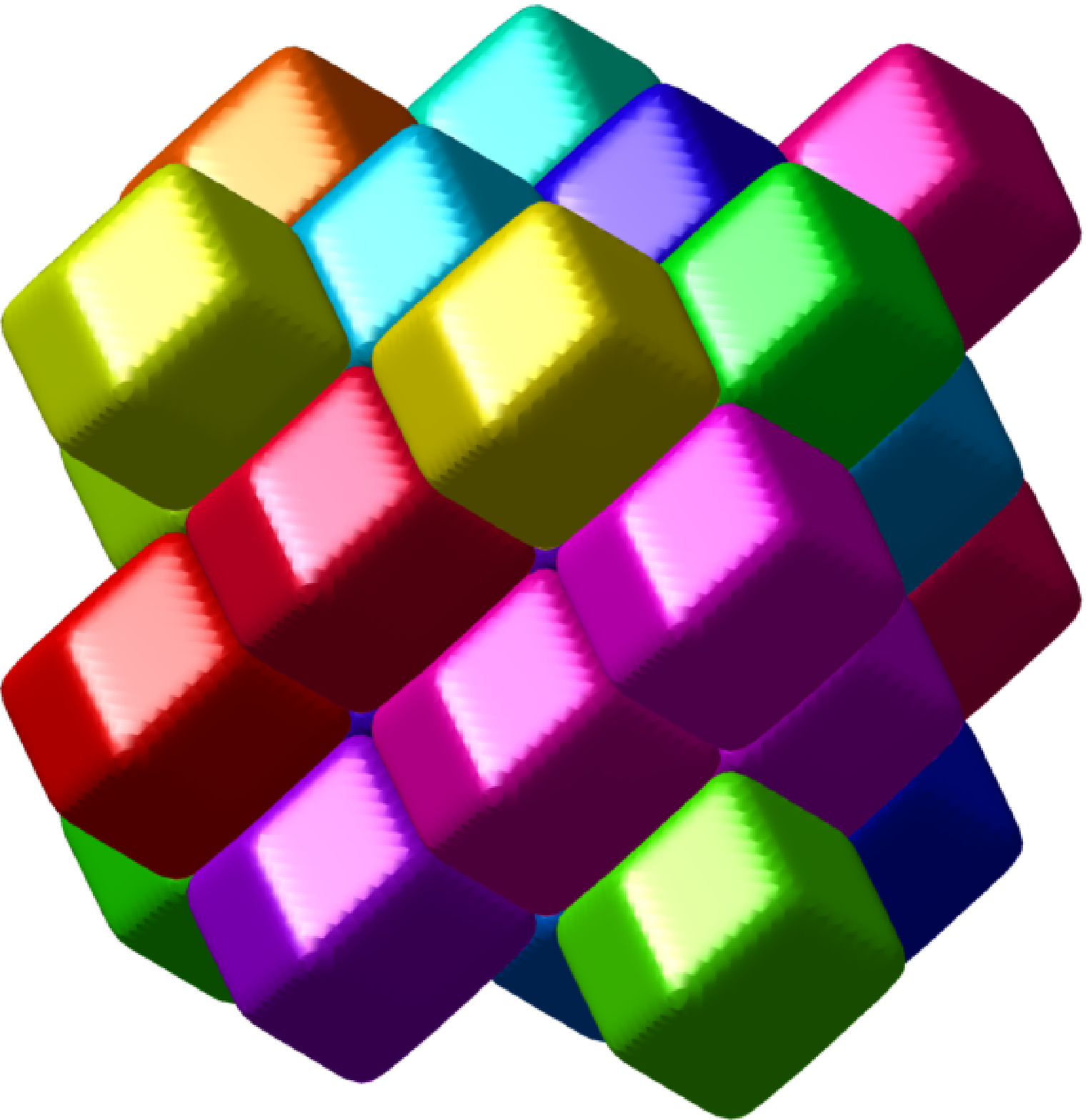}~
\caption{Some computations on the periodic cube: the Weaire-Phelan structure for $n=8$, the Kelvin structure for $n=16$ and the rhombic dodecahedron structure for $n=32$.}
\label{3D_percube}
\end{figure}
Optimal candidates in the case of the cube for various values of $n$ are shown in Figure \ref{3D_cube}. Some computations for the case of the ball are presented in Figure \ref{3D_ball}. We also consider the case of a regular tetrahedron in Figure \ref{3D_tetra}. In order to simulate partitions of $\Bbb{R}^3$ we may consider periodic conditions for a cube. Some results in the periodic case are presented in Figure \ref{3D_percube}.

Concerning partitions in a small number of cells, we make some global remarks. For $2$ cells, it seems that for the cube and the ball, the partitions which minimize the sum of eigenvalues are formed of congruent cells. As already noted in the two dimensional case (see \cite{BoBN16}), working on symmetric domains does not guarantee that optimal partitions with only 2 cells will be made of equal cells. We observe the same thing in 3D, when working on the regular tetrahedron. Like on the equilateral triangle (see \cite{BoBN16}, for example) we see that one of the two cells is near a corner as can be seen in Figure \ref{3D_tetra}. 

If we consider four cells we observe that the cube, the ball and the regular tetrahedron all seem to have optimal candidates made out of equal cells. This is not too surprising, since the regular tetrahedron has the right symmetries in order to decomposed into four congruent, symmetric pieces. On the other hand, the cube can be obtained from the regular tetrahedron by adding four regular pyramids, one on each face. 

In the case of the cube, we note that for $6$ cells we obtain again a type of partition made out of equal cells. For $9$ cells we obtain a first partition which has a cell which does not touch the boundary. Having in mind the results in \cite{HK15}, we know that the partition into $8$ cubes is not minimal for the max, since the $8$th eigenvalue of the cube is not Courant-sharp. For more details about the terminology regarding minimal partitions for the largest eigenvalue you may look in the article cited above or in \cite{BoBN16}. Since the partition in $8$ cubes is not optimal for the max, it cannot be optimal for the sum. This is due to the fact that the optimal partition for the max gives an upper bound for the one minimizing the sum. In our computations for $k=8$ we do not manage to obtain something better than the partition into $8$ cubes. This is obviously a good candidate to be a local minimizer. It is possible that in order to obtain a better candidate we would need to have a more fine computational grid.

\begin{figure}
\centering
\includegraphics[width=0.2\textwidth]{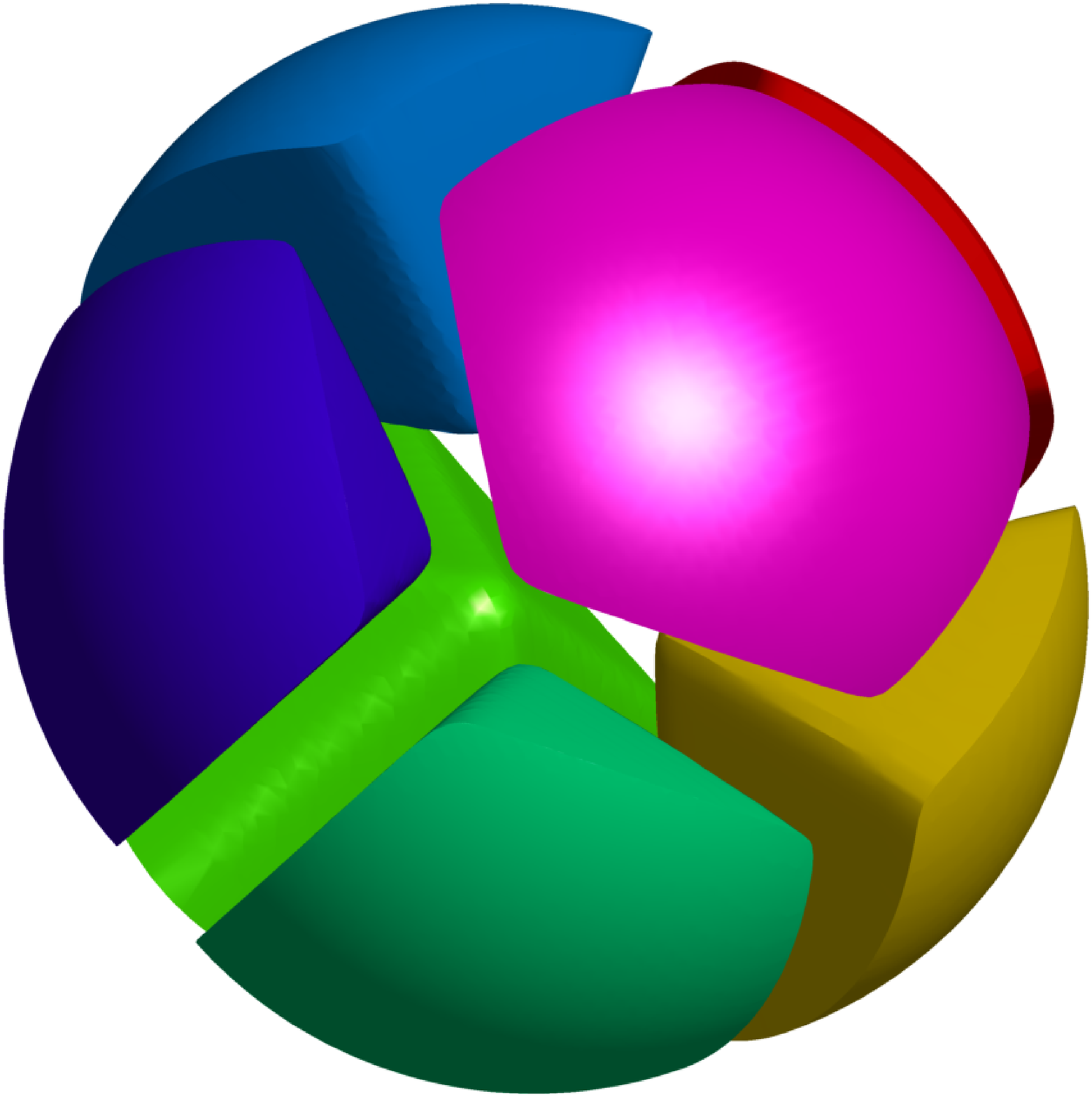}~
\includegraphics[width=0.2\textwidth]{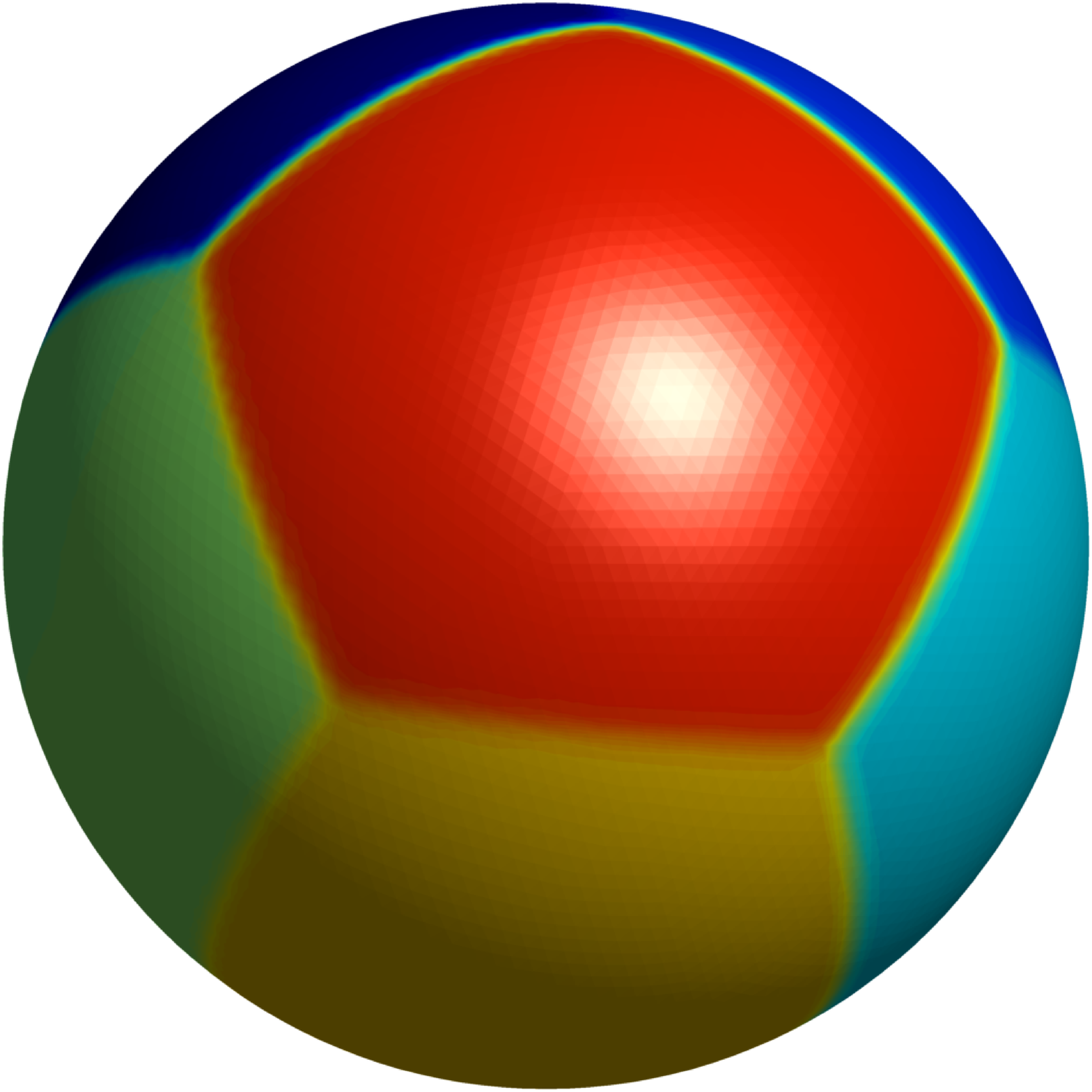}~
\includegraphics[width=0.2\textwidth]{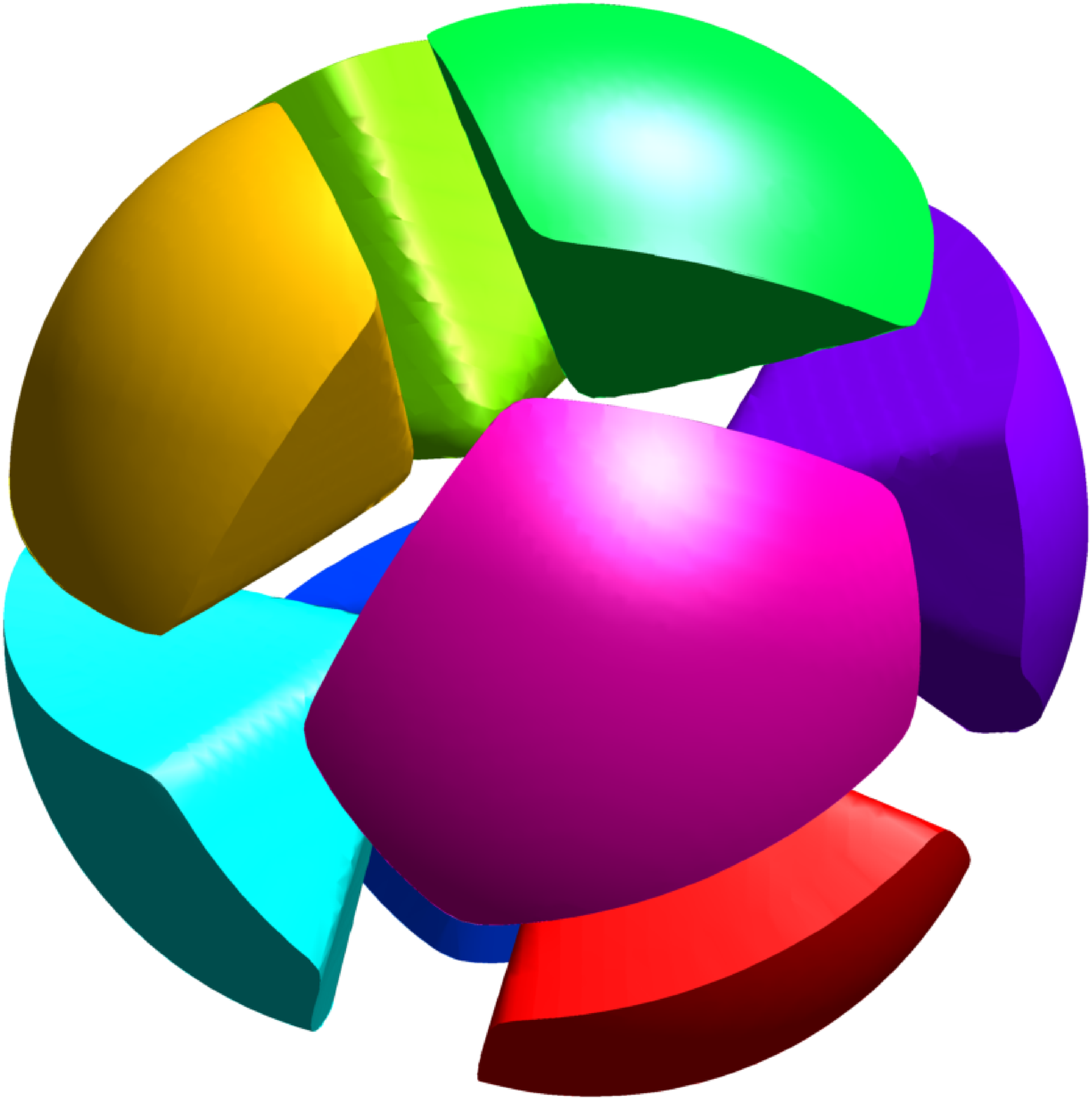}~
\includegraphics[width=0.2\textwidth]{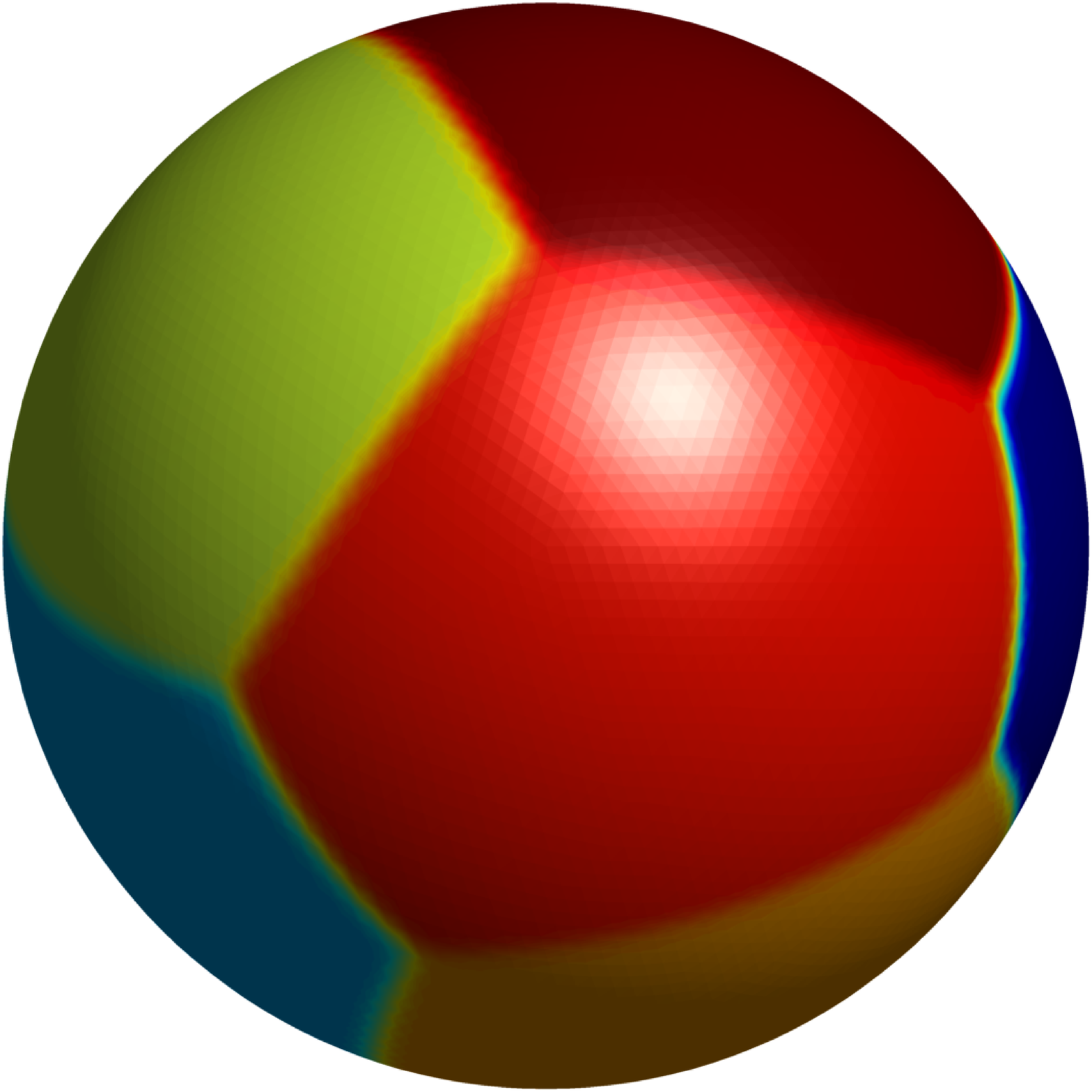}~
\caption{Presentation of optimal partitions in the ball and on the sphere for $n=7,8$}
\label{ball_bel}
\end{figure}

For the ball we obtain that for $n=2,3$ the optimal candidates are partitions into two half balls and three equal slices of angle $2\pi/3$. Moreover, symmetric partitions corresponding to regular polyhedra are obtained for $n=4,6,12$. For all $n \leq 12$ we observe an interesting connection between the partitions obtained on the ball and those obtained for the Laplace-Beltrami operator in the unit sphere. It seems that the cells of the optimal candidates for the ball are \emph{cones} determined by sets on the sphere. More precisely, the sets observed in partitions of the ball with less than $12$ cells are of the form $\omega \times [0,1] = \{ \lambda x : x \in \omega, \lambda \in [0,1]\}$ where $\omega$ is a subset of $\Bbb{S}^2$. Moreover, in these cases the structure of the partitions on the ball resembles the one on the sphere. For some examples see Figure \ref{ball_bel} (As a parallel, we recall that in 2D the situation is similar on the disk. For $n \leq 5$ optimal partitions seem to be made out of equal sectors, which correspond to optimal partitions of the Laplace-Beltrami on the circle. Indeed, noting that the one dimensional Laplace-Beltrami eigenvalues are explicit in terms of the length of the curves, it is possible to prove that the partition in to arcs of circles of equal length is optimal on the circle.)

\begin{figure}
\centering
\includegraphics[width=0.2\textwidth]{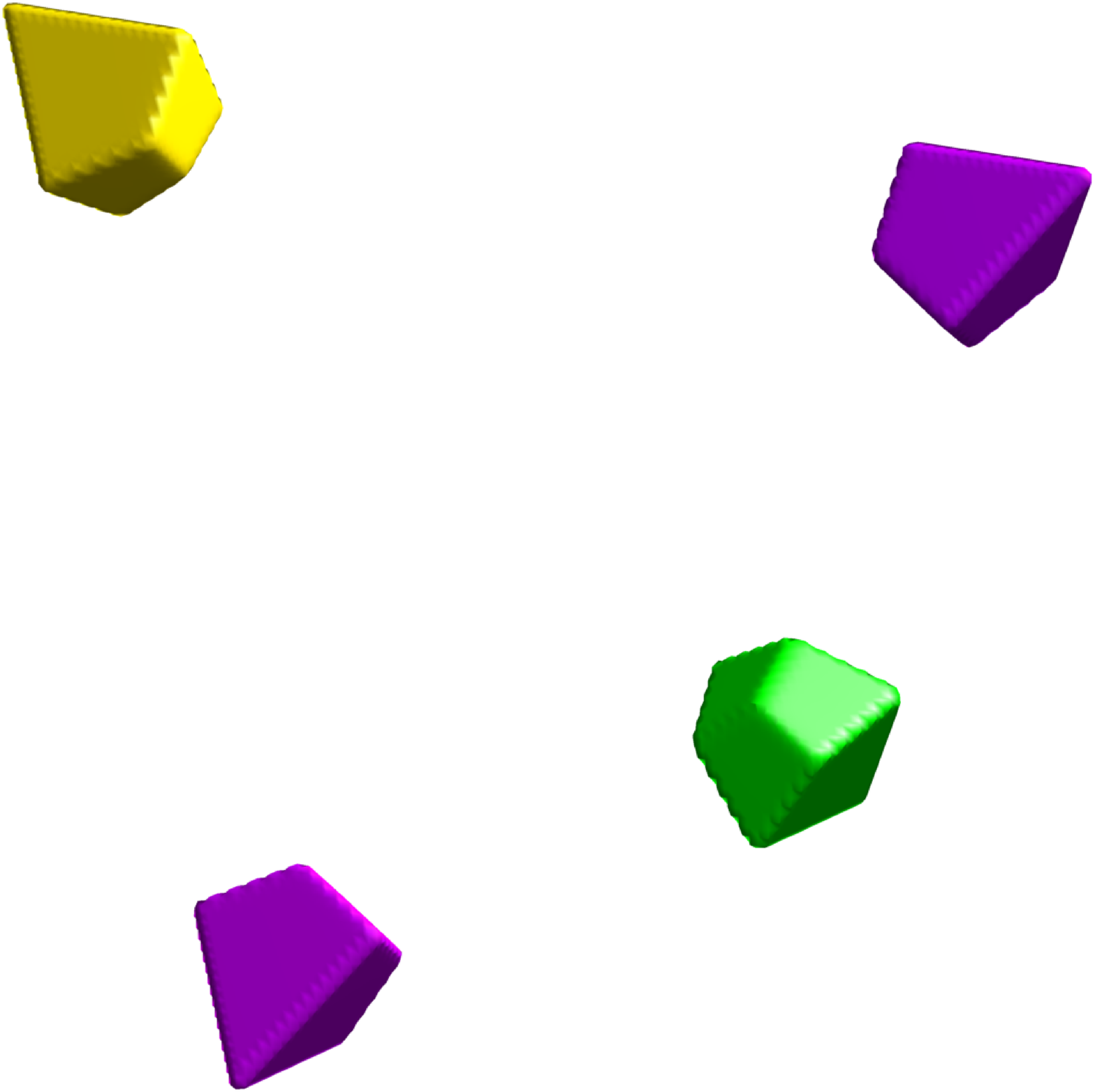}~
\includegraphics[width=0.2\textwidth]{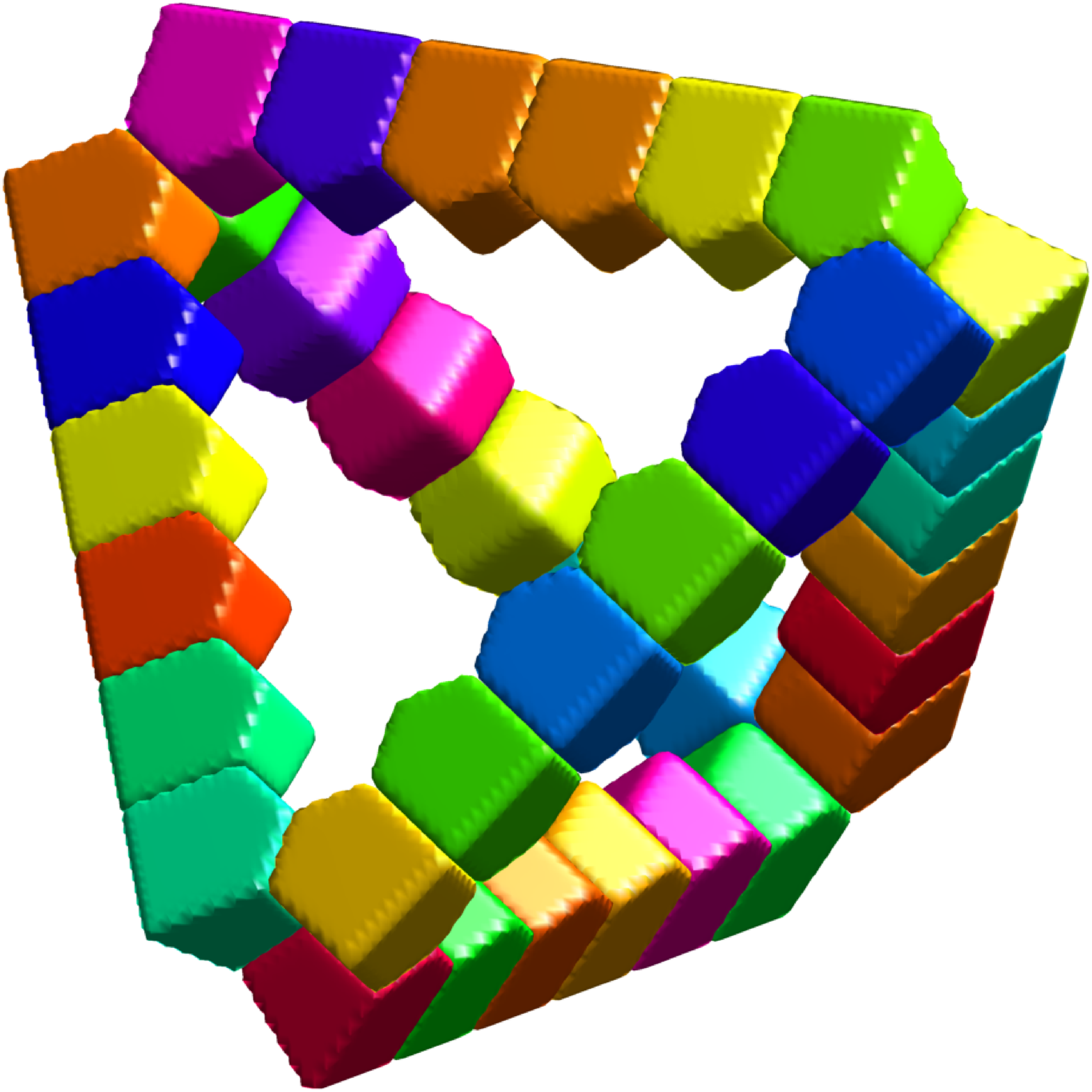}~
\includegraphics[width=0.2\textwidth]{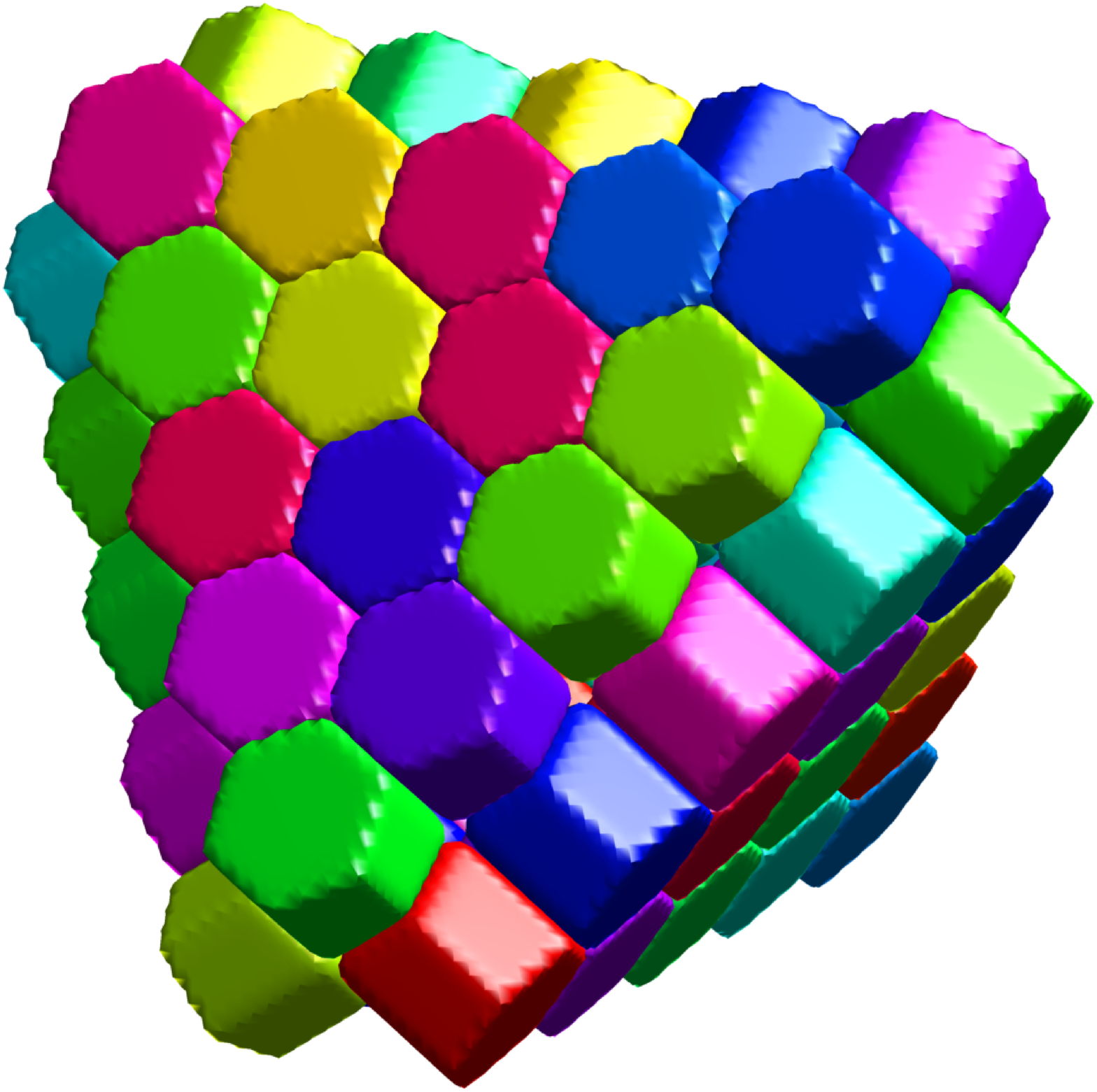}~
\includegraphics[width=0.2\textwidth]{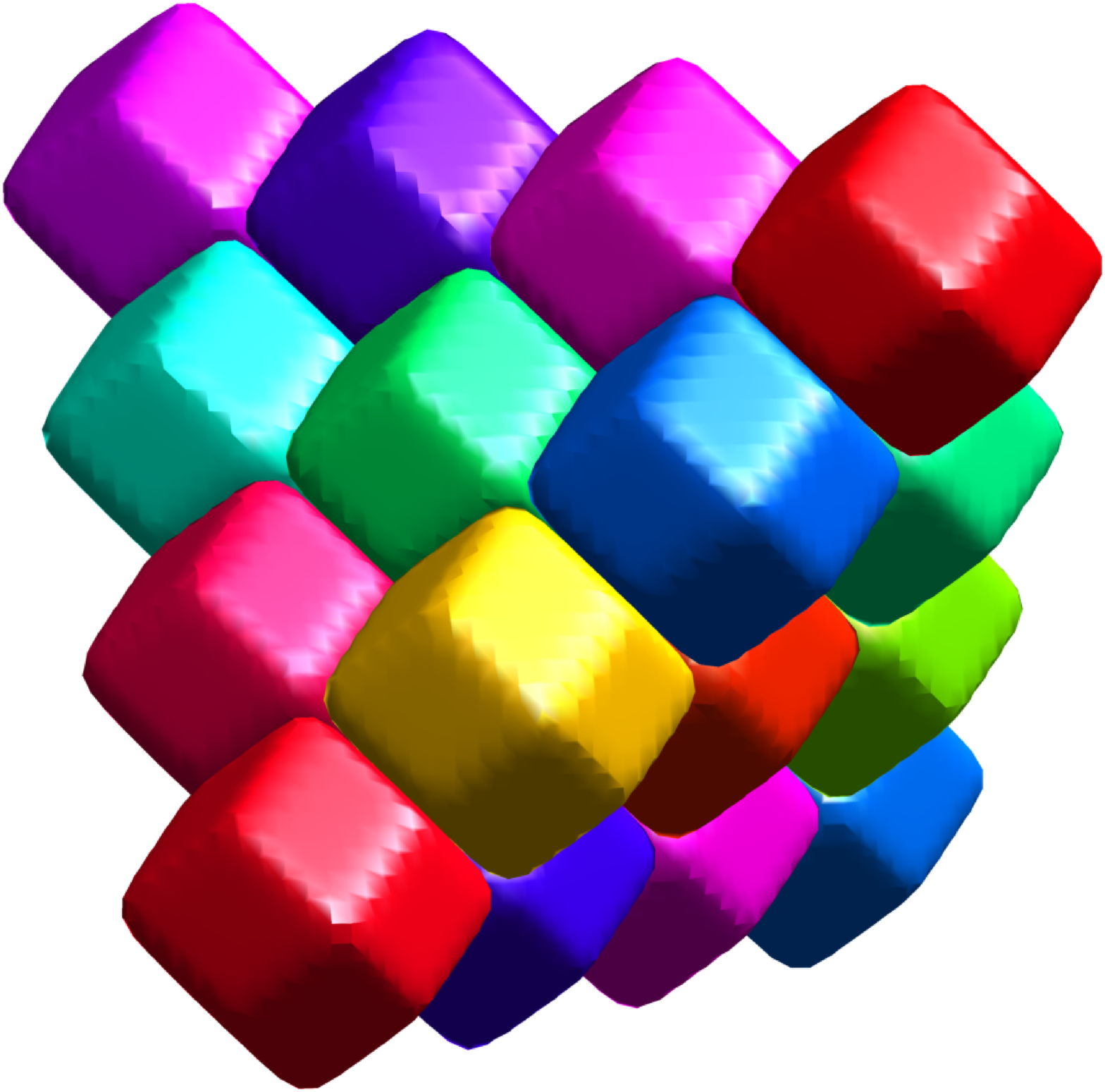}
\caption{Automatic classification using the spectrum of the Laplace-Beltrami operator for the cells of the partition of the tetrahedron into $120$ cells: $4$ corner cells, $36$ cells along the edges, $60$ cells corresponding to the faces, $20$ cells in the interior.}
\label{class_tetra120}
\end{figure}

In the case of the regular tetrahedron, as in the case of the equilateral triangle in the plane, there are some numbers of cells which give rise to interesting structures. We consider numbers which are sums of the first $k$ triangular numbers $T_i = i(i+1)/2$. These are called \emph{pyramidal numbers} and we observe that they give rise to partitions which seem polyhedral. There is an explicit formula for such numbers, given by $P_n = n(n+1)(n+2)/6$. The first terms in this sequence are $1,4,10,20,35,56,84,120$. We observe that partitions of the regular tetrahedron into $P_n$ cells seem to be polyhedral. Moreover, the traces of these partitions on the sides of the regular tetrahedron have the same structure as the partition of the equilateral triangle into $T_n$ cells. For pyramidal numbers greater than $35$ there are also cells which are in the interior of the regular tetrahedron. A closer inspection of these cells shows that they seem to be \emph{rhombic dodecahedra}. This case provides an excellent test case for Algorithm \ref{detect_top3D}. Running the algorithm on the cells of the partition of the tetrahedron for $n=120$ provides the results given in Figure \ref{class_tetra120}. We observe that the algorithm detects correctly cells which are similar.

We may wonder what is the 3D equivalent of the spectral honeycomb problem in the plane. In order to simulate partitions of $\Bbb{R}^3$ we may work on the periodic cube. For the $1\times 1 \times 1$ cube and $n=16$ we obtain the Kelvin structure where all cells are equal and resemble a truncated octahedron. For $n=8$ we obtain the Weaire-Phelan structure consisting of pieces of two types: $6$ identical polyhedral cells with faces consisting of $12$ pentagons and $2$ hexagons and $2$ irregular dodecahedra having $12$ pentagons as faces. Running the Algorithm \ref{detect_top3D} on cells of Figure \ref{3D_percube} allows us to automatically see that for $n=16$ we have $16$ similar cells, while for $n=8$ we obtain $6$ cells of one type and $2$ cells of another type, corresponding to the Weaire-Phelan structure \cite{weaire-phelan}. Note that the simulations performed in the case of the regular tetrahedron presented above and shown in Figure \ref{class_tetra120} give another candidate, which is the partition into rhombic dodecahedra. Indeed, numerical results on the periodic cube for $n=4$ and $n=32$ consist of tilings with rhombic dodecahedra. For $n=12$ and $n=20$ we also obtain tilings of the cube with congruent cells. In each of these cases Algorithm \ref{detect_top3D} was used in order to automatically classify the cells of the partition. For $n \in \{12,16,20,32\}$ the algorithm concluded that we have only one class of congruent cells, while for $n=8$ we have two classes, one containing $6$ cells and the other the remaining $2$ cells. 

It seems that simulations in the periodic case do not give enough information about the three dimensional optimal partition in the asymptotic case. The results obtained depend on the number of cells and on the sizes of the rectangular box considered in the computation. Nevertheless, the fact that the rhombic dodecahedra partition can be obtained as numerical optimizer in the periodic case and that this structure seems to appear in the optimal partition of the regular tetrahedron for pyramidal numbers of cells makes it a good candidate. Moreover, the rhombic dodecahedron partition is the Voronoi diagram associated to an optimal sphere packing in 3D, as shown in \cite{hales_kepler}. If we make again the connection between the sphere packing problem and the minimization of spectral partitions via the multiphase problem \eqref{multiphase}, the rhombic dodecahedron partition has similar properties as the honeycomb partition in 2D, which makes it a suitable candidate to be an optimizer in 3D. 

\begin{figure}
\centering
\setlength{\tabcolsep}{2pt}
\begin{tabular}{|c|c|c|c|c|c|}
\hline
\multicolumn{2}{|c|}{Weaire-Phelan} & & Kelvin & & \parbox{70pt}{Rhombic dodecahedron}\\ 
\hline
 $8.a$ $(\times 6)$ & $8.b$ $(\times 2)$  & $12$ & $16$ & $20$ & $32$ \\
 \hline
 \hline
 \includegraphics[height=0.13\textwidth]{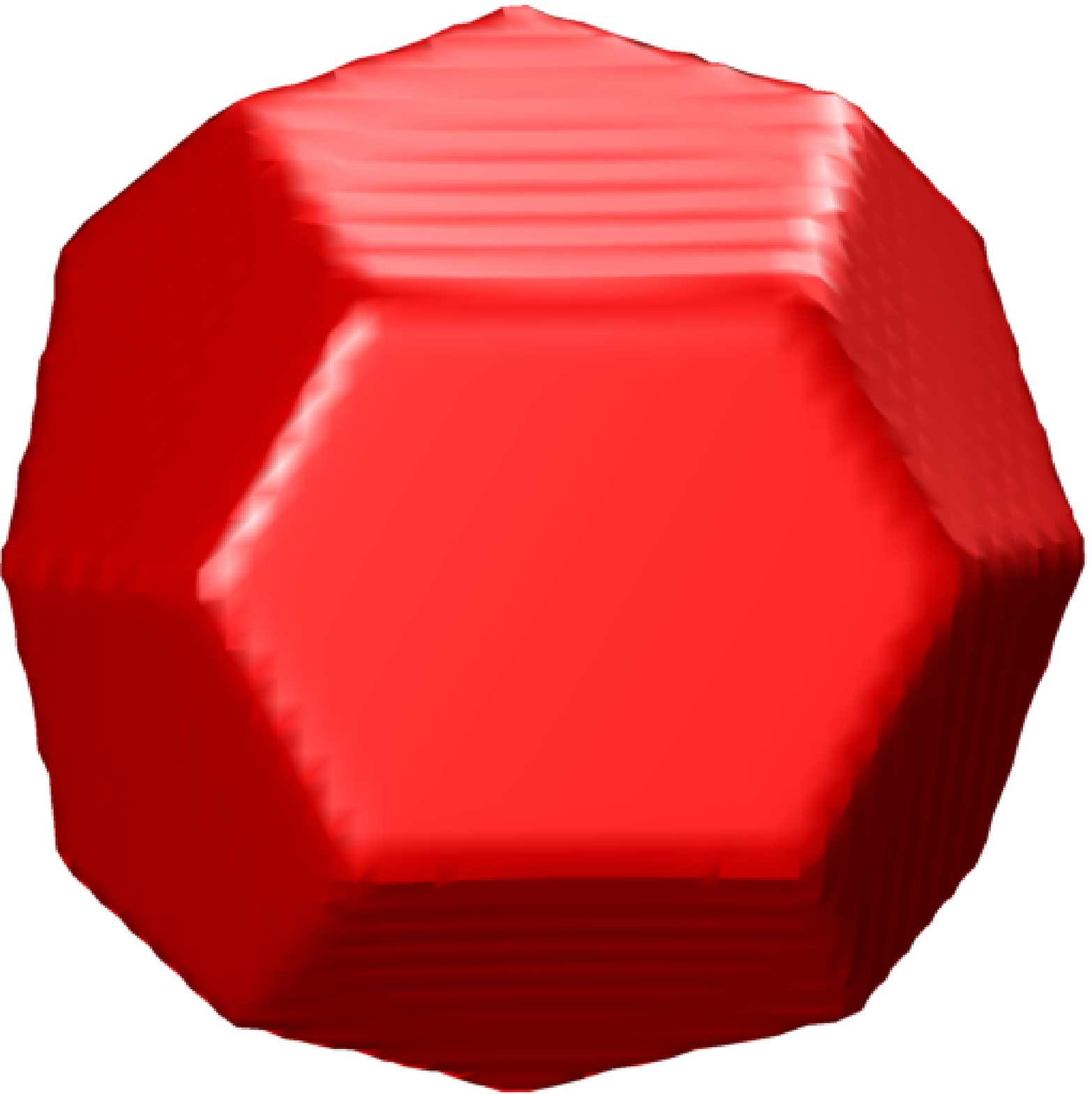} & \includegraphics[height=0.13\textwidth]{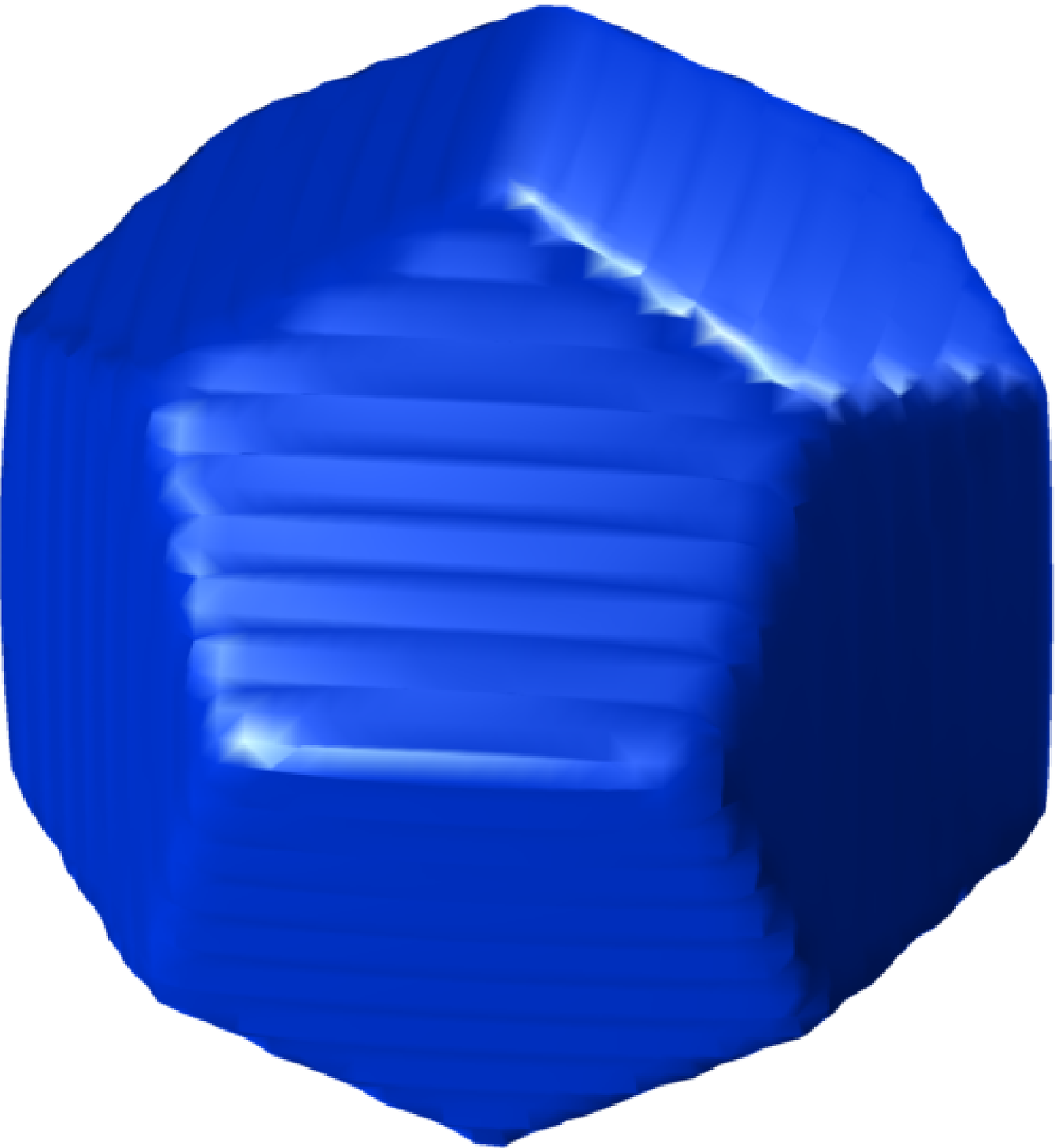} & \includegraphics[height=0.13\textwidth]{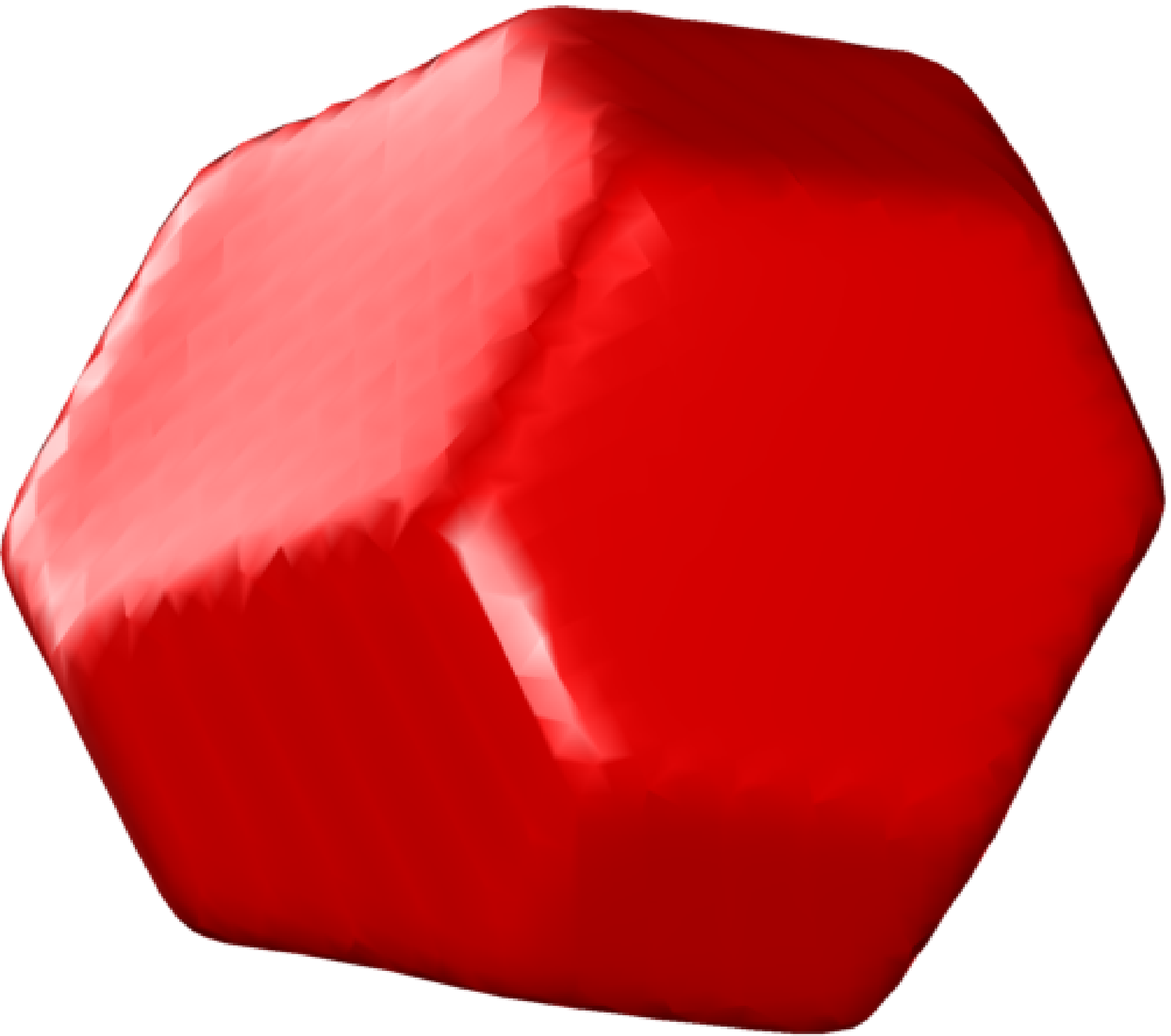} & \includegraphics[height=0.13\textwidth]{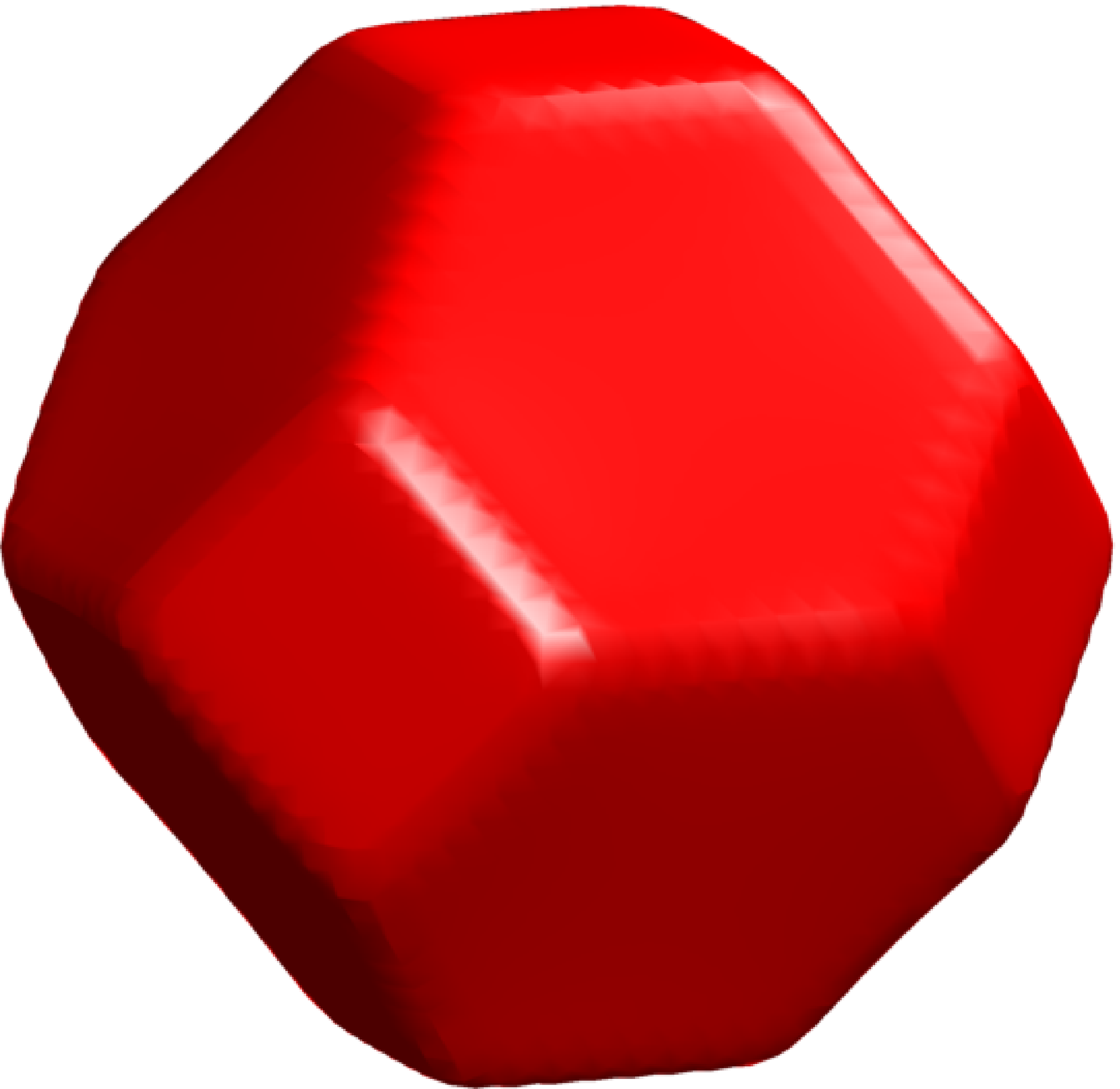} & \includegraphics[height=0.13\textwidth]{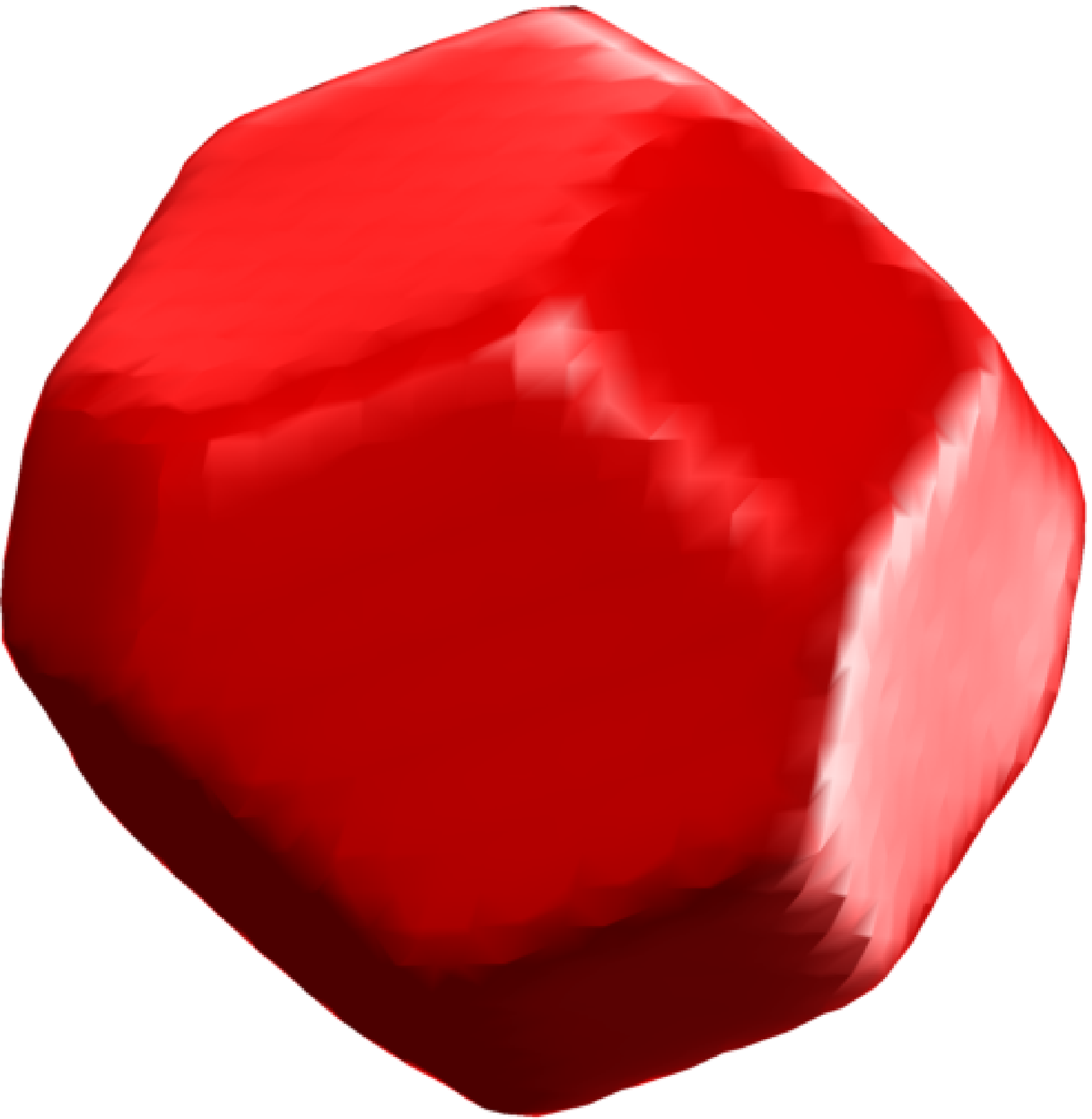} & \includegraphics[height=0.13\textwidth]{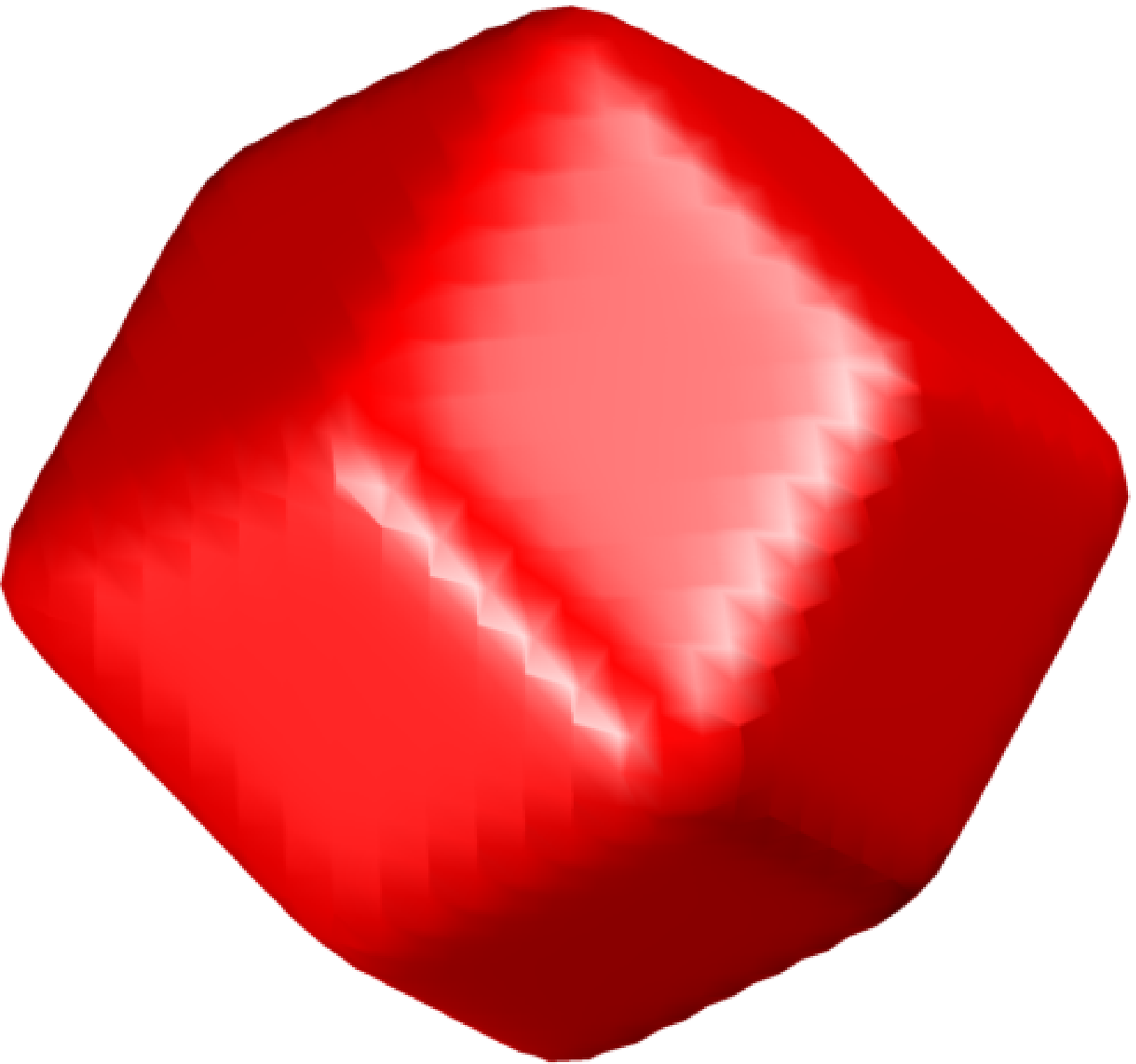} \\
 \hline
 \hline
 
$45045$p: $27.33$ & $41656$p: $26.97$  & $35905$p: $27.76$  & $126183$p: $27.11$ & $105967$p: $27.00$ & $114162$p: $26.94$ \\ 
$196122$p: $27.25$ & $219074$p: $26.88$ & $153788$p: $27.61$ & $198993$p: $27.13$ & $221404$p: $26.96$ & $296884$p: $26.91$ \\
\hline 
\end{tabular}
\caption{Eigenvalues of cells obtained in the periodic case}
\label{periodic-summary}
\end{figure}

We give a brief analysis of the cells obtained in the case of the periodic cube. Since we obtain different structures when considering different values of $n$, we should compare the value of the functional associated to each of these partition. This is difficult to do in the relaxed case, since the finite difference grid is not fine enough to give enough precision. Also, the error may depend on the discretization, the number of cells. Moreover, since we are looking for an asymptotic estimate on the number of cells, we should renormalize the eigenvalues so that the cells have the same volume. In order to achieve this, we do the following:
\begin{itemize}[leftmargin=10pt]
\item construct a triangulation of the surface of the cell $\omega$ using \texttt{isosurface} or Distmesh \cite{distmesh}.
\item build a tetrahedral mesh using iso2mesh \cite{iso2mesh}.
\item use $P1$ finite elements on the 3D mesh to compute the eigenvalue $\lambda_1(\omega)$ more precisely.
\item compute the scale invariant quantity $\lambda_1(\omega)\text{Vol}(\omega)^{2/3}$.
\end{itemize}
The shapes obtained in the periodic case as well as the values of their rescaled eigenvalues are presented in Figure \ref{periodic-summary}. We present two values for each cell, for different discretization parameters. In the table you can find the number of points and the corresponding eigenvalue. Even though the lowest eigenvalue was obtained for one of the cells in the Weaire-Phelan configuration, it is the one which repeats two times among the $8$ cell configuration. The lowest average eigenvalue is obtained for the rhombic dodecahedron configuration. Thus, heuristics presented in the previous paragraph and the computations presented in Figure \ref{periodic-summary} suggest that the rhombic dodecahedron partition might be a good optimal candidate in 3D. Note that all shapes presented Figure \ref{periodic-summary} were obtained directly from their density functions. This makes them have smoothed edges and vertices, therefore the values presented in Figure \ref{periodic-summary} may be a bit different than the ones corresponding to the exact shapes.

Since the conclusion of our analysis is the fact that the rhombic dodecahedron partition is the optimal numeric candidate, we computed the scale invariant eigenvalues for an exact rhombic dodecahedron for three meshes: $81216$ elements: $\lambda_1 = 27.0461$, $148433$ elements: $\lambda_1 = 27.0385$, $249523$ elements: $\lambda_1 = 27.0348$.

\subsection{Computation complexity}
\label{complexity}

We present in the following the parameters for some of the computations, which give an idea on the complexity of the problems the algorithm can solve. We also present the observed limits which need to be respected in order to have an efficient computation. The final parameters for some of the computations are presented in Table \ref{parameters}. Their signification is given below:
\begin{itemize}[leftmargin=*]
\item \emph{Geometry}: the domain for which we compute the partition
\item $n$: the number of cells in the partition
\item \emph{grid}: discretization parameters
\item \emph{degrees of freedom}: number of variables in the optimization problem. Can also be computed as $n \times (\text{number of points in the discretization grid})$. When we work on a non-rectangular domain $D$, only points inside $D$ are considered. 
\item \emph{average size of the eigenvalue problem}: the size of the matrix for which we compute the eigenvalue in the numerical algorithm. Also equal to the average number of points taken in the computational neighborhoods for each of the cells.
\end{itemize}

\begin{table}[!ht]
\centering
\begin{tabular}{|c|c|c|c|c|}
\hline
Geometry & $n$ & grid & deg. of freedom & average eig. pb. size \\
\hline
\hline
square  & $3$ & $480\times 480$ & $691\ 200$ & $81000$ \\
\hline
hex. union & $61$ & $400\times 400$ & $5\ 138\ 884$ & $2250$ \\
\hline
square & $1000$ & $1000\times 1000$ & $1\ 000 \ 000 \ 000$ & $2000$ \\
\hline 
cube & $3$ & $60 \times 60 \times 60$ & $648 \ 000$ & $70000$ \\
\hline 
cube & $4$ & $60 \times 60 \times 60$ & $864\ 000$  & $65500$ \\
\hline 

cube & $14$ & $80 \times 80 \times 80$  & $6\ 643\ 728$  & $50000$ \\
\hline
cube & $16$ & $70 \times 70 \times 70$  & $5\ 488\ 000$  & $35000$ \\
\hline 
tetrahedron & $35$ & $100 \times 100 \times 100$ & $11\ 648\ 280$  & $18100$ \\
\hline
tetrahedron & $84$ & $100 \times 100 \times 100$ & $27\ 955\ 872$  & $8500$ \\
\hline
tetrahedron & $120$ & $100\times 100\times 100$ & $39\ 936\ 960$ & $7500$ \\
\hline 
sphere & $3$ & $10242$ pts & $30\ 726$ & $4400$ \\
\hline 
sphere & $32$ & $163842$ pts & $5\ 242\ 944$ & $7500$ \\
\hline
sphere & $120$ & $163842$ pts & $19\ 661\ 040$ &  $2700$  \\
\hline

\end{tabular}

\caption{Parameters for various computations}
\label{parameters}
\end{table}

Let us give a few remarks regarding the information presented in Table \ref{parameters}. First of all, we see that the number of degrees of freedom can get really large, depending on the discretization and on the number of cells. For the case of $1000$ cells on the grid of size $1000\times 1000$ we have an optimization problem with one billion variables. What makes the optimization problem solvable numerically, even for so many variables, is the fact that for each grid point only few of the cells in the partition have a density function which is strictly positive there. Thus, even though all variables are free, the number of non-zero variables is much less than the number of degrees of freedom. 

Let us note a second aspect regarding the limits imposed on the grid discretization parameters and the average eigenvalue problem size. We observe that for the same discrete grid, having more cells decreases the size of the eigenvalue problem size. This is due to the fact that more cells need to be distributed on the same number of points, which makes the computational neighborhoods smaller. In practice, given a certain number of cells, one can consider finer and finer discretizations until the size of the eigenvalue problem reaches around $80000^2$ in 2D or $70000^2$ in 3D without having a too great computational cost for one \emph{cell}. This computational cost scales with the number of cells considered in the partition. Larger grid sizes may lead to a more important memory consumption and this should also be taken into consideration, as it was already discussed in the beginning of Section \ref{simulations}.

\section{Conclusions}

The algorithm described in previous sections shows significant improvement compared to \cite{buboou}. This allows great gains in computational cost for 2D and surface partitions and it allows us to work in 3D, for the first time on sufficiently large resolutions. The computations presented in 2D reinforce the Caffarelli-Lin conjecture that the honeycomb partition asymptotically minimizes the sum of the eigenvalues of the cells in the plane. 

Computations made on surfaces show that a similar behavior is present: optimal partitions contain patches of regular hexagons, while some other hapes may be necessary to fulfil Euler's relation regarding the number of faces, edges and vertices in the configuration.

Numerical results in 3D show interesting behaviour. We made computations on the cube, the ball and the regular tetrahedron, but the optimization algorithm does not pose any problem if we want to adapt it to other type of cells. Tests in the periodic case show different structures which may be candidates as optimizers in the asymptotic case in $\Bbb{R}^3$. Among these we have the rhombic dodecahedral structure which also appears in the interior of the partition of the regular tetrahedron into a pyramidal number of cells. Further analysis of the value of the eigenvalues of the cells obtained in the periodic case suggests that the partition into rhombic dodecahedra is likely to be optimal as asymptotic minimizer for the sum of the eigenvalues for a tiling of $\Bbb{R}^3$.  

\section*{Acknowledgements} 

This work was supported by the following grants of the French National Research Agency (ANR): LabEx Sciences Mathématiques de Paris ANR-10-LABX-0098 and ANR Optiform ANR-12-BS01-0007-02.

\bibliographystyle{alpha}
\bibliography{../master.bib}

\end{document}